\documentclass{aims}
\usepackage{amsmath}
  \usepackage{paralist}
  \usepackage{graphics} 
  \usepackage{epsfig} 
\usepackage{graphicx}  \usepackage{epstopdf}
 \usepackage[colorlinks=true]{hyperref}
\hypersetup{urlcolor=blue, citecolor=red}

\usepackage[title]{appendix}
\usepackage[list=true]{subcaption}

\def\currentvolume{X}
 \def\currentissue{X}
  \def\currentyear{200X}
   \def\currentmonth{XX}
    \def\ppages{X--XX}
     \def\DOI{10.3934/xx.xx.xx.xx}

\newtheorem{theorem}{Theorem}[section]

\newtheorem{lemma}[theorem]{Lemma}

\theoremstyle{definition}

\newtheorem{remark}{Remark}[section]

\title[On a porous-elastic system with damping and source terms] 
      {Existence, blow-up and exponential decay of solutions for a porous-elastic system with damping and source terms}

\author[Vo Anh Khoa, Le Thi Phuong Ngoc and Nguyen Thanh Long]{}


 \email{vakhoa.hcmus@gmail.com}
 \email{ngoc1966@gmail.com}
 \email{longnt2@gmail.com}

\thanks{This research is funded by Vietnam National University Ho Chi Minh City (VNU-HCM) under
	Grant no. \textbf{B2017-18-04}. The work of the first author was partly supported by  a postdoctoral fellowship of the Research Foundation-Flanders (FWO)}

\thanks{$^*$ Corresponding author: Nguyen Thanh Long}

\begin{document}
\maketitle

\centerline{\scshape Vo Anh Khoa}
\medskip
{\footnotesize
 \centerline{Institute for Numerical and Applied Mathematics}
   \centerline{University of Goettingen, Lotzestra\ss e 16-18, 37083 Goettingen, Germany}
   \centerline{Department of Mathematics and Computer Science}
   \centerline{VNUHCM-University of Science, 227 Nguyen Van Cu Str., Dist. 5, Ho Chi Minh City, Vietnam}
   \centerline{Faculty of Sciences, Hasselt University}
   \centerline{Campus Diepenbeek, Agoralaan
   	Building D, BE3590 Diepenbeek, Belgium}
} 

\medskip

\centerline{\scshape Le Thi Phuong Ngoc}
\medskip
{\footnotesize
	\centerline{University of Khanh Hoa, 01 Nguyen Chanh Str., Nha Trang City, Vietnam}
}

\medskip

\centerline{\scshape Nguyen Thanh Long$^*$}
\medskip
{\footnotesize
 \centerline{Department of Mathematics and Computer Science}
   \centerline{VNUHCM-University of Science, 227 Nguyen Van Cu Str., Dist. 5, Ho Chi Minh City, Vietnam}
}

\bigskip

 \centerline{(Communicated by the associate editor name)}

\begin{abstract}
In this paper we consider a porous-elastic system
consisting of nonlinear boundary/interior damping and 
nonlinear boundary/interior sources. Our interest lies in the
theoretical understanding of the existence, finite time blow-up of solutions
and their exponential decay using non-trivial adaptations of well-known techniques. First, we apply the conventional Faedo-Galerkin
method with standard arguments of density on the regularity of initial
conditions to establish two local existence theorems of weak solutions.
Moreover, we detail the uniqueness result in some specific cases.
In the second theme, we prove that any weak solution possessing negative initial energy
has the latent blow-up in finite time. Finally, we obtain the
so-called exponential decay estimates for the global solution under the
construction of a suitable Lyapunov functional. In order to corroborate our
theoretical decay, a numerical example is provided.
\end{abstract}

\section{Introduction}

This paper is concerned with the following polynomially damped system
of wave equations%
\begin{equation}
\begin{cases}
u_{tt}-u_{xx}+\lambda _{1}\left\vert u_{t}\right\vert
^{r_{1}-2}u_{t}=f_{1}\left( u,v\right) +F_{1}\left( x,t\right) , \\ 
v_{tt}-v_{xx}+\lambda _{2}\left\vert v_{t}\right\vert
^{r_{2}-2}v_{t}=f_{2}\left( u,v\right) +F_{2}\left( x,t\right) .%
\end{cases}
\tag{1.1}  \label{main1}
\end{equation}
This prototypical system of wave equations arises naturally within
frameworks of material science and quantum field theory. Accounting for the
Kirchhoff-Love plate theory in shear deformations, the system (\ref{main1})
is closely related to the Reissner-Mindlin plate equations (see \cite{8}),
  structured by three coupled wave and wave-like equations involving
the influence of nonlinear damping and source terms. Mathematically, systems
of wave equations have been extensively studied by many authors, see \cite%
{1,4,6,13,13a,14} and references therein where the existence, regularity and
the asymptotic behavior of solutions are investigated.

In \cite{6}, Guo et al. considered the local and global
well-posedness of a general system%
\begin{equation}
\begin{cases}
u_{tt}-u_{xx}+g_{1}\left( u_{t}\right) =f_{1}\left( u,v\right) , \\ 
v_{tt}-v_{xx}+g_{2}\left( v_{t}\right) =f_{2}\left( u,v\right) ,%
\end{cases}
\tag{1.2}  \label{eq:sample}
\end{equation}%
in a bounded domain of $\mathbb{R}^{n}$ with a nonlinear Robin boundary
condition on $u$ and a zero boundary conditions on $v$. The nonlinearities 
$f_{1}\left( u,v\right) $ and $f_{2}\left( u,v\right) $ are supercritical
exponents representing strong sources, while $g_{1}\left( u_{t}\right) $ and 
$g_{2}\left( v_{t}\right) \ $act as damping. These damping terms are
assumed to be continuous and monotone increasing functions vanishing at the
origin and satisfying restrictions in growing up at infinity. By employing the 
nonlinear semigroups and the theory of monotone operators, the
well-posedness of (\ref{eq:sample}) is moderately investigated. An important
result obtained in \cite{6} and further in \cite{13} is that every weak solution
blows up in finite time, provided the initial energy is negative and the
sources are more dominant than the damping involved in the system.

In \cite{3,4}, Cavalcanti et al. studied the existence of global
solutions, and showed the relation between the asymptotic behavior of the
energy and the degenerate system of wave equations with boundary conditions
of memory type. Constructing a suitable Lyapunov functional, the authors
proved that the energy decays exponentially. The same method was also used
in \cite{13a,14} to study the asymptotic behavior of the solutions to a
coupled system having integral convolutions as the memory terms. They proved
that the solution decays uniformly in time with rates depending on the speed
of decay of the convolutions kernel.

In recent years, various types of wave equations with linear or nonlinear
damping and sources have been solved by using the Galerkin approximation
(see, for instance, \cite{2,11,13,13a,15}). Based on \emph{a priori}
estimates, weak convergence, and compactness techniques, and via the
construction of a suitable Lyapunov functional, the existence, regularity,
blow-up, and exponential decay estimates of solutions for such typical wave
equations have been proved in \cite{2,11,15}. On the other side, the finite
time blow-up of any weak solutions with negative initial energy is obtained
in \cite{15}.

In light of the aforementioned works, we put ourselves into the study of the
existence, blow-up, and exponential decay estimate for the system (\ref%
{main1}).

Let $\Omega =\left( 0,1\right) $ and $Q_{T}=\Omega \times \left( 0,T\right) $
for $T>0$, a couple of real unknown functions $\left( u,v\right) $ is sought
for $\left( x,t\right) \in \overline{Q_{T}}=\left[ 0,1\right] \times \left[
0,T\right] $. The problem we consider here is made of (\ref{main1}) and the
following nonlinear boundary conditions%
\begin{equation}
\begin{cases}
u\left( 0,t\right) =0,\text{ }-u_{x}\left( 1,t\right) +K_{1}\left\vert
u\left( 1,t\right) \right\vert ^{p_{1}-2}u\left( 1,t\right) =\mu
_{1}\left\vert u_{t}\left( 1,t\right) \right\vert ^{q_{1}-2}u_{t}\left(
1,t\right) , \\ 
v_{x}\left( 0,t\right) +K_{2}\left\vert v\left( 0,t\right) \right\vert
^{p_{2}-2}v\left( 0,t\right) =\mu _{2}\left\vert v_{t}\left( 0,t\right)
\right\vert ^{q_{2}-2}v_{t}\left( 0,t\right) ,\text{ }v\left( 1,t\right) =0,%
\end{cases}
\tag{1.3}  \label{eq:main2}
\end{equation}%
and the initial conditions%
\begin{equation}
\begin{cases}
u\left( x,0\right) =\tilde{u}_{0}\left( x\right) ,\text{ }u_{t}\left(
x,0\right) =\tilde{u}_{1}\left( x\right) , \\ 
v\left( x,0\right) =\tilde{v}_{0}\left( x\right) ,\text{ }v_{t}\left(
x,0\right) =\tilde{v}_{1}\left( x\right) .%
\end{cases}
\tag{1.4}  \label{main3}
\end{equation}

As short-hand explanation for physical parameters in this system, the
constants $\lambda _{i}>0$ ($i=1,2$) are usually called as the friction terms,
while the constants $r_{i}\geq 2$ play a central role in deciding the order
of damped parts. The functions $f_{i}$ are known as the interior sources, while $%
F_{i}$ are the external functions. Moreover, we have on the boundary the
presence of the constants $K_{i}>0,$ $\mu _{i}>0,$ $p_{i}\geq 2,$ $q_{i}\geq
2$ ($i=1,2$) as well as given
functions $\tilde{u}_{i},$ $\tilde{v}_{i}$ ($i=0,1$) that will be specified later.

It is worth mentioning that the nonlinear boundary condition (\ref{eq:main2}%
) is the main difficulty we face, although the approach we use here is already 
analyzed in many simpler models. As far as we know, a nonlinear
wave equation with the two-point boundary conditions has been considered in 
\cite{11,15}. Nevertheless, the circumstance for the coupled system (\ref%
{main1})-(\ref{main3}) is still open and its rigorous treatment is technically demanding in this direction. In \cite{13}, the authors have solved the similar coupled system where the results are controlled not only by the damping orders $r_1$, $r_2$, but also by the involved parameters on the boundary. Note that compared again to \cite{13} with the homogeneous Dirichlet boundary condition for $u$, our paper needs a careful adaptation to handle several different parameters at the same time. On top of that, we would like to see how necessary assumptions on such input data will be established. Thus, now is the
moment we discover the answer.

\section{Preliminaries}
\subsection{Abstract settings}
Let us denote the usual functional spaces used in this paper by $%
C^{m}\left( \overline{\Omega }\right) ,$ $W^{m,p}=W^{m,p}\left( \Omega
\right) ,$ $L^{p}=W^{0,p}\left( \Omega \right) ,$ $H^{m}=W^{m,2}\left(
\Omega \right) $ for $1\leq p\leq \infty $ and $m\in \mathbb{N}$. Let $%
\langle \cdot ,\cdot \rangle $ be either the scalar product in $L^{2}$ or
the dual pairing of a continuous linear functional and an element of a
functional space. The notation $\left\Vert \cdot \right\Vert $ stands for
the norm in $L^{2}$ and we denote by $\left\Vert \cdot \right\Vert _{X}$ the
norm in the Banach space $X$. We call $X^{\prime }$ the dual space of $X$
and denote by $L^{p}\left( 0,T;X\right) ,$ $1\leq p\leq \infty $ for the
Banach space of the real functions $u:\left( 0,T\right) \rightarrow X$
measurable, such that%
\begin{equation*}
\left\Vert u\right\Vert _{L^{p}\left( 0,T;X\right) }=\left(
\int_{0}^{T}\left\Vert u\left( t\right) \right\Vert _{X}^{p}dt\right)
^{1/p}<\infty ,\quad \mbox{for}\;1\leq p<\infty ,
\end{equation*}%
and%
\begin{equation*}
\left\Vert u\right\Vert _{L^{\infty }\left( 0,T;X\right) }=\mbox{ess}%
\sup_{0<t<T}\left\Vert u\left( t\right) \right\Vert _{X},\quad \mbox{for}%
\;p=\infty .
\end{equation*}

Since the domain of interest is one-dimensional, let $u\left( t\right) ,\ u^{\prime } =u_{t} ,\
u^{\prime \prime }\left( t\right) =u_{tt}\left( t\right) ,\ \nabla
u=u_{x} ,\ \Delta u=u_{xx} \ $denote $%
u\left( x,t\right) ,$ $\frac{\partial u}{\partial t} ,\ 
\frac{\partial ^{2}u}{\partial t^{2}} ,$ $\frac{\partial u}{%
	\partial x} ,\ \frac{\partial ^{2}u}{\partial x^{2}},$ respectively.
In $H^{1}$, we use the following norm:%
\begin{equation*}
\left\Vert u\right\Vert _{H^{1}}=\left( \left\Vert u\right\Vert
^{2}+\left\Vert u_{x}\right\Vert ^{2}\right) ^{1/2}.
\end{equation*}
We define%
\begin{equation*}
\mathbb{V}_{1}=\left\{ v\in H^{1}:v\left( 0\right) =0\right\} ,\quad \mathbb{%
	V}_{2}=\left\{ v\in H^{1}:v\left( 1\right) =0\right\}
\end{equation*}%
two closed subspaces of $H^{1}$. Moreover, the following standard lemmas
read as the imbedding $H^{1}$ into $C^{0}\left( \overline{\Omega }\right) $,
and the equivalence between the norms $\left\Vert v_{x}\right\Vert $ and $%
\left\Vert v\right\Vert _{H^{1}}$ in $\mathbb{V}_{1}$ and $\mathbb{V}_{2}$.

\begin{lemma}
	The imbedding $H^{1}\hookrightarrow C\left( 
	\overline{\Omega }\right) $ is compact and the following inequality
	holds
	\begin{equation*}
	\left\Vert v\right\Vert _{C\left( \overline{\Omega }\right) }\leq \sqrt{2%
	}\left\Vert v\right\Vert _{H^{1}}\quad \mbox{for all}\;v\in H^{1}.
	\end{equation*}
\end{lemma}

\begin{lemma}\label{lem:2.2}
	On $\mathbb{V}_{1}$ \textit{and} $\mathbb{V}%
	_{2} $ two norms $v\mapsto \left\Vert v_{x}\right\Vert $ and $v\mapsto \left\Vert v\right\Vert _{H^{1}}$ are equivalent.
	Furthermore, it holds
	\begin{equation*}
	\left\Vert v\right\Vert _{C\left( \overline{\Omega }\right) }\leq
	\left\Vert v_{x}\right\Vert \quad \mbox{for all}\;v\in \mathbb{V}_{1}\;%
	\mbox{and}\;\mathbb{V}_{2}.
	\end{equation*}
\end{lemma}

For the sake of simplicity, we refer $\left( P\right) $ to the problem (\ref%
{main1}) endowed with the conditions (\ref{eq:main2})-(\ref{main3}). In
addition, we denote the damping terms and also possibly related functions by 
$\Psi _{r}\left( z\right) =\left\vert z\right\vert ^{r-2}z$ where $r$ is a
given constant.

\subsection{Weak formulation of $(P)$}
The weak formulation of the initial-boundary value problem $\left(P\right)$
can be given in the following manner:

For $T>0$, find a pair of real unknown solutions $\left( u,v\right) $ belonging to the
following functional space%
\begin{equation*}
\begin{tabular}{l}
$\mathbb{W}=\{\left( u,v\right) \in L^{\infty }\left( 0,T;\left( \mathbb{V}%
_{1}\cap H^{2}\right) \times \left( \mathbb{V}_{2}\cap H^{2}\right)
\right) :\left( u_{t},v_{t}\right) \in L^{\infty }\left( 0,T;\mathbb{V}%
_{1}\times \mathbb{V}_{2}\right) ,\smallskip $ \\ 
$\ \ \ \ \ \ \ \ \ \ \ \ \ \ \ \ \ \ \ \ \ \ \ \ \ \ \ \ \ \ \ \ \ \ \ \ \ \
\ \ \ \ \ \ \ \ \ \ \ \ \ \ \ \ \ \ \ \ \ \ \ \ \ \
\left( u_{tt},v_{tt}\right) \in L^{\infty }\left( 0,T;L^{2}\times
L^{2}\right) \},$%
\end{tabular}%
\end{equation*}%
such that $\left( u,v\right) $ satisfies the variational equations%
\begin{equation}
\begin{cases}
\langle u_{tt}\left( t\right) ,\phi \rangle +\langle u_{x}\left( t\right)
,\phi _{x}\rangle +\lambda _{1}\langle \Psi _{r_{1}}\left( u_{t}\left(
t\right) \right) ,\phi \rangle +\mu _{1}\Psi _{q_{1}}\left( u_{t}\left(
1,t\right) \right) \phi \left( 1\right) \smallskip \\ 
\text{ \ \ \ \ \ \ \ \ \ \ \ \ \ \ \ \ \ \ \ \ \ \ \ \ \ \ \ \ \ \ \ }%
=K_{1}\Psi _{p_{1}}\left( u\left( 1,t\right) \right) \phi \left( 1\right)
+\langle f_{1}\left( u,v\right) ,\phi \rangle +\langle F_{1}\left( t\right)
,\phi \rangle ,\smallskip \\ 
\langle v_{tt}\left( t\right) ,\tilde{\phi}\rangle +\langle v_{x}\left(
t\right) ,\tilde{\phi}_{x}\rangle +\lambda _{2}\langle \Psi _{r_{2}}\left(
v_{t}\left( t\right) \right) ,\tilde{\phi}\rangle +\mu _{2}\Psi
_{q_{2}}\left( v_{t}\left( 0,t\right) \right) \tilde{\phi}\left( 0\right)
\smallskip \\ 
\text{ \ \ \ \ \ \ \ \ \ \ \ \ \ \ \ \ \ \ \ \ \ \ \ \ \ \ \ \ \ \ \ }%
=K_{2}\Psi _{p_{2}}\left( v\left( 0,t\right) \right) \tilde{\phi}\left(
0\right) +\langle f_{2}\left( u,v\right) ,\tilde{\phi}\rangle +\langle
F_{2}\left( t\right) ,\tilde{\phi}\rangle ,%
\end{cases}
\tag{2.1}  \label{b1}
\end{equation}%
for all $(\phi ,\tilde{\phi})\in \mathbb{V}_{1}\times \mathbb{V}_{2}$ and for almost all $t\in (0,T)$. This system is endowed with the initial conditions%
\begin{equation}
\begin{tabular}{l}
$\left( u\left( 0\right) ,u_{t}\left( 0\right) \right) =\left( \tilde{u}_{0},%
\tilde{u}_{1}\right) ,\text{ \ }\left( v\left( 0\right) ,v_{t}\left(
0\right) \right) =\left( \tilde{v}_{0},\tilde{v}_{1}\right) .$%
\end{tabular}
\tag{2.2}  \label{b2}
\end{equation}

\section{The existence and uniqueness of a weak solution}
We now pose the following assumptions:

\textbf{(}$\mathbf{A}_{\mathbf{1}}$\textbf{)} $\left( \tilde{u}_{0},\tilde{u}%
_{1}\right) \in \left( \mathbb{V}_{1}\cap H^{2}\right) \times \mathbb{V}_{1}$
and $\left( \tilde{v}_{0},\tilde{v}_{1}\right) \in \left( \mathbb{V}_{2}\cap
H^{2}\right) \times \mathbb{V}_{2}$;

\textbf{(}$\mathbf{A}_{\mathbf{2}}$\textbf{)} $F_{1},F_{2}\in L^{2}\left(
Q_{T}\right) $ such that $F_{1}^{\prime },$ $F_{2}^{\prime }\in L^{1}\left(
0,T;L^{2}\right) $;

\textbf{(}$\mathbf{A}_{\mathbf{3}}$\textbf{)} there exists a  $C^{2}$-function $\mathcal{F}:%
\mathbb{R}^{2}\rightarrow \mathbb{R}$  such that%
\begin{equation}
\frac{\partial \mathcal{F}}{\partial u}\left( u,v\right) =f_{1}\left(
u,v\right) ,\text{ \ }\frac{\partial \mathcal{F}}{\partial v}\left(
u,v\right) =f_{2}\left( u,v\right) ,%
\tag{2.3}  \label{b3}
\end{equation}%
and there also exists the constants $\alpha ,$ $\beta >2$ and $C>0$ such that%
\begin{equation}
\begin{tabular}{l}
$\mathcal{F}\left( u,v\right) \leq C\left( 1+\left\vert u\right\vert
^{\alpha }+\left\vert v\right\vert ^{\beta }\right) ,\quad \mbox{for all}%
\;u,v\in \mathbb{R};$%
\end{tabular}
\tag{2.4}
\end{equation}

\begin{remark}\label{rem:3.1}
	There are several examples in which the functions $f_{1}
	$ and $f_{2}$ satisfy  \textbf{(}$\mathbf{A}_{\mathbf{3}}$\textbf{)%
	}, see e.g. \cite{1,13}. In particular, the authors in \cite{1}
	considered%
	\begin{equation*}
	\begin{tabular}{l}
	$\mathcal{F}\left( u,v\right) =\alpha \left\vert u+v\right\vert
	^{p+1}+2\beta \left\vert uv\right\vert ^{\frac{p+1}{2}},$%
	\end{tabular}%
	\end{equation*}%
	where $p\geq 3,$ $\alpha >1$ and $\beta >0$. In \cite{13}, the authors
	exploited another type%
	\begin{equation}
	\begin{tabular}{l}
	$\mathcal{F}\left( u,v\right) =\gamma _{1}\left( \left\vert u\right\vert
	^{\alpha }+\left\vert v\right\vert ^{\beta }\right) +\gamma _{2}\left\vert
	u\right\vert ^{\frac{\alpha }{2}}\left\vert v\right\vert ^{\frac{\beta }{2}%
	}, $%
	\end{tabular}
	\tag{2.5}  \label{b5}
	\end{equation}%
	where $\alpha ,$ $\beta ,$ $\gamma _{1}$ and $\gamma _{2}$ are positive
	constants with $\gamma _{2}<2\gamma _{1}$.
\end{remark}
In the following, we claim the existence and uniqueness of a weak solution. Depending on data assumptions, we consider the weak setting (\ref{b1})-(\ref{b2}) to show that the problem $(P)$ has a solution in some given function spaces with a small time length $T$. In some subcases of $q_1,q_2$ and $p_1,p_2$, we obtain the uniqueness result in such a small time. 

\begin{theorem}\label{thm:3.1}
	Suppose that \emph{\textbf{(}$\mathbf{A}_{\mathbf{1}%
	}$\textbf{)-(}$\mathbf{A}_{\mathbf{3}}$\textbf{)}} hold and the
	initial data obey the compatibility relation
	\begin{equation}
	\begin{cases}
	-\tilde{u}_{0x}\left( 1\right) +K_{1}\Psi _{p_{1}}\left( \tilde{u}_{0}\left(
	1\right) \right) =\mu _{1}\Psi _{q_{1}}\left( \tilde{u}_{1}\left( 1\right)
	\right) ,\smallskip \\ 
	\tilde{v}_{0x}\left( 0\right) +K_{2}\Psi _{p_{2}}\left( \tilde{v}_{0}\left(
	0\right) \right) =\mu _{2}\Psi _{q_{2}}\left( \tilde{v}_{1}\left( 0\right)
	\right) .%
	\end{cases}
	\tag{3.1}  \label{c1}
	\end{equation}
	If $ p_{1},p_{2},q_{1},q_{2} $ are such that
	\begin{equation*}
	\begin{cases}
	p_{1},\text{ }p_{2}\geq 2, \\ 
	2\leq q_{1},\text{ }q_{2}\leq 4,%
	\end{cases}%
	\mbox{or }\;%
	\begin{cases}
	p_{1},\text{ }p_{2}\in \left\{ 2\right\} \cup \left[ 3,\infty \right) , \\ 
	q_{1},\text{ }q_{2}>4,%
	\end{cases}%
	\end{equation*}%
	and $r_{1},$ $r_{2}\geq 2,$ then there exists a local weak
	solution $\left( u,v\right) $ of the problem $\left( P\right) $ such that
	\begin{equation}
	\begin{cases}
	\left( u,v\right) \in L^{\infty }\left( 0,T_{\ast };\left( \mathbb{V}%
	_{1}\cap H^{2}\right) \times \left( \mathbb{V}_{2}\cap H^{2}\right) \right) ,
	\\ 
	\left( u_{t},v_{t}\right) \in L^{\infty }\left( 0,T_{\ast };\mathbb{V}%
	_{1}\times \mathbb{V}_{2}\right) , \\ 
	\left( u_{tt},v_{tt}\right) \in L^{\infty }\left( 0,T_{\ast};L^{2}\times
	L^{2}\right) , \\ 
	\left\vert u_{t}\right\vert ^{\frac{r_{1}}{2}-1}u_{t},\text{ }\left\vert
	v_{t}\right\vert ^{\frac{r_{2}}{2}-1}v_{t}\in H^{1}\left( Q_{T_{\ast
	}}\right) , \\ 
	\left\vert u_{t}\left( 1,\cdot \right) \right\vert ^{\frac{q_{1}}{2}%
		-1}u_{t}\left( 1,\cdot \right) ,\text{ }\left\vert v_{t}\left( 0,\cdot
	\right) \right\vert ^{\frac{q_{2}}{2}-1}v_{t}\left( 0,\cdot \right) \in
	H^{1}\left( 0,T_{\ast }\right) ,%
	\end{cases}
	\tag{3.2}  \label{c2}
	\end{equation}%
	for $T_{\ast }>0$ sufficiently small. Furthermore, if $q_{1}=q_{2}=2$ and $p_{1},p_{2}\geq 2$, the obtained
	solution is unique. 
\end{theorem}

\begin{proof}[Proof of Theorem \ref{thm:3.1}]
	$\;$\\
	\textbf{Step 1. The Faedo-Galerkin approximation.} Let $\{(\phi _{i},\tilde{\phi}_{j})\}$ be a denumerable base of $\left( 
	\mathbb{V}_{1}\cap H^{2}\right) \times \left( \mathbb{V}_{2}\cap
	H^{2}\right) $. The approximate solution of $\left( P\right) $ is a sequence 
	$\left\{ \left( u_{m},v_{m}\right) \right\}_{m\in \mathbb{N}} $ structured as%
	\begin{equation*}
	u_{m}\left( t\right) =\sum_{j=1}^{m}c_{mj}\left( t\right) \phi
	_{j},\quad v_{m}\left( t\right) =\sum_{j=1}^{m}d_{mj}\left(
	t\right) \tilde{\phi}_{j},%
	\end{equation*}%
	where the time-dependent coefficient functions $\left( c_{mj},d_{mj}\right) $
	satisfy the following system%
	\begin{equation}
	\begin{cases}
	\langle u_{m}^{\prime \prime }\left( t\right) ,\phi _{j}\rangle +\langle
	u_{mx}\left( t\right) ,\phi _{jx}\rangle +\lambda _{1}\langle \Psi
	_{r_{1}}\left( u_{m}^{\prime }\left( t\right) \right) ,\phi _{j}\rangle +\mu
	_{1}\Psi _{q_{1}}\left( u_{m}^{\prime }\left( 1,t\right) \right) \phi
	_{j}\left( 1\right) \smallskip \\ 
	\text{ \ \ \ \ \ \ \ \ \ \ \ \ \ \ }=K_{1}\Psi _{p_{1}}\left( u_{m}\left(
	1,t\right) \right) \phi _{j}\left( 1\right) +\langle f_{1}\left(
	u_{m},v_{m}\right) ,\phi _{j}\rangle +\langle F_{1}\left( t\right) ,\phi
	_{j}\rangle ,\smallskip \\ 
	\langle v_{m}^{\prime \prime }\left( t\right) ,\tilde{\phi}_{j}\rangle
	+\langle v_{mx}\left( t\right) ,\tilde{\phi}_{jx}\rangle +\lambda
	_{2}\langle \Psi _{r_{2}}\left( v_{m}^{\prime }\left( t\right) \right) ,%
	\tilde{\phi}_{j}\rangle +\mu _{2}\Psi _{q_{2}}\left( v_{m}^{\prime }\left(
	0,t\right) \right) \tilde{\phi}_{j}\left( 0\right) \smallskip \\ 
	\text{ \ \ \ \ \ \ \ \ \ \ \ \ \ \ }=K_{2}\Psi _{p_{2}}\left( v_{m}\left(
	0,t\right) \right) \tilde{\phi}_{j}\left( 0\right) +\langle f_{2}\left(
	u_{m},v_{m}\right) ,\tilde{\phi}_{j}\rangle +\langle F_{2}\left( t\right) ,%
	\tilde{\phi}_{j}\rangle ,\smallskip \\ 
	\left( u_{m}\left( 0\right) ,u_{m}^{\prime }\left( 0\right) \right) =\left( 
	\tilde{u}_{0},\tilde{u}_{1}\right) ,\text{ }\left( v_{m}\left( 0\right)
	,v_{m}^{\prime }\left( 0\right) \right) =\left( \tilde{v}_{0},\tilde{v}%
	_{1}\right) ,%
	\end{cases}
	\tag{3.3}  \label{c3}
	\end{equation}%
	for $1\leq j\leq m$. A combination of assumptions of this theorem is the
	direct argument to gain the existence of solution $\left( u_{m},v_{m}\right) 
	$ for the system (\ref{c3}) on an interval $\left[ 0,T_{m}\right] \subset %
	\left[ 0,T\right] $.
	
	$\;$\\
	\textbf{Step 2. A priori estimates.}\\
	\textit{Step 2.1. The first estimate.} Multiplying the $j$-th system of (\ref%
	{c3}) (specifically, we multiply the first equation by $c_{mj}^{\prime
	}\left( t\right) $ and the second equation by $d_{mj}^{\prime }\left(
	t\right) $), summing up to $m$ with respect to $j$, and then integrating the
	resulting equation with respect to the time variable from 0 to $t$, we
	obtain the following equation%
	\begin{align}
	 \label{c4} \tag{3.4}
	   & \mathcal{S}_{m}\left( t\right)  \\  &=\mathcal{S}_{m}\left( 0\right)
	+K_{1}\int_{0}^{t}\Psi _{p_{1}}\left( u_{m}\left( 1,s\right) \right)
	u_{m}^{\prime }\left( 1,s\right) ds+K_{2}\int_{0}^{t}\Psi _{p_{2}}\left(
	v_{m}\left( 0,s\right) \right) v_{m}^{\prime }\left( 0,s\right) ds 
   \nonumber \\
	&+2\int_{0}^{t}\left[ \langle \frac{\partial \mathcal{F}}{\partial u}\left(
	u_{m}\left( s\right) ,v_{m}\left( s\right) \right) ,u_{m}^{\prime }\left(
	s\right) \rangle +\langle \frac{\partial \mathcal{F}}{\partial v}\left(
	u_{m}\left( s\right) ,v_{m}\left( s\right) \right) ,v_{m}^{\prime }\left(
	s\right) \rangle \right] ds  \nonumber \\
	&+2\int_{0}^{t}\left[ \langle F_{1}\left( s\right) ,u_{m}^{\prime }\left(
	s\right) \rangle +\langle F_{2}\left( s\right) ,v_{m}^{\prime }\left(
	s\right) \rangle \right] ds  \nonumber \\
	&= \mathcal{S}_{m}\left( 0\right) +\mathcal{I}_{1}+\mathcal{I}_{2}+\mathcal{I%
	}_{3}+\mathcal{I}_{4},  \nonumber
	\end{align}
	where we have denoted by%
	\begin{align}
	&  \mathcal{S}_{m}\left( t\right) \label{c5} \tag{3.5}\\ &:=\left\Vert u_{m}^{\prime }\left( t\right)
	\right\Vert ^{2}+\left\Vert v_{m}^{\prime }\left( t\right) \right\Vert
	^{2}+\left\Vert u_{mx}\left( t\right) \right\Vert ^{2}+\left\Vert
	v_{mx}\left( t\right) \right\Vert ^{2}+2\lambda _{1}\int_{0}^{t}\left\Vert
	u_{m}^{\prime }\left( s\right) \right\Vert _{L^{r_{1}}}^{r_{1}}ds 
	 \nonumber   \\
	&+2\lambda _{2}\int_{0}^{t}\left\Vert v_{m}^{\prime }\left( s\right)
	\right\Vert _{L^{r_{2}}}^{r_{2}}ds+2\mu _{1}\int_{0}^{t}\left\vert
	u_{m}^{\prime }\left( 1,s\right) \right\vert ^{q_{1}}ds+2\mu
	_{2}\int_{0}^{t}\left\vert v_{m}^{\prime }\left( 0,s\right) \right\vert
	^{q_{2}}ds.  \notag
	\end{align}
	By the structures of (\ref{c5}) and the third equation in (\ref{c3}), there
	exists a positive constant $\mathcal{S}_{0}$ such that that for all $m\in 
	\mathbb{N}$%
	\begin{equation*}
	\begin{tabular}{l}
	$\mathcal{S}_{m}\left( 0\right) = \left\Vert \tilde{u}_{1}\right\Vert
	^{2}+\left\Vert \tilde{v}_{1}\right\Vert ^{2}+\left\Vert \tilde{u}%
	_{0x}\right\Vert ^{2}+\left\Vert \tilde{v}_{0x}\right\Vert ^{2}\equiv 
	\mathcal{S}_{0}.$%
	\end{tabular}%
	\end{equation*}
	
	From now on, we estimate from above the integrals $\mathcal{I}_{k}$ for $k=%
	\overline{1,4}$ in (\ref{c4}). To do this, we need the following elementary
	inequalities.
	\begin{remark}[Young-type inequality]
		Let $\delta >0$ and $%
		a,b\geq 0$ be arbitrarily real numbers and given $q,q^{\prime }>1$ real
		constants which are H\"{o}lder conjugates of each other. The following
		inequality holds%
		\begin{equation}
		ab\leq \frac{1}{q}\delta ^{q}a^{q}+\frac{1}{q^{\prime }}\delta ^{-q^{\prime
		}}b^{q^{\prime }}.
		\tag{3.6}  \label{c6}
		\end{equation}
	\end{remark}
	
	\begin{remark}\label{rem:3.2}
		Given $N=\frac{1}{2}\max \left\{2;\alpha ;\beta ;\frac{%
			q_{1}\left( p_{1}-1\right) }{q_{1}-1};\frac{q_{2}\left( p_{2}-1\right) }{%
			q_{2}-1}\right\}$, then for all $s\geq 0$, the inequality $s^{\gamma }\leq 1+s^{N}$
		holds for all $\gamma \in \left( 0,N\right] $.
	\end{remark}

	Observe that the estimates for $\mathcal{I}_{1}$ and $\mathcal{I}_{2}$
	are similar when using the above inequalities. By choosing $\delta >0$
	in such a way that%
	\begin{equation*}
	\delta =\min \left\{ \sqrt[q_{1}]{\frac{\mu _{1}q_{1}}{2K_{1}}};\sqrt[q_{2}]%
	{\frac{\mu _{2}q_{2}}{2K_{2}}}\right\} ,
	\end{equation*}%
	it naturally arises that for $C_{T}>0$ only dependent of $T$, it holds%
	\begin{equation}
	\mathcal{I}_{1}+\mathcal{I}_{2}\leq \frac{1}{2}\mathcal{S}_{m}\left(
	t\right) +C_{T}\int_{0}^{t}\left( 1+\mathcal{S}_{m}^{N}\left( s\right)
	\right) ds.  \tag{3.7}  \label{c7}
	\end{equation}
	
	Additionally, for $C_{0}>0$ depending only on the initial data $\tilde{u}%
	_{0},$ $\tilde{v}_{0},$ $\tilde{u}_{1},$ $\tilde{v}_{1}$, and the constants $%
	\alpha ,$ $\beta $, we can prove that%
	\begin{equation*}
	\left\Vert u_{m}\left( t\right) \right\Vert _{L^{\alpha }}^{\alpha
	}+\left\Vert v_{m}(t)\right\Vert _{L^{\beta }}^{\beta }\leq C_{0}+\left(
	\alpha +\beta \right) \int_{0}^{t}\left( 1+\mathcal{S}_{m}^{N}\left(
	s\right) \right) ds,
	\end{equation*}%
	and then in the same spirit of \cite{13} the integral $\mathcal{I}_{3}$ is essentially estimated by%
	\begin{align}
	\mathcal{I}_{3} &\leq 2\sup_{\left\vert y\right\vert ,\left\vert
		z\right\vert \leq \sqrt{C_{0}}}\left\vert \mathcal{F}\left( y,z\right)
	\right\vert +2C_{1}+2C_{1}\left[ C_{0}+\left( \alpha +\beta \right)
	\int_{0}^{t}\left( 1+\mathcal{S}_{m}^{N}\left( s\right) \right) ds\right] 
	\tag{3.8}  \label{c8} \\
	&\leq C_{0}+C_{0}\int_{0}^{t}\left( 1+\mathcal{S}_{m}^{N}\left( s\right)
	\right) ds,  \notag  \label{eq:esI3}
	\end{align}%
	where we have also used the assumption \textbf{(}$\mathbf{A}_{\mathbf{3}}$%
	\textbf{)}. The last integral can be bounded by the standard
	Cauchy-Schwartz inequality, i.e.%
	\begin{equation}
	\mathcal{I}_{4}\leq \left\Vert F_{1}\right\Vert _{L^{2}\left( Q_{T}\right)
	}^{2}+\left\Vert F_{2}\right\Vert _{L^{2}\left( Q_{T}\right)
	}^{2}+\int_{0}^{t}\mathcal{S}_{m}\left( s\right) ds\leq
	C_{T}+\int_{0}^{t}\left( 1+\mathcal{S}_{m}^{N}\left( s\right) \right) ds. 
	\tag{3.9}  \label{c9}
	\end{equation}
	
	We hence obtain from (\ref{c7})-(\ref{c9}) that for $0\leq t\leq T_{m}$%
	\begin{equation*}
	\mathcal{S}_{m}\left( t\right) \leq C_{T}+C_{T}\int_{0}^{t}\left( 1+\mathcal{%
		S}_{m}^{N}\left( s\right) \right) ds.
	\end{equation*}
	
	Cf. \cite{9}, it should be noted that there exists a constant $T_{\ast }>0$
	depending on $T$ (but independent of $m$) such that%
	\begin{equation}
	\mathcal{S}_{m}\left( t\right) \leq C_{T},\text{ }\forall m\in \mathbb{N},%
	\text{ }\forall t\in \left[ 0,T_{\ast }\right] .  \tag{3.10}  \label{c10}
	\end{equation}
	Consequently, this result allows us to take $T_{m}=T_{\ast }$ for all $m$.\\
	\textit{Step 2.2. The second estimate.}  Consider the first equation
	of (\ref{c3}). Letting $t\rightarrow 0^{+}$, then multiplying the equation
	by $c_{mj}^{\prime \prime }\left( 0\right) $ with summing up to $m$ with
	respect to $j$, and using the first compatibility relation (\ref{c1}), we
	thus obtain%
	\begin{align*}
	 \left\Vert u_{m}^{\prime \prime }\left( 0\right) \right\Vert ^{2}&-\langle
	\Delta u_{m}\left( 0\right) ,u_{m}^{\prime \prime }\left( 0\right) \rangle
	+\lambda _{1}\langle \left\vert \tilde{u}_{1}\right\vert ^{r_{1}-2}\tilde{u}%
	_{1},u_{m}^{\prime \prime }\left( 0\right) \rangle \\ &  =\langle f_{1}\left( 
	\tilde{u}_{0},\tilde{v}_{0}\right) ,u_{m}^{\prime \prime }\left( 0\right)
	\rangle +\langle F_{1}\left( 0\right) ,u_{m}^{\prime \prime }\left( 0\right)
	\rangle .
	\end{align*}
	
	Relying on the classical inequalities, we have%
	\begin{equation}
	\left\Vert u_{m}^{\prime \prime }\left( 0\right) \right\Vert \leq \left\Vert
	\Delta \tilde{u}_{0}\right\Vert +\lambda _{1}\left\Vert \left\vert \tilde{u}%
	_{1}\right\vert ^{r_{1}-1}\right\Vert +\left\Vert f_{1}\left( \tilde{u}_{0},%
	\tilde{v}_{0}\right) \right\Vert +\left\Vert F_{1}\left( 0\right)
	\right\Vert ,  \tag{3.11}  \label{c11}
	\end{equation}%
	then state that there exists $\mathcal{C}_{1}>0$ such that 
	$\left\Vert u_{m}^{\prime \prime }\left( 0\right) \right\Vert \leq \mathcal{C%
	}_{1}$ for all $m\in \mathbb{N}$.
	
	For the second equation of (\ref{c3}), one also proves without difficulty
	using similar arguments that there exists $\mathcal{C}_{2}>0$ in which it
	bounds $\left\Vert v_{m}^{\prime \prime }\left( 0\right) \right\Vert $, i.e. 
	$\left\Vert v_{m}^{\prime \prime }\left( 0\right) \right\Vert \leq \mathcal{C%
	}_{2}$ for all $m\in \mathbb{N}$.
	
	Next, we differentiate (\ref{c3}) with respect to $t$. This way the first equation
	becomes%
	\begin{align}\tag{3.12}  \label{c12}
	\langle u_{m}^{\prime \prime \prime }\left( t\right) ,\phi _{j}\rangle&  +\langle u_{mx}^{\prime }\left( t\right) ,\phi _{jx}\rangle\\
	&  +\lambda
	_{1}\langle \Psi _{r_{1}}^{\prime }\left( u_{m}^{\prime }\left( t\right)
	\right) u_{m}^{\prime \prime }\left( t\right) ,\phi _{j}\rangle +\mu
	_{1}\Psi _{q_{1}}^{\prime }\left( u_{m}^{\prime }\left( 1,t\right) \right)
	u_{m}^{\prime \prime }\left( 1,t\right) \phi _{j}\left( 1\right) 	 \nonumber \\ 
	&=K_{1}\Psi _{p_{1}}^{\prime }\left( u_{m}\left( 1,t\right) \right)
	u_{m}^{\prime }\left( 1,t\right) \phi _{j}\left( 1\right) \nonumber\\ &+\langle \frac{%
		\partial ^{2}\mathcal{F}}{\partial u^{2}}\left( u_{m},v_{m}\right)
	u_{m}^{\prime }+\frac{\partial ^{2}\mathcal{F}}{\partial u\partial v}\left(
	u_{m},v_{m}\right) v_{m}^{\prime },\phi _{j}\rangle +\langle F_{1}^{\prime
	}\left( t\right) ,\phi _{j}\rangle , \nonumber %
	\end{align}
	and for $1\leq j\leq m$, the second equation is%
	\begin{align}	\tag{3.13}  \label{c13}
	\langle v_{m}^{\prime \prime \prime }\left( t\right) ,\tilde{\phi}%
	_{j}\rangle & +\langle v_{mx}^{\prime }\left( t\right) ,\tilde{\phi}%
	_{jx}\rangle \nonumber \\&+\lambda _{2}\langle \Psi _{r_{2}}^{\prime }\left(
	v_{m}^{\prime }\left( t\right) \right) v_{m}^{\prime \prime }\left( t\right)
	,\tilde{\phi}_{j}\rangle +\mu _{2}\Psi _{q_{2}}^{\prime }\left(
	v_{m}^{\prime }\left( 0,t\right) \right) v_{m}^{\prime \prime }\left(
	0,t\right) \tilde{\phi}_{j}\left( 0\right) \nonumber \\ 
	& =K_{2}\Psi _{p_{2}}^{\prime }\left( v_{m}\left( 0,t\right) \right)
	v_{m}^{\prime }\left( 0,t\right) \tilde{\phi}_{j}\left( 0\right) \nonumber\\ &+\langle 
	\frac{\partial ^{2}\mathcal{F}}{\partial v\partial u}\left(
	u_{m},v_{m}\right) u_{m}^{\prime }+\frac{\partial ^{2}\mathcal{F}}{\partial
		v^{2}}\left( u_{m},v_{m}\right) v_{m}^{\prime },\tilde{\phi}_{j}\rangle
	+\langle F_{2}^{\prime }\left( t\right) ,\tilde{\phi}_{j}\rangle .\nonumber %
	\end{align}%
	
	Multiplying the $j $-th equation of (\ref{c12}) and (\ref{c13}),
	respectively, by $c_{mj}^{\prime \prime }\left( t\right) $ and $%
	d_{mj}^{\prime \prime }\left( t\right) $, summing with respect to $j$ up to $%
	m$, and then integrating with respect to the time variable from $0$ to $t$,
	we obtain%
	\begin{align}\tag{3.14}  \label{c14}
	& \mathcal{P}_{m}\left( t\right) \\&=\mathcal{P}_{m}\left( 0\right)
	+2\int_{0}^{t}\langle \frac{\partial ^{2}\mathcal{F}}{\partial u^{2}}\left(
	u_{m},v_{m}\right) u_{m}^{\prime }\left( s\right) +\frac{\partial ^{2}%
		\mathcal{F}}{\partial u\partial v}\left( u_{m},v_{m}\right) v_{m}^{\prime
	}\left( s\right) ,u_{m}^{\prime \prime }\left( s\right) \rangle ds  \nonumber \\
	&+2\int_{0}^{t}\langle \frac{\partial ^{2}\mathcal{F}}{\partial v\partial u}%
	\left( u_{m},v_{m}\right) u_{m}^{\prime }\left( s\right) +\frac{\partial ^{2}%
		\mathcal{F}}{\partial v^{2}}\left( u_{m},v_{m}\right) v_{m}^{\prime }\left(
	s\right) ,v_{m}^{\prime \prime }\left( s\right) \rangle ds  \notag \\
	&+2\int_{0}^{t}\left[ \langle F_{1}^{\prime }\left( s\right) ,u_{m}^{\prime
		\prime }\left( s\right) \rangle +\langle F_{2}^{\prime }\left( s\right)
	,v_{m}^{\prime \prime }\left( s\right) \rangle \right] ds  \notag \\
	&+2\int_{0}^{t}\left[ K_{1}\Psi _{p_{1}}^{\prime }\left( u_{m}\left(
	1,s\right) \right) u_{m}^{\prime }\left( 1,s\right) u_{m}^{\prime \prime
	}\left( 1,s\right) +K_{2}\Psi _{p_{2}}^{\prime }\left( v_{m}\left(
	0,s\right) \right) v_{m}^{\prime }\left( 0,s\right) v_{m}^{\prime \prime
	}\left( 0,s\right) \right] ds  \notag \\
	&=\mathcal{P}_{m}\left( 0\right) +\mathcal{J}_{1}+\mathcal{J}_{2}+\mathcal{J%
	}_{3}+\mathcal{J}_{4},  \notag 
	\end{align}%
	where we have denoted by%
	\begin{align}\tag{3.15}  \label{c15}
	\mathcal{P}_{m}\left( t\right) &:=\left\Vert u_{m}^{\prime \prime }\left(
	t\right) \right\Vert ^{2}+\left\Vert v_{m}^{\prime \prime }\left( t\right)
	\right\Vert ^{2}+\left\Vert u_{mx}^{\prime }\left( t\right) \right\Vert
	^{2}+\left\Vert v_{mx}^{\prime }\left( t\right) \right\Vert ^{2} \\
	&+\frac{8\lambda _{1}\left( r_{1}-1\right) }{r_{1}^{2}}\int_{0}^{t}\left%
	\Vert \frac{\partial }{\partial s}\left( \left\vert u_{m}^{\prime }\left(
	s\right) \right\vert ^{\frac{r_{1}}{2}-1}u_{m}^{\prime }\left(s\right)
	\right) \right\Vert ^{2}ds  \notag \\
	&+\frac{8\lambda _{2}\left( r_{2}-1\right) }{r_{2}^{2}}\int_{0}^{t}\left%
	\Vert \frac{\partial }{\partial s}\left( \left\vert v_{m}^{\prime }\left(
	s\right) \right\vert ^{\frac{r_{2}}{2}-1}v_{m}^{\prime }\left(s\right)
	\right) \right\Vert ^{2}ds  \notag \\
	&+\frac{8\mu _{1}\left( q_{1}-1\right) }{q_{1}^{2}}\int_{0}^{t}\left\vert 
	\frac{\partial }{\partial s}\left( \left\vert u_{m}^{\prime }\left(
	1,s\right) \right\vert ^{\frac{q_{1}}{2}-1}u_{m}^{\prime }\left( 1,s\right)
	\right) \right\vert ^{2}ds  \notag \\
	&+\frac{8\mu _{2}\left( q_{2}-1\right) }{q_{2}^{2}}\int_{0}^{t}\left\vert 
	\frac{\partial }{\partial s}\left( \left\vert v_{m}^{\prime }\left(
	0,s\right) \right\vert ^{\frac{q_{2}}{2}-1}v_{m}^{\prime }\left( 0,s\right)
	\right) \right\vert ^{2}ds.  \notag 
	\end{align}
	
	Combining the arguments from the boundedness of $\left\Vert u_{m}^{\prime
		\prime }\left( 0\right) \right\Vert ,$ $\left\Vert v_{m}^{\prime \prime
	}\left( 0\right) \right\Vert $ and the third equation of (\ref{c3}) to (\ref%
	{c15}), there exists a positive constant $\mathcal{P}_{0}$ that bounds $%
	\mathcal{P}_{m}\left( 0\right) $ for all $m\in \mathbb{N}$, i.e.%
	\begin{equation}
	\mathcal{P}_{m}\left( 0\right) =\left\Vert u_{m}^{\prime \prime }\left(
	0\right) \right\Vert ^{2}+\left\Vert v_{m}^{\prime \prime }\left( 0\right)
	\right\Vert ^{2}+\left\Vert \tilde{u}_{1x}\right\Vert ^{2}+\left\Vert \tilde{%
		v}_{1x}\right\Vert ^{2}\leq \mathcal{P}_{0}.  \tag{3.16}
	\end{equation}
	Here the constant $\mathcal{P}_{0}$ is dependent of the initial data $\tilde{%
		u}_{0},$ $\tilde{v}_{0},$ $\tilde{u}_{1},$ $\tilde{v}_{1}$, the interior
	sources $f_{1},$ $f_{2}$, the functions $F_{1},$ $F_{2}$ and the given
	constants $r_{1},$ $r_{2},$ $\lambda _{1},$ $\lambda _{2}$.
	
	Estimating $\mathcal{P}_{m}\left( t\right) $ is almost similar to what have
	been done in \cite{13}, so we claim that by putting 
	$
	\mathcal{K}\left( T,\mathcal{F}\right) =\sup_{\left\vert y\right\vert ,%
		\text{ }\left\vert z\right\vert \leq \sqrt{C_{T}},\text{ }\left\vert \alpha
		\right\vert =2}\left\vert D^{\alpha }\mathcal{F}\left( y,z\right)
	\right\vert ,$
	there are possibilities to estimate the first three integrals, i.e. one can
	show that%
	\begin{equation}
	\mathcal{J}_{1}+\mathcal{J}_{2}\leq C_{T}+\int_{0}^{t}\mathcal{P}%
	_{m}\left( s\right) ds,%
	\tag{3.17}  \label{c17}
	\end{equation}%
	and by the standard Cauchy-Schwartz inequality, we obtain%
	\begin{equation}
	\mathcal{J}_{3}\leq C_{T}+\int_{0}^{t}\left( \left\Vert
	F_{1}^{\prime }\left( s\right) \right\Vert +\left\Vert F_{2}^{\prime }\left(
	s\right) \right\Vert \right) \mathcal{P}_{m}\left( s\right) ds.%
	\tag{3.18}  \label{c18}
	\end{equation}
	
	We can go through the last integral $\mathcal{J}_{4}$ by the following lemma.
	\begin{lemma}\label{lem:3.2}
		If one of the following cases is valid, which is
		\begin{equation*}
		\begin{cases}
		p_{1},\text{ }p_{2}\geq 2, \\ 
		2\leq q_{1},\text{ }q_{2}\leq 4,%
		\end{cases}%
		\mbox{or}\;%
		\begin{cases}
		p_{1},\text{ }p_{2}\in \left\{ 2\right\} \cup \left[ 3,\infty \right) , \\ 
		q_{1},\text{ }q_{2}>4,%
		\end{cases}%
		\end{equation*}%
		the integral $\mathcal{J}_{4}$ given by (\ref{c14}) can be bounded by
		\begin{equation}
		\mathcal{J}_{4}\leq C_{T}+\frac{1}{2}%
		\mathcal{P}_{m}\left( t\right) .%
		\tag{3.19}  \label{c19}
		\end{equation}
	\end{lemma}

	\begin{proof}[Proof of Lemma \ref{lem:3.2}]
		Since $\mathcal{J}_{4}$ can be divided into two integrals  and those can
		be separately estimated, we put ourselves into the following cases:\\
		\textbf{Case 1}: $2\leq q_{1},$ $q_{2}\leq 4,$ $p_{1}=p_{2}=2$ and $1\leq
		q_{1},$ $q_{2}\leq 4,$ $p_{1},$ $p_{2}>2$.\\
		\textbf{Case 2}: $q_{1},$ $q_{2}>4,$ $p_{1},$ $p_{2}\geq 3$ and $q_{1},$ $%
		q_{2}>4,$ $p_{1}=p_{2}=2$.
		
		Accordingly, our proof is presented below case by case.\\
		\textbf{Case 1.1:} $2\leq q_{1},$ $q_{2}\leq 4,$ $p_{1}=p_{2}=2.$
		We have%
		\begin{align}\tag{3.20} \label{c20}
		\mathcal{J}_{4} 
		&\leq \frac{4K_{1}}{q_{1}}\int_{0}^{t}\left\vert u_{m}^{\prime }\left(
		1,s\right) \right\vert ^{2-\frac{q_{1}}{2}}\left\vert \frac{\partial }{%
			\partial s}\left( \left\vert u_{m}^{\prime }\left( 1,s\right) \right\vert ^{%
			\frac{q_{1}}{2}-1}u_{m}^{\prime }\left( 1,s\right) \right) \right\vert ds 
		\notag \\
		&+\frac{4K_{2}}{q_{2}}\int_{0}^{t}\left\vert v_{m}^{\prime }\left(
		0,s\right) \right\vert ^{2-\frac{q_{2}}{2}}\left\vert \frac{\partial }{%
			\partial s}\left( \left\vert v_{m}^{\prime }\left( 0,s\right) \right\vert ^{%
			\frac{q_{2}}{2}-1}v_{m}^{\prime }\left( 0,s\right) \right) \right\vert ds. 
		\notag 
		\end{align}
		
		Due to (\ref{c5}) and (\ref{c15}) which read as%
		\begin{align} \tag{3.21}  \label{c21}
		\mathcal{S}_{m}\left( t\right) & \geq 2\mu _{1}\int_{0}^{t}\left\vert
		u_{m}^{\prime }\left( 1,s\right) \right\vert ^{q_{1}}ds+2\mu
		_{2}\int_{0}^{t}\left\vert v_{m}^{\prime }\left( 0,s\right) \right\vert
		^{q_{2}}ds,  
		\\ 
		\mathcal{P}_{m}\left( t\right) &\geq \frac{8\mu _{1}\left( q_{1}-1\right) }{%
			q_{1}^{2}}\int_{0}^{t}\left\vert \frac{\partial }{\partial s}\left(
		\left\vert u_{m}^{\prime }\left( 1,s\right) \right\vert ^{\frac{q_{1}}{2}%
			-1}u_{m}^{\prime }\left( 1,s\right) \right) \right\vert ^{2}ds  \tag{3.22}
		\label{c22}  \\
		&+\frac{8\mu _{2}\left( q_{2}-1\right) }{q_{2}^{2}}\int_{0}^{t}\left\vert 
		\frac{\partial }{\partial s}\left( \left\vert v_{m}^{\prime }\left(
		0,s\right) \right\vert ^{\frac{q_{2}}{2}-1}v_{m}^{\prime }\left( 0,s\right)
		\right) \right\vert ^{2}ds,  \notag
		\end{align}%
		together with (\ref{c10}), we bound $\mathcal{J}_{4}$ by
		\begin{align}\tag{3.23}  \label{c23}
		\mathcal{J}_{4} &\leq \frac{2K_{1}}{q_{1}}\int_{0}^{t}\left[ \frac{1%
		}{\delta _{1}}\left\vert u_{m}^{\prime }\left( 1,s\right) \right\vert
		^{4-q_{1}}+\delta _{1}\left\vert \frac{\partial }{\partial s}\left(
		\left\vert u_{m}^{\prime }\left( 1,s\right) \right\vert ^{\frac{q_{1}}{2}%
			-1}u_{m}^{\prime }\left( 1,s\right) \right) \right\vert ^{2}\right] ds \\ &
		+\frac{2K_{2}}{q_{2}}\int_{0}^{t}\left[ \frac{1}{\delta _{2}%
		}\left\vert v_{m}^{\prime }\left( 0,s\right) \right\vert ^{4-q_{2}}+\delta
		_{2}\left\vert \frac{\partial }{\partial s}\left( \left\vert v_{m}^{\prime
		}\left( 0,s\right) \right\vert ^{\frac{q_{2}}{2}-1}v_{m}^{\prime }\left(
		0,s\right) \right) \right\vert ^{2}\right] ds%
		\notag
		\\&\leq \frac{2K_{1}}{q_{1}\delta _{1}}\int_{0}^{t}\left[ 1+\left\vert
		u_{m}^{\prime }\left( 1,s\right) \right\vert ^{q_{1}}\right] ds+\frac{2K_{2}%
		}{q_{2}\delta _{2}}\int_{0}^{t}\left[ 1+\left\vert v_{m}^{\prime }\left(
		0,s\right) \right\vert ^{q_{2}}\right] ds \notag \\
		&+\left( \frac{\delta _{1}K_{1}q_{1}}{4\mu _{1}\left( q_{1}-1\right) }+%
		\frac{\delta _{2}K_{2}q_{2}}{4\mu _{2}\left( q_{2}-1\right) }\right) 
		\mathcal{P}_{m}\left( t\right) \notag\\
		&\leq \frac{2K_{1}}{q_{1}\delta _{1}}\left[ T+\frac{1}{2\mu _{1}}\mathcal{S}%
		_{m}\left( t\right) \right] +\frac{2K_{2}}{q_{2}\delta _{2}}\left[ T+\frac{1%
		}{2\mu _{2}}\mathcal{S}_{m}\left( t\right) \right] \notag\\
		&+\left( \frac{\delta _{1}K_{1}q_{1}}{4\mu _{1}\left( q_{1}-1\right) }+%
		\frac{\delta _{2}K_{2}q_{2}}{4\mu _{2}\left( q_{2}-1\right) }\right) 
		\mathcal{P}_{m}\left( t\right) \notag\\
		&\leq \left( \frac{1}{\delta _{1}}+\frac{1}{\delta _{2}}\right)
		C_{T}+\left( \frac{\delta _{1}K_{1}q_{1}}{4\mu _{1}\left( q_{1}-1\right) }+%
		\frac{\delta _{2}K_{2}q_{2}}{4\mu _{2}\left( q_{2}-1\right) }\right) 
		\mathcal{P}_{m}\left( t\right) . \notag
		\end{align}%
		
		We use $a^{4-q}\leq 1+a^{q}$ for all $a\geq 0,$ $2\leq q\leq 4$ and $2ab\leq
		\delta a^{2}+\delta ^{-1}b^{2}$ for all $a,$ $b\geq 0$ and $\delta >0$.
		Thus, to deduce (\ref{c19}) from (\ref{c23}) we choose%
		\begin{equation*}
		\delta _{1}=\delta _{2}\leq \frac{2\mu _{1}\mu _{2}\left( q_{1}-1\right)
			\left( q_{2}-1\right) }{K_{1}\mu _{2}q_{1}\left( q_{2}-1\right) +K_{2}\mu
			_{1}q_{2}\left( q_{1}-1\right) }.%
		\end{equation*}
		\textbf{Case 1.2:} $2\leq q_{1},$ $q_{2}\leq 4,$ $p_{1},$ $p_{2}>2.$ By simple computations, we have
		\begin{align*}
		& \mathcal{J}_{4}\\
		&=2K_{1}\left( p_{1}-1\right) \int_{0}^{t}\left\vert u_{m}\left( 1,s\right)
		\right\vert ^{p_{1}-2}\left\vert u_{m}^{\prime }\left( 1,s\right)
		\right\vert ^{1-\frac{q_{1}}{2}}u_{m}^{\prime }\left( 1,s\right) \left\vert
		u_{m}^{\prime }\left( 1,s\right) \right\vert ^{\frac{q_{1}}{2}%
			-1}u_{m}^{\prime \prime }\left( 1,s\right) ds \\
		&+2K_{2}\left( p_{2}-1\right) \int_{0}^{t}\left\vert v_{m}\left( 0,s\right)
		\right\vert ^{p_{2}-2}\left\vert v_{m}^{\prime }\left( 0,s\right)
		\right\vert ^{1-\frac{q_{2}}{2}}v_{m}^{\prime }\left( 0,s\right) \left\vert
		v_{m}^{\prime }\left( 0,s\right) \right\vert ^{\frac{q_{2}}{2}%
			-1}v_{m}^{\prime \prime }\left( 0,s\right) ds \\
		&\le \frac{4}{q_{1}}K_{1}\left( p_{1}-1\right) C_{T}^{%
			\frac{p_{1}}{2}-1}\int_{0}^{t}\left\vert u_{m}^{\prime }\left( 1,s\right)
		\right\vert ^{2-\frac{q_{1}}{2}}\left\vert \frac{\partial }{\partial s}%
		\left( \left\vert u_{m}^{\prime }\left( 1,s\right) \right\vert ^{\frac{q_{1}%
			}{2}-1}u_{m}^{\prime }\left( 1,s\right) \right) \right\vert ds \\
		&+\frac{4}{q_{2}}K_{2}\left( p_{2}-1\right) C_{T}^{\frac{p_{2}}{2}%
			-1}\int_{0}^{t}\left\vert v_{m}^{\prime }\left( 0,s\right) \right\vert ^{2-%
			\frac{q_{2}}{2}}\left\vert \frac{\partial }{\partial s}\left( \left\vert
		v_{m}^{\prime }\left( 0,s\right) \right\vert ^{\frac{q_{2}}{2}%
			-1}v_{m}^{\prime }\left( 0,s\right) \right) \right\vert ds,
		\end{align*}
		using the fact that $\left\vert u_{m}\left( 1,s\right) \right\vert +\left\vert
		v_{m}\left( 0,s\right) \right\vert \leq \sqrt{C_{T}}$ (cf. Lemma \ref{lem:2.2}), (\ref%
		{c5}) and (\ref{c10}).
		This way gives us back to (\ref{c20}) in the previous case. Thus, (\ref{c19})
		holds.\\
		\textbf{Case 2.1:} $q_{1},$ $q_{2}>4,$ $p_{1},$ $p_{2}\geq 3.$ We start by using the integration by parts
		\begin{align*}
		\mathcal{J}_{4} 
		&=K_{1}\left( p_{1}-1\right) \left\vert u_{m}\left( 1,t\right) \right\vert
		^{p_{1}-2}\left\vert u_{m}^{\prime }\left( 1,t\right) \right\vert
		^{2}-K_{1}\left( p_{1}-2\right) \left\vert \tilde{u}_{0}\left( 1\right)
		\right\vert ^{p_{1}-2}\tilde{u}_{1}^{2}\left( 1\right) \\
		&-K_{1}\left( p_{1}-1\right) \left( p_{1}-2\right) \int_{0}^{t}\left\vert
		u_{m}\left( 1,s\right) \right\vert ^{p_{1}-4}u_{m}\left( 1,s\right) \left(
		u_{m}^{\prime }\left( 1,s\right) \right) ^{3}ds \\
		&+K_{2}\left( p_{2}-1\right) \left\vert v_{m}\left( 0,t\right) \right\vert
		^{p_{2}-2}\left\vert v_{m}^{\prime }\left( 0,t\right) \right\vert
		^{2}-K_{2}\left( p_{2}-2\right) \left\vert \tilde{v}_{0}\left( 0\right)
		\right\vert ^{p_{2}-2}\tilde{v}_{1}^{2}\left( 0\right) \\
		&-K_{2}\left( p_{2}-1\right) \left( p_{2}-2\right) \int_{0}^{t}\left\vert
		v_{m}\left( 0,s\right) \right\vert ^{p_{2}-4}v_{m}\left( 0,s\right) \left(
		v_{m}^{\prime }\left( 0,s\right) \right) ^{3}ds.
		\end{align*}
		
		Based on (\ref{c5}) and (\ref{c10}), we first estimate $\mathcal{J}_{4} $ as
		follows:
		\begin{align*}
		\mathcal{J}_{4} &\leq K_{1}\left( p_{1}-1\right) \mathcal{S}_{m}^{\frac{%
				p_{1}}{2}-1}\left\vert u_{m}^{\prime }\left( 1,t\right) \right\vert
		^{2}+K_{1}\left( p_{1}-1\right) \left( p_{1}-2\right) \int_{0}^{t}\mathcal{S}%
		_{m}^{\frac{p_{1}-3}{2}}\left( s\right) \left\vert u_{m}^{\prime }\left(
		1,s\right) \right\vert ^{3}ds \\
		&+K_{2}\left( p_{2}-1\right) \mathcal{S}_{m}^{\frac{p_{2}}{2}-1}\left\vert
		v_{m}^{\prime }\left( 0,t\right) \right\vert ^{2}+K_{2}\left( p_{2}-1\right)
		\left( p_{2}-2\right) \int_{0}^{t}\mathcal{S}_{m}^{\frac{p_{2}-3}{2}}\left(
		s\right) \left\vert v_{m}^{\prime }\left( 0,s\right) \right\vert ^{3}ds \\
		&\leq C_{T}\left( \left\vert u_{m}^{\prime }\left( 1,t\right) \right\vert
		^{2}+\left\vert v_{m}^{\prime }\left( 0,t\right) \right\vert ^{2}\right)
		+C_{T}\int_{0}^{t}\left[ \left\vert u_{m}^{\prime }\left( 1,s\right)
		\right\vert ^{3}+\left\vert v_{m}^{\prime }\left( 0,s\right) \right\vert ^{3}%
		\right] ds.
		\end{align*}
		
		Secondly, it is immediate to see that using the inequality $a^{3}\leq 1+a^{q}
		$ for all $a\geq 0$ and $q\geq 3$, together with (\ref{c10}) and (\ref{c21}%
		), we arrive at%
		\begin{equation*}
		\int_{0}^{t}\left[ \left\vert u_{m}^{\prime }\left( 1,s\right) \right\vert
		^{3}+\left\vert v_{m}^{\prime }\left( 0,s\right) \right\vert ^{3}\right]
		ds\leq T+\frac{1}{2}\left( \frac{1}{\mu _{1}}+\frac{1}{\mu _{2}}\right) 
		\mathcal{S}_{m}\left( t\right) \leq C_{T}.
		\end{equation*}
		
		One also deduces from%
		\begin{align*}
		& \left\vert u_{m}^{\prime }\left( 1,t\right) \right\vert ^{\frac{q_{1}}{2}%
			-1}u_{m}^{\prime }\left( 1,t\right) \\&=\left\vert \tilde{u}_{1m}\left(
		1\right) \right\vert ^{\frac{q_{1}}{2}-1}\tilde{u}_{1m}\left( 1\right)
		+\int_{0}^{t}\frac{\partial }{\partial s}\left( \left\vert u_{m}^{\prime
		}\left( 1,s\right) \right\vert ^{\frac{q_{1}}{2}-1}u_{m}^{\prime }\left(
		1,s\right) \right) ds, \\
		& \left\vert v_{m}^{\prime }\left( 0,t\right) \right\vert ^{\frac{q_{2}}{2}%
			-1}v_{m}^{\prime }\left( 0,t\right) \\&=\left\vert \tilde{v}_{1m}\left(
		0\right) \right\vert ^{\frac{q_{2}}{2}-1}\tilde{v}_{1m}\left( 0\right)
		+\int_{0}^{t}\frac{\partial }{\partial s}\left( \left\vert v_{m}^{\prime
		}\left( 0,s\right) \right\vert ^{\frac{q_{2}}{2}-1}v_{m}^{\prime }\left(
		0,s\right) \right) ds,
		\end{align*}%
		and from the elementary inequality $\left( a+b\right) ^{2}\leq 2\left(
		a^{2}+b^{2}\right) $ for all $a,$ $b\geq 0$ with H\"{o}lder's inequality and
		(\ref{c22}) that
		\begin{align*}
		\left\vert u_{m}^{\prime }\left( 1,t\right) \right\vert ^{q_{1}}& \leq
		2\left\vert \tilde{u}_{1m}\left( 1\right) \right\vert
		^{q_{1}}+2t\int_{0}^{t}\left\vert \frac{\partial }{\partial s}\left(
		\left\vert u_{m}^{\prime }\left( 1,s\right) \right\vert ^{\frac{q_{1}}{2}%
			-1}u_{m}^{\prime }\left( 1,s\right) \right) \right\vert ^{2}ds \\ &
		\leq 2\left\vert \tilde{u}_{1m}\left( 1\right)
		\right\vert ^{q_{1}}+\frac{q_{1}^{2}T}{4\mu _{1}\left( q_{1}-1\right) }%
		\mathcal{P}_{m}\left( t\right) .%
		\end{align*}
		
		In the same vein, we get%
		\begin{equation*}
		\left\vert v_{m}^{\prime }\left( 0,t\right) \right\vert ^{q_{2}}\leq
		2\left\vert \tilde{v}_{1m}\left( 0\right) \right\vert ^{q_{2}}+\frac{%
			q_{2}^{2}T}{4\mu _{2}\left( q_{2}-1\right) }\mathcal{P}_{m}\left( t\right) .%
		\end{equation*}
		
		Using the inequalities%
		\begin{align*}
		&\left(a+b\right) ^{\frac{2}{q}}\leq a^{\frac{2}{q}}+b^{\frac{2}{q}},\text{
			\ }\forall a,b\geq 0,\text{ }\forall q\geq 2, \\ &
		ab\leq \left(1-\frac{2}{q}\right) \delta ^{-\frac{q}{q-2}}a^{\frac{q}{q-2}%
		}+\frac{2}{q}\delta ^{\frac{q}{2}}b^{\frac{q}{2}},\text{ \ }\forall a,b\geq 0,\text{ }\forall q>2,\text{ }\delta >0,%
		\end{align*}%
		we therefore obtain%
		\begin{align*}
		&C_{T}\left( \left\vert u_{m}^{\prime }\left( 1,t\right) \right\vert
		^{2}+\left\vert v_{m}^{\prime }\left( 0,t\right) \right\vert ^{2}\right)\\&
		\leq C_{T}\left( 2\left\vert \tilde{u}_{1m}\left( 1\right) \right\vert
		^{q_{1}}+\frac{q_{1}^{2}T}{4\mu _{1}\left( q_{1}-1\right) }\mathcal{P}%
		_{m}\left( t\right) \right) ^{\frac{2}{q_{1}}} \\ &
		+C_{T}\left( 2\left\vert \tilde{v}_{1m}\left( 0\right) \right\vert
		^{q_{2}}+\frac{q_{2}^{2}T}{4\mu _{2}\left( q_{2}-1\right) }\mathcal{P}%
		_{m}\left( t\right) \right) ^{\frac{2}{q_{2}}} \\ &
		\leq C_{T}\left[ 2^{\frac{2}{q_{1}}}\left\vert \tilde{u}_{1m}\left(
		1\right) \right\vert ^{2}+2^{\frac{2}{q_{2}}}\left\vert \tilde{v}_{1m}\left(
		0\right) \right\vert ^{2}\right]  \\ &
		+C_{T}\left[ \left( \frac{q_{1}^{2}T}{4\mu _{1}\left( q_{1}-1\right) }%
		\right) ^{\frac{2}{q_{1}}}\mathcal{P}_{m}^{\frac{2}{q_{1}}}\left( t\right)
		+\left( \frac{q_{2}^{2}T}{4\mu _{2}\left( q_{2}-1\right) }\right) ^{\frac{2}{%
				q_{2}}}\mathcal{P}_{m}^{\frac{2}{q_{2}}}\left( t\right) \right]  \\ &
		\leq C_{0}+C_{T}\left[ \left( 1-\frac{2}{q_{1}}\right) \delta _{1}^{-%
			\frac{q_{1}}{q_{1}-2}}\left( \frac{q_{1}^{2}T}{4\mu _{1}\left(
			q_{1}-1\right) }\right) ^{\frac{4}{q_{1}\left( q_{1}-2\right) }}+\frac{2}{%
			q_{1}}\delta _{1}^{\frac{q_{1}}{2}}\mathcal{P}_{m}\left( t\right) \right] 
		\\ &
		+C_{T}\left[ \left( 1-\frac{2}{q_{2}}\right) \delta _{2}^{-\frac{q_{2}%
			}{q_{2}-2}}\left( \frac{q_{2}^{2}T}{4\mu _{2}\left( q_{2}-1\right) }\right)
		^{\frac{4}{q_{2}\left( q_{2}-2\right) }}+\frac{2}{q_{2}}\delta _{2}^{\frac{%
				q_{2}}{2}}\mathcal{P}_{m}\left( t\right) \right]  \\ &
		\leq C_{T}\left( \delta _{1},\delta _{2}\right) +2\left( \frac{\delta
			_{1}^{\frac{q_{1}}{2}}}{q_{1}}+\frac{\delta _{2}^{\frac{q_{2}}{2}}}{q_{2}}%
		\right) \mathcal{P}_{m}\left( t\right) .%
		\end{align*}
		
		Hence, to deduce (\ref{c19}) we choose $\delta =\delta _{1}=\delta _{2}>0$
		such that $q_{2}\delta ^{\frac{q_{1}}{2}}+q_{1}\delta ^{\frac{q_{2}}{2}}\leq 
		\frac{q_{1}q_{2}}{2}$.\\
		\textbf{Case 2.2:} $q_{1},$ $q_{2}>4,$ $p_{1}=p_{2}=2.$ By the same arguments exploited in the previous case, here we can state that%
		\begin{align*}
		\mathcal{J}_{4} 
		&=K_{1}\int_{0}^{t}\frac{d}{ds}\left( \left\vert u_{m}^{\prime }\left(
		1,s\right) \right\vert ^{2}\right) ds+K_{2}\int_{0}^{t}\frac{d}{ds}\left(
		\left\vert v_{m}^{\prime }\left( 0,s\right) \right\vert ^{2}\right) ds \\
		&=K_{1}\left( \left\vert u_{m}^{\prime }\left( 1,t\right) \right\vert ^{2}-%
		\tilde{u}_{1m}^{2}\left( 1\right) \right) +K_{2}\left( \left\vert
		v_{m}^{\prime }\left( 0,t\right) \right\vert ^{2}-\tilde{v}_{1m}^{2}\left(
		0\right) \right) \\
		&\leq K_{1}\left\vert u_{m}^{\prime }\left( 1,t\right) \right\vert
		^{2}+K_{2}\left\vert v_{m}^{\prime }\left( 0,t\right) \right\vert ^{2},
		\end{align*}%
		also leads to (\ref{c19}).
		
		Hence, we complete the proof of Lemma \ref{lem:3.2}.	
	\end{proof}

	Now, combining (\ref{c17}), (\ref{c18}), and (\ref{c19}) we are in a great
	position to obtain%
	\begin{equation*}
	\mathcal{P}_{m}\left( t\right) \leq 2\mathcal{P}_{0}+6C_{T}+2\int%
	_{0}^{t}\left( 1+\left\Vert F_{1}^{\prime }\left( s\right)
	\right\Vert +\left\Vert F_{2}^{\prime }\left( s\right) \right\Vert \right) 
	\mathcal{P}_{m}\left( s\right) ds.%
	\end{equation*}
	
	Thanks to Gronwall's inequality, we conclude that%
	\begin{equation}
	\mathcal{P}_{m}\left( t\right) \leq \left( 2\mathcal{P}_{0}+6C_{T}\right) %
	\mbox{exp}\left[ \int_{0}^{t}\left( 1+\left\Vert F_{1}^{\prime
	}\left( s\right) \right\Vert +\left\Vert F_{2}^{\prime }\left( s\right)
	\right\Vert \right) ds\right] \leq C_{T},%
	\tag{3.24}  \label{c24}
	\end{equation}%
	for all $m\in \mathbb{N}$ and $t\in \left[ 0,T_{\ast }\right] $.\\
	\textit{Step 3. Passing to the limit.} The existence of solution in the
	interval $\left[ 0,T_{\ast }\right] $ is now approaching.
	To summarize,
	using the Banach-Alaoglu theorem (see, e.g., \cite{5}), the uniform
	bounds with respect to $m$, as stated in the above results (\ref{c5}), (\ref%
	{c10}), (\ref{c15}), and (\ref{c24}), imply that one can extract a further
	subsequence (which we relabel with the index $m$ if necessary) such that
	\begin{equation}
	\left( u_{m},v_{m}\right) \rightarrow \left( u,v\right) \;\mbox{weak-*
		in}\;L^{\infty }\left( 0,T_{\ast };\mathbb{V}_{1}\times \mathbb{V}%
	_{2}\right) ,  \tag{3.25}  \label{c25}
	\end{equation}%
	\begin{equation}
	\left( u_{m}^{\prime },v_{m}^{\prime }\right) \rightarrow \left( u^{\prime
	},v^{\prime }\right) \;\mbox{weak in}\;L^{r_{1}}\left( Q_{T_{\ast
	}}\right) \times L^{r_{2}}\left( Q_{T_{\ast }}\right), \mbox{weak-* in}\;L^{\infty }\left( 0,T_{\ast };\mathbb{V}_{1}\times \mathbb{V}%
	_{2}\right) ,  \tag{3.26}  \label{c26}
	\end{equation}%
	\begin{equation}
	\left( u_{m}^{\prime \prime },v_{m}^{\prime \prime }\right) \rightarrow
	\left( u^{\prime \prime },v^{\prime \prime }\right) \;\mbox{weak-* in}%
	\;L^{\infty }\left( 0,T_{\ast };L^{2}\times L^{2}\right) ,  \tag{3.27}
	\label{c27}
	\end{equation}%
	\begin{equation}
	\left( u_{m}\left( 1,\cdot \right) ,v_{m}\left( 0,\cdot \right) \right)
	\rightarrow \left( u\left( 1,\cdot \right) ,v\left( 0,\cdot \right) \right)
	\;\mbox{weak in}\;W^{1,q_{1}}\left( 0,T_{\ast }\right) \times
	W^{1,q_{2}}\left( 0,T_{\ast }\right) ,  \tag{3.28}
	\end{equation}%
	\begin{equation}
	\left( u_{m}^{\prime }\left( 1,\cdot \right) ,v_{m}^{\prime }\left( 0,\cdot
	\right) \right) \rightarrow \left( u^{\prime }\left( 1,\cdot \right)
	,v^{\prime }\left( 0,\cdot \right) \right) \;\mbox{weak in}%
	\;L^{q_{1}}\left( 0,T_{\ast }\right) \times L^{q_{2}}\left( 0,T_{\ast
	}\right) ,  \tag{3.29}  \label{c29}
	\end{equation}%
	\begin{equation}
	\left( \left\vert u_{m}^{\prime }\left( 1,\cdot \right) \right\vert ^{\frac{%
			q_{1}}{2}-1}u_{m}^{\prime }\left( 1,\cdot \right) ,\left\vert
	v_{m}^{\prime }\left( 0,\cdot \right) \right\vert ^{\frac{q_{2}}{2}%
		-1}v_{m}^{\prime }\left( 0,\cdot \right) \right) \rightarrow \left( \chi
	_{1},\chi _{2}\right) \;\mbox{weak in}\;[H^{1}\left( 0,T_{\ast }\right)]^2,  \tag{3.30}
	\end{equation}%
	\begin{equation}
	\left( \frac{\partial }{\partial t}\left( \left\vert u_{m}^{\prime
	}\right\vert ^{\frac{r_{1}}{2}-1}u_{m}^{\prime }\right) ,\frac{%
		\partial }{\partial t}\left( \left\vert v_{m}^{\prime }\right\vert ^{\frac{%
			r_{2}}{2}-1}v_{m}^{\prime }\right) \right) \rightarrow \left( \chi _{3},\chi
	_{4}\right) \;\mbox{weak in}\;[L^{2}\left( Q_{T_{\ast }}\right)]^2 .  \tag{3.31}  \label{c31}
	\end{equation}
	
	Furthermore, by the Aubin-Lions compactness theorem in combination with the
	imbeddings $H^{2}\left( 0,T_{\ast }\right) \hookrightarrow C^{1}\left( \left[
	0,T_{\ast }\right] \right) ,$ $H^{1}\left( 0,T_{\ast }\right)
	\hookrightarrow C\left( \left[ 0,T_{\ast }\right] \right) ,$ $%
	W^{1,q_{1}}\left( 0,T_{\ast }\right) \hookrightarrow C\left( \left[
	0,T_{\ast }\right] \right) ,$ $W^{1,q_{2}}\left( 0,T_{\ast }\right)
	\hookrightarrow C\left( \left[ 0,T_{\ast }\right] \right) $, it is
	straightforward to go on extracting from the weak convergence results (\ref%
	{c25})-(\ref{c31}) a subsequence $\left\{ \left( u_{m},v_{m}\right) \right\} 
	$ such that%
	\begin{equation}
	\left( u_{m},v_{m}\right) \rightarrow \left( u,v\right) \;\mbox{strong in}%
	\;[L^{2}\left( Q_{T_{\ast }}\right)]^2 \;%
	\mbox{and almost everywhere in \ensuremath{Q_{T_{*}}}},  \tag{3.32}
	\label{c32}
	\end{equation}%
	\begin{equation}
	\left( u_{m}^{\prime },v_{m}^{\prime }\right) \rightarrow \left( u^{\prime
	},v^{\prime }\right) \;\mbox{strong in}\;[L^{2}\left( Q_{T_{\ast }}\right)]^2 \;%
	\mbox{and
		almost everywhere in \ensuremath{Q_{T_{*}}}},  \tag{3.33}  \label{c33}
	\end{equation}%
	\begin{equation}
	\left( u_{m}\left( 1,\cdot \right) ,v_{m}\left( 0,\cdot \right) \right)
	\rightarrow \left( u\left( 1,\cdot \right) ,v\left( 0,\cdot \right) \right)
	\;\mbox{strong in}\;[C\left( \left[ 0,T_{\ast }\right] \right)]^2,  \tag{3.34}  \label{c34}
	\end{equation}%
	\begin{equation}
	\left( \left\vert u_{m}^{\prime }\left( 1,\cdot \right) \right\vert ^{\frac{%
			q_{1}}{2}-1}u_{m}^{\prime }\left( 1,\cdot \right) ,\left\vert
	v_{m}^{\prime }\left( 0,\cdot \right) \right\vert ^{\frac{q_{2}}{2}%
		-1}v_{m}^{\prime }\left( 0,\cdot \right) \right) \rightarrow \left( \chi
	_{1},\chi _{2}\right) \;\mbox{strong in}\;[C\left( \left[ 0,T_{\ast }%
	\right] \right)]^2 . 
	\tag{3.35}  \label{c35}
	\end{equation}
	
	Now we have to show the convergence of the nonlinear terms including damping
	and interior sources. In fact, using the continuity of $f_{1}$, one
	deduces that%
	\begin{equation*}
	f_{1}\left( u_{m},v_{m}\right) \rightarrow f_{1}\left( u,v\right) \;%
	\mbox{almost everywhere in}\;Q_{T_{\ast }}.
	\end{equation*}
	Obverse that $\left\Vert f_{1}\left( u_{m},v_{m}\right) \right\Vert
	_{L^{2}\left( Q_{T_{\ast }}\right) }$ is bounded by $\sqrt{T_{\ast }}{%
		\sup_{\left\vert y\right\vert ,\text{ }\left\vert z\right\vert
			\leq \sqrt{C_{T}}}\left\vert f_{1}\left( y,z\right) \right\vert }$ which
	cannot go to infinity, and together with \cite[Lemma 1.3]{10}, one continues
	to obtain the following%
	\begin{equation}
	f_{1}\left( u_{m},v_{m}\right) \rightarrow f_{1}\left( u,v\right) \;%
	\mbox{weak in}\;L^{2}\left( Q_{T_{\ast }}\right) ,  \tag{3.36}  \label{c36}
	\end{equation}%
	and%
	\begin{equation}
	f_{2}\left( u_{m},v_{m}\right) \rightarrow f_{2}\left( u,v\right) \;%
	\mbox{weak in}\;L^{2}\left( Q_{T_{\ast }}\right) .  \tag{3.37}  \label{c37}
	\end{equation}
	
	It remains to see the weak convergence of damping terms. Thanks to the
	inequality%
	\begin{equation}
	\left\vert \Psi _{r}\left( z_{1}\right) -\Psi _{r}\left( z_{2}\right)
	\right\vert \leq \left( r-1\right) C^{r-2}\left\vert z_{1}-z_{2}\right\vert ,%
	\text{ \ }\forall z_{1},z_{2}\in \left[ -C,C\right] ,\text{ }C>0,%
	\text{ }r\geq 2  \tag{3.38}  \label{c38}
	\end{equation}%
	in accordance with (\ref{c15}), (\ref{c24}) and (\ref{c33}), one easily
	obtains%
	\begin{equation}
	\left( \Psi _{r_{1}}\left( u_{m}^{\prime }\right) ,\Psi _{r_{2}}\left(
	v_{m}^{\prime }\right) \right) \rightarrow \left( \Psi _{r_{1}}\left(
	u^{\prime }\right) ,\Psi _{r_{2}}\left( v^{\prime }\right) \right) \;%
	\mbox{strong in}\;[L^{2}\left( Q_{T_{\ast }}\right)]^2.  \tag{3.39}  \label{c39}
	\end{equation}
	
	It is then worthwhile to mention that (\ref{c34}) gives%
	\begin{equation}
	\left( \Psi _{p_{1}}\left( u_{m}\left( 1,\cdot \right) \right) ,\Psi
	_{p_{2}}\left( v_{m}\left( 0,\cdot \right) \right) \right) \rightarrow
	\left( \Psi _{p_{1}}\left( u\left( 1,\cdot \right) \right) ,\Psi
	_{p_{2}}\left( v\left( 0,\cdot \right) \right) \right) \;\mbox{strong in}%
	\;[C\left( \left[ 0,T_{\ast }\right] \right)]^2,  \tag{3.40}  \label{c40}
	\end{equation}%
	by the continuity of $\Psi _{p_{i}}$ for $i=1,2$, and (\ref{c35})
	provides that%
	\begin{equation}
	\left( u_{m}^{\prime }\left( 1,\cdot \right) ,v_{m}^{\prime }\left( 0,\cdot
	\right) \right) \rightarrow \left( \left\vert \chi _{1}\right\vert ^{\frac{2%
		}{q_{1}}-1}\chi _{1},\left\vert \chi _{2}\right\vert ^{\frac{2}{q_{2}}%
		-1}\chi _{2}\right) \;\mbox{strong in}\;[C\left( \left[ 0,T_{\ast }
	\right] \right)]^2. 
	\tag{3.41}  \label{c41}
	\end{equation}
	
	Thus, we take (\ref{c29}) and (\ref{c41}) to gain%
	\begin{equation}
	\left( \left\vert \chi _{1}\right\vert ^{\frac{2}{q_{1}}-1}\chi
	_{1},\left\vert \chi _{2}\right\vert ^{\frac{2}{q_{2}}-1}\chi _{2}\right)
	=\left( u^{\prime }\left( 1,\cdot \right) ,v^{\prime }\left( 0,\cdot \right)
	\right) %
	\tag{3.42}  \label{c42}
	\end{equation}%
	by virtue of the uniqueness of convergence.
	
	In the same vein, we claim%
	\begin{equation}
	\left( \Psi _{q_{1}}\left( u_{m}^{\prime }\left( 1,\cdot \right) \right)
	,\Psi _{q_{2}}\left( v_{m}^{\prime }\left( 0,\cdot \right) \right) \right)
	\rightarrow \left( \Psi _{q_{1}}\left( u^{\prime }\left( 1,\cdot \right)
	\right) ,\Psi _{q_{2}}\left( v^{\prime }\left( 0,\cdot \right) \right)
	\right) \;\mbox{strong in}\;[C\left( \left[ 0,T_{\ast }\right] \right)]^2  \tag{3.43}  \label{c43}
	\end{equation}%
	from (\ref{c41}) and (\ref{c42}).
	
	From here on, combining (\ref{c25})-(\ref{c27}), (\ref{c32})-(\ref{c37}), (%
	\ref{c39}), (\ref{c40}), (\ref{c43}) is sufficient to pass to the limit in (\ref%
	{c3}) and then to show that $\left( u,v\right) $ satisfies the problem $%
	\left( P\right) $. In addition, one can use (\ref{c25})-(\ref{c27}), (\ref%
	{b1}) and \textbf{(}$\mathbf{A}_{2}$\textbf{)} to prove that%
	\begin{equation*}
	\begin{cases}
	u_{xx}=u_{tt}+\lambda _{1}\Psi _{r_{1}}\left( u_{t}\right) -f_{1}\left(
	u,v\right) -F_{1}\in L^{\infty }\left( 0,T_{\ast };L^{2}\right) , \\ 
	v_{xx}=v_{tt}+\lambda _{2}\Psi _{r_{2}}\left( v_{t}\right) -f_{2}\left(
	u,v\right) -F_{2}\in L^{\infty }\left( 0,T_{\ast };L^{2}\right),%
	\end{cases}%
	\end{equation*}%
	which verifies $\left( u,v\right) \in L^{\infty }\left( 0,T_{\ast};\left( \mathbb{V}%
	_{1}\cap H^{2}\right) \times \left( \mathbb{V}_{2}\cap H^{2}\right) \right) $
	and completes the proof of the existence of a local weak solution.\\
	\textit{Step 4. Uniqueness of the solution. } Suppose $\left(
	u_{1},v_{1}\right) $ and $\left( u_{2},v_{2}\right) $ are two solutions to $%
	\left( P\right) $ in the interval $[0,T_{\ast }]$, which is devoted to the
	case $q_{1}=q_{2}=2$ and $p_{1},$ $p_{2}\geq 2$. Going along with the
	same initial data $\left( \tilde{u}_{0},\tilde{u}_{1}\right) $ and $\left( 
	\tilde{v}_{0},\tilde{v}_{1}\right) $, we prove that these solutions must be equal.
	
	Define $\left( u,v\right) :=\left( u_{1}-u_{2},v_{1}-v_{2}\right) $ and
	based on (\ref{b1}) and (\ref{b2}), these quantities satisfy the following
	system:
	\begin{align}	\tag{3.44}  \label{c44}
	\langle u^{\prime \prime }\left( t\right) ,\phi \rangle +\langle u_{x}\left(
	t\right) ,\phi _{x}\rangle +\lambda _{1}\langle \Psi _{r_{1}}\left(
	u_{1}^{\prime }\left( t\right) \right) -\Psi _{r_{1}}\left( u_{2}^{\prime
	}\left( t\right) \right) ,\phi \rangle +\mu _{1}u^{\prime }\left( 1,t\right)
	\phi \left( 1\right) \\ 
	=K_{1}\left[ \Psi _{p_{1}}\left( u_{1}\left( 1,t\right) \right) -\Psi
	_{p_{1}}\left( u_{2}\left( 1,t\right) \right) \right] \phi \left( 1\right)
	+\langle f_{1}\left( u_{1},v_{1}\right) -f_{2}\left( u_{2},v_{2}\right)
	,\phi \rangle ,%
	\nonumber 
	\end{align}%
	\begin{align}\tag{3.45}  \label{c45}
	\langle v^{\prime \prime }\left( t\right) ,\tilde{\phi}\rangle +\langle
	v_{x}\left( t\right) ,\tilde{\phi}_{x}\rangle +\lambda _{2}\langle \Psi
	_{r_{2}}\left( v_{1}^{\prime }\left( t\right) \right) -\Psi _{r_{2}}\left(
	v_{2}^{\prime }\left( t\right) \right) ,\tilde{\phi}\rangle +\mu
	_{2}v^{\prime }\left( 0,t\right) \tilde{\phi}\left( 0\right) \\ 
	=K_{2}\left[ \Psi _{p_{2}}\left( v_{1}\left( 0,t\right) \right) -\Psi
	_{p_{2}}\left( v_{2}\left( 0,t\right) \right) \right] \tilde{\phi}\left(
	0\right) +\langle f_{2}\left( u_{1},v_{1}\right) -f_{2}\left(
	u_{2},v_{2}\right) ,\tilde{\phi}\rangle ,%
	\nonumber 
	\end{align}%
	for all $(\phi ,\tilde{\phi})\in \mathbb{V}_{1}\times \mathbb{V}_{2}$. We
	endow this system with the initial conditions%
	\begin{equation*}
	u\left( 0\right) =v\left( 0\right) =u^{\prime }\left( 0\right) =v^{\prime
	}\left( 0\right) =0.
	\end{equation*}
	
	Taking into account $(\phi ,\tilde{\phi})=\left( u^{\prime },v^{\prime
	}\right) $ in (\ref{c44}) and (\ref{c45}), then integrating with respect to $%
	t$, we obtain the following:%
	\begin{align}\tag{3.46}  \label{c46}
	&\mathcal{W}\left( t\right) \\ &=2\int_{0}^{t}\langle f_{1}\left(
	u_{1},v_{1}\right) -f_{1}\left( u_{2},v_{2}\right) ,u^{\prime }\left(
	s\right) \rangle ds+2\int_{0}^{t}\langle f_{2}\left( u_{1},v_{1}\right)
	-f_{2}\left( u_{2},v_{2}\right) ,v^{\prime }\left( s\right) \rangle
	ds  \nonumber \\ &
	 +2K_{1}\int_{0}^{t}\left[ \Psi _{p_{1}}\left( u_{1}\left(
	1,s\right) \right) -\Psi _{p_{1}}\left( u_{2}\left( 1,s\right) \right) %
	\right] u^{\prime }\left( 1,s\right) ds \nonumber \\ 
	&+2K_{2}\int_{0}^{t}\left[ \Psi _{p_{2}}\left( v_{1}\left(
	0,s\right) \right) -\Psi _{p_{2}}\left( v_{2}\left( 0,s\right) \right) %
	\right] v^{\prime }\left( 0,s\right) ds \nonumber \\ &
	 =\mathcal{K}_{1}+\mathcal{K}_{2}+\mathcal{K}_{3}+\mathcal{K}%
	_{4}, \nonumber%
	\end{align}%
	where we have denoted by
	\begin{align}\tag{3.47}  \label{c47}
	\mathcal{W}\left( t\right) &:=\left\Vert u^{\prime }\left( t\right)
	\right\Vert ^{2}+\left\Vert v^{\prime }\left( t\right) \right\Vert
	^{2}+\left\Vert u_{x}\left( t\right) \right\Vert ^{2}+\left\Vert v_{x}\left(
	t\right) \right\Vert ^{2}  \\
	&+2\lambda _{1}\int_{0}^{t}\langle \Psi _{r_{1}}\left( u_{1}^{\prime
	}\left( s\right) \right) -\Psi _{r_{1}}\left( u_{2}^{\prime }\left( s\right)
	\right) ,u^{\prime }\left( s\right) \rangle ds+2\mu
	_{1}\int_{0}^{t}\left\vert u^{\prime }\left( 1,s\right) \right\vert ^{2}ds 
	\notag \\
	&+2\lambda _{2}\int_{0}^{t}\langle \Psi _{r_{2}}\left( v_{1}^{\prime
	}\left( s\right) \right) -\Psi _{r_{2}}\left( v_{2}^{\prime }\left( s\right)
	\right) ,v^{\prime }\left( s\right) \rangle ds+2\mu
	_{2}\int_{0}^{t}\left\vert v^{\prime }\left( 0,s\right) \right\vert ^{2}ds. 
	\notag
	\end{align}
	
	Essentially, our procedure below is similar to the above parts: attempt to
	estimate $\mathcal{K}_{i}$ for $i=1,4$ to derive the uniform boundedness of $%
	\mathcal{W}\left( t\right) $ for which we can use Gronwall's inequality,
	then the proof of uniqueness is self-contained. To handle this, we first
	state the following inequality: for all $r\geq 2$, there exists $\bar{C}%
	_{r}>0$ such that%
	\begin{equation*}
	\left( \Psi _{r}\left( z_{1}\right) -\Psi _{r}\left( z_{2}\right) \right)
	\left( z_{1}-z_{2}\right) \geq \bar{C}_{r}\left\vert z_{1}-z_{2}\right\vert
	^{r},\text{ \ }\forall z_{1},z_{2}\in \mathbb{R}.
	\end{equation*}
	
	It therefore leads to the fact that%
	\begin{align} \tag{3.48}  \label{c48}
	\mathcal{W}\left( t\right) &\geq \left\Vert u^{\prime }\left( t\right)
	\right\Vert ^{2}+\left\Vert v^{\prime }\left( t\right) \right\Vert
	^{2}+\left\Vert u_{x}\left( t\right) \right\Vert ^{2}+\left\Vert v_{x}\left(
	t\right) \right\Vert ^{2}   \\
	&+2\bar{C}_{r_{1}}\lambda _{1}\int_{0}^{t}\left\Vert u^{\prime }\left(
	s\right) \right\Vert _{L^{r_{1}}}^{r_{1}}ds+\bar{C}_{r_{2}}\lambda
	_{2}\int_{0}^{t}\left\Vert v^{\prime }\left( s\right) \right\Vert
	_{L^{r_{2}}}^{r_{2}}ds  \notag \\
	&+2\mu _{1}\int_{0}^{t}\left\vert u^{\prime }\left( 1,s\right) \right\vert
	^{2}ds+2\mu _{2}\int_{0}^{t}\left\vert v^{\prime }\left( 0,s\right)
	\right\vert ^{2}ds.  \notag
	\end{align}
	
	Second, we introduce%
	\begin{equation*}
	M=\max_{i=1,2}\left( \left\Vert u_{ix}\right\Vert _{L^{\infty }\left(
		0,T_{\ast };H^{1}\right) }+\left\Vert v_{ix}\right\Vert _{L^{\infty }\left(
		0,T_{\ast };H^{1}\right) }\right) ,
	\end{equation*}%
	\begin{equation*}
	\mathcal{\tilde{C}}\left( M\right) =\max_{i=1,2}\sup_{\left\vert
		y\right\vert ,\left\vert z\right\vert \leq M}\left( \left\vert \frac{%
		\partial f_{i}}{\partial y}\left( y,z\right) \right\vert +\left\vert \frac{%
		\partial f_{i}}{\partial z}\left( y,z\right) \right\vert \right) ,\quad
	i=1,2.
	\end{equation*}
	Then it is sufficient to estimate $\mathcal{K}_{i}$ for $i=\overline{1,4}$. Indeed,
	we apply the Cauchy-Schwartz inequality to have%
	\begin{align}	\tag{3.49}  \label{c49}
	&\mathcal{K}_{1}+\mathcal{K}_{2} \\ &\leq 2\int_{0}^{t}\left[ \left\Vert
	f_{1}\left( u_{1},v_{1}\right) -f_{1}\left( u_{2},v_{2}\right) \right\Vert
	\left\Vert u^{\prime }\left( s\right) \right\Vert +\left\Vert f_{2}\left(
	u_{1},v_{1}\right) -f_{2}\left( u_{2},v_{2}\right) \right\Vert \left\Vert
	v^{\prime }\left( s\right) \right\Vert \right] ds \nonumber \\ 
	&\leq 2\int_{0}^{t}\left( \left\Vert u\left( s\right)
	\right\Vert +\left\Vert v\left( s\right) \right\Vert \right) \left( \mathcal{%
		\tilde{C}}\left( M\right) \left\Vert u^{\prime }\left( s\right) \right\Vert +%
	\mathcal{\tilde{C}}\left( M\right) \left\Vert v^{\prime }\left( s\right)
	\right\Vert \right) ds \nonumber\\ 
	&\leq 2\int_{0}^{t}\left( \left\Vert u_{x}\left(
	s\right) \right\Vert +\left\Vert v_{x}\left( s\right) \right\Vert \right)
	\left( \mathcal{\tilde{C}}\left( M\right) \left\Vert u^{\prime }\left(
	s\right) \right\Vert +\mathcal{\tilde{C}}\left( M\right) \left\Vert
	v^{\prime }\left( s\right) \right\Vert \right) ds  \nonumber\\& 
	 \leq 4\mathcal{\tilde{C}}\left( M\right) \int_{0}^{t}%
	\mathcal{W}\left( s\right) ds. \nonumber%
	\end{align}
	
	To estimate $\mathcal{K}_{3}$ and $\mathcal{K}_{4}$, we only need to
	consider two cases, $p_{1}=p_{2}=2$ and $p_{1},$ $p_{2}>2$. First, one may easily show that for $p_{1}=p_{2}=2$,
	\begin{align}\tag{3.50}  \label{c50}
	&\mathcal{K}_{3}+\mathcal{K}_{4} \\&=2K_{1}\int_{0}^{t}u\left( 1,s\right)
	u^{\prime }\left( 1,s\right) ds+2K_{2}\int_{0}^{t}v\left( 0,s\right)
	v^{\prime }\left( 0,s\right) ds   \notag\\
	&\leq \frac{K_{1}^{2}}{\mu _{1}}\int_{0}^{t}u^{2}\left( 1,s\right) ds+\mu
	_{1}\int_{0}^{t}\left\vert u^{\prime }\left( 1,s\right) \right\vert ^{2}ds+%
	\frac{K_{2}^{2}}{\mu _{2}}\int_{0}^{t}v^{2}\left( 0,s\right) ds+\mu
	_{2}\int_{0}^{t}\left\vert v^{\prime }\left( 0,s\right) \right\vert ^{2}ds 
	\notag \\
	&\leq \left( \frac{K_{1}^{2}}{\mu _{1}}+\frac{K_{2}^{2}}{\mu _{2}}\right)
	\int_{0}^{t}\mathcal{W}\left( s\right) ds+\frac{1}{2}\mathcal{W}\left(
	t\right) ,  \notag
	\end{align}%
	where we have followed from (\ref{c48}) the inequality%
	\begin{equation*}
	\mathcal{W}\left( t\right) \geq u^{2}\left( 1,t\right) +v^{2}\left(
	0,t\right) +2\left( \mu _{1}\int_{0}^{t}\left\vert u^{\prime
	}\left( 1,s\right) \right\vert ^{2}ds+\mu
	_{2}\int_{0}^{t}\left\vert v^{\prime }\left( 0,s\right)
	\right\vert ^{2}ds\right) .%
	\end{equation*}
	
	For $p_{1},$ $p_{2}>2$, we have%
	\begin{align}	\tag{3.51}  \label{c51}
	\mathcal{K}_{3}+\mathcal{K}_{4}
	& \leq 2\left[ K_{1}\left( p_{1}-1\right)
	M^{p_{1}-2}\int_{0}^{t}\left\vert u\left( 1,s\right) \right\vert \left\vert
	u^{\prime }\left( 1,s\right) \right\vert ds \right. \\
	& \left.
	+K_{2}\left( p_{2}-1\right)
	M^{p_{2}-2}\int_{0}^{t}\left\vert v\left( 0,s\right) \right\vert \left\vert
	v^{\prime }\left( 0,s\right) \right\vert ds\right]  \nonumber \\ 
	&\leq \left[ \frac{K_{1}^{2}}{\mu _{1}}\left(
	p_{1}-1\right) ^{2}M^{2p_{1}-4}+\frac{K_{2}^{2}}{\mu _{2}}\left(
	p_{2}-1\right) ^{2}M^{2p_{2}-4}\right] \int_{0}^{t}\mathcal{W}\left(
	s\right) ds+\frac{1}{2}\mathcal{W}\left( t\right), \nonumber%
	\end{align}%
	where we have recalled the inequality (\ref{c38}).
	
	Combining (\ref{c50}) and (\ref{c51}), we claim that there exists $\eta
	=\eta \left( p_{1},p_{2}\right) >0$ depending on $p_{1},$ $p_{2}$ such that%
	\begin{equation}
	\mathcal{K}_{3}+\mathcal{K}_{4}\leq \eta \left( p_{1},p_{2}\right)
	\int_{0}^{t}\mathcal{W}\left( s\right) ds+\frac{1}{2}\mathcal{W}%
	\left( t\right) .
	\tag{3.52}  \label{c52}
	\end{equation}
	
	Therefore, (\ref{c49}) and (\ref{c52}) together with (\ref{c46}) -- (\ref%
	{c48}) imply that%
	\begin{equation*}
	\mathcal{W}\left( t\right) \leq 2\left[ 4\mathcal{\tilde{C}}\left( M\right)
	+\eta \left( p_{1},p_{2}\right) \right] \int_{0}^{t}\mathcal{W}%
	\left( s\right) ds.%
	\end{equation*}
	Thanks to Gronwall's inequality, we have $\mathcal{W}\left( t\right)
	\equiv 0$ that indicates the uniqueness of solution. Hence, this completes the
	proof of the theorem.
\end{proof}

\begin{remark}
	It is worth noting in the above theorem that the
	existence of a strong solution can be obtained from the regularity of weak
	solutions. In fact, (\ref{c2}) allows us to show that there exists a pair of
	strong solutions $\left( u,v\right) $ to the problem $\left( P\right) $,
	which satisfies%
	\begin{equation}
	\begin{cases}
	\left( u,v\right) \in L^{\infty }\left( 0,T_{\ast };\left( \mathbb{V}%
	_{1}\cap H^{2}\right) \times \left( \mathbb{V}_{2}\cap H^{2}\right) \right)
	\cap C^{0}\left( \left[ 0,T_{\ast }\right] ;\mathbb{V}_{1}\times \mathbb{V}%
	_{2}\right) \\ \qquad \qquad \qquad \qquad \qquad \qquad \qquad \qquad \qquad  \qquad \cap C^{1}\left( \left[ 0,T_{\ast }\right] ;L^{2}\times
	L^{2}\right) , \\ 
	\left( u^{\prime },v^{\prime }\right) \in L^{\infty }\left( 0,T_{\ast };%
	\mathbb{V}_{1}\times \mathbb{V}_{2}\right) \cap C^{0}\left( \left[ 0,T_{\ast
	}\right] ;L^{2}\times L^{2}\right) , \\ 
	\left( u^{\prime \prime },v^{\prime \prime }\right) \in L^{\infty }\left(
	0,T;L^{2}\times L^{2}\right) , \\ 
	\left\vert u^{\prime }\right\vert ^{\frac{r_{1}}{2}-1}u^{\prime },\text{ \ }%
	\left\vert v^{\prime }\right\vert ^{\frac{r_{2}}{2}-1}v^{\prime }\in
	H^{1}\left( Q_{T_{\ast }}\right) , \\ 
	\left\vert u^{\prime }\left( 1,\cdot \right) \right\vert ^{\frac{q_{1}}{2}%
		-1}u^{\prime }\left( 1,\cdot \right) ,\text{ \ }\left\vert v^{\prime }\left(
	0,\cdot \right) \right\vert ^{\frac{q_{2}}{2}-1}v^{\prime }\left( 0,\cdot
	\right) \in H^{1}\left( 0,T_{\ast }\right) .%
	\end{cases}
	\tag{3.53}  \label{c53}
	\end{equation}
\end{remark}

Set the following assumptions:

\textbf{(}$\mathbf{B}_{\mathbf{1}}$\textbf{)} $\left( \tilde{u}_{0},\tilde{u}%
_{1}\right) \in \mathbb{V}_{1}\times L^{2}$ and $\left( \tilde{v}_{0},\tilde{%
	v}_{1}\right) \in \mathbb{V}_{2}\times L^{2}$;

\textbf{(}$\mathbf{B}_{\mathbf{2}}$\textbf{)} $F_{1},F_{2}\in L^{2}\left(
Q_{T}\right) $.

In the following theorem, the existence
and uniqueness of a local weak solution are also obtainable using \textbf{(}$\mathbf{B}_{\mathbf{1}}$\textbf{)}-\textbf{(}$\mathbf{B}_{\mathbf{2}}$\textbf{)} instead of \textbf{(}$\mathbf{A}_{\mathbf{1}}$\textbf{)}-\textbf{(}$\mathbf{A}_{\mathbf{2}}$\textbf{)}. In this scenario, we remark that the regularity of initial data is lower (compared to \textbf{(}$\mathbf{A}_{\mathbf{1}}$\textbf{)}), while the external functions lack of information (compared to \textbf{(}$\mathbf{A}_{\mathbf{2}}$\textbf{)}).

\begin{theorem}\label{thm:3.3}
	Let $q_{1}=q_{2}=2$ and $p_{1},p_{2}\geq 2$ in the problem $(P)$. Assume that \emph{\textbf{(}$\mathbf{A}_{\mathbf{3}}$\textbf{)}} and \emph{\textbf{(}$\mathbf{B}_{\mathbf{1}}$\textbf{)}-\textbf{(}$\mathbf{B}_{\mathbf{2}}$\textbf{)}} hold. Moreover, the initial data obey the
	compatibility conditions (\ref{c1}). Then the problem $(P)$ admits a unique local solution $\left( u,v\right) $ such that
	\begin{equation*}
	\left\{ 
	\begin{tabular}{l}
	$\left( u,v\right) \in C^{0}\left( \left[ 0,T_{\ast }\right] ;\mathbb{V}%
	_{1}\times \mathbb{V}_{2}\right) \cap C^{1}\left( \left[ 0,T_{\ast }\right]
	;L^{2}\times L^{2}\right) ,$ \\ 
	$\left( u^{\prime },v^{\prime }\right) \in L^{r_{1}}\left( Q_{T_{\ast
	}}\right) \times L^{r_{2}}\left( Q_{T_{\ast }}\right) ,u\left( 1,\cdot
	\right) ,\text{ }v\left( 0,\cdot \right) \in H^{1}\left( 0,T_{\ast }\right)
	, $%
	\end{tabular}%
	\right.
	\end{equation*}%
	for $T_{\ast }>0$ sufficiently small.
\end{theorem}

\begin{proof}[Proof of Theorem \ref{thm:3.3}]
	In this proof, we establish several sequences $\left\{ \left(
	u_{0m},u_{1m}\right) \right\} \in C_{0}^{\infty }\left( \overline{\Omega }%
	\right) \times C_{0}^{\infty }\left( \overline{\Omega }\right) ,$ $\left\{
	\left( v_{0m},v_{1m}\right) \right\} \in C_{0}^{\infty }\left( \overline{%
		\Omega }\right) \times C_{0}^{\infty }\left( \overline{\Omega }\right) $
	and $\left\{ \left( F_{1m},F_{2m}\right) \right\} \subset C_{0}^{\infty
	}\left( \overline{Q_{T}}\right) \times C_{0}^{\infty }\left( \overline{Q_{T}}%
	\right) $ satisfying%
	\begin{equation*}
	\left( u_{0m},u_{1m}\right) \rightarrow \left( \tilde{u}_{0},\tilde{u}%
	_{1}\right) \;\mbox{strong in}\;\mathbb{V}_{1}\times L^{2},
	\end{equation*}%
	\begin{equation*}
	\left( v_{0m},v_{1m}\right) \rightarrow \left( \tilde{v}_{0},\tilde{v}%
	_{1}\right) \;\mbox{strong in}\;\mathbb{V}_{2}\times L^{2},
	\end{equation*}%
	\begin{equation*}
	\left( F_{1m},F_{2m}\right) \rightarrow \left( F_{1},F_{2}\right) \;%
	\mbox{strong in}\;[L^{2}\left( Q_{T}\right)]^2
	.
	\end{equation*}
	
	Note that the sequences $\left\{ \left( u_{0m},u_{1m}\right) \right\} $ and $%
	\left\{ \left( v_{0m},v_{1m}\right) \right\} $, as a result, satisfy
	themselves the compatibility relation for all $m\in \mathbb{N}$. So what we
	can deduce next is the existence of a pair of unique functions $\left(
	u_{m},v_{m}\right) $ for each $m$ that makes the conditions in the
	aforementioned theorem self-propelling. Thus, one easily verifies that such
	functions $\left( u_{m},v_{m}\right) $ for each $m$ satisfy the variational
	problem (\ref{b1})-(\ref{b2}), i.e.%
	\begin{equation}
	\begin{cases}
	\langle u_{m}^{\prime \prime }\left( t\right) ,\phi \rangle +\langle
	u_{mx}\left( t\right) ,\phi _{x}\rangle +\lambda _{1}\langle \Psi
	_{r_{1}}\left( u_{m}^{\prime }\left( t\right) \right) ,\phi \rangle +\mu
	_{1}u_{m}^{\prime }\left( 1,t\right) \phi \left( 1\right) \\ 
	\text{ \ \ \ \ \ \ \ \ \ \ \ \ \ }=K_{1}\Psi _{p_{1}}\left( u_{m}\left(
	1,t\right) \right) \phi \left( 1\right) +\langle f_{1}\left(
	u_{m},v_{m}\right) ,\phi \rangle +\langle F_{1m}\left( t\right) ,\phi
	\rangle , \\ 
	\langle v_{m}^{\prime \prime }\left( t\right) ,\tilde{\phi}\rangle +\langle
	v_{mx}\left( t\right) ,\tilde{\phi}_{x}\rangle +\lambda _{2}\langle \Psi
	_{r_{2}}\left( v_{m}^{\prime }\left( t\right) \right) ,\tilde{\phi}\rangle
	+\mu _{2}v_{m}^{\prime }\left( 0,t\right) \tilde{\phi}\left( 0\right) \\ 
	\text{ \ \ \ \ \ \ \ \ \ \ \ \ \ }=K_{2}\Psi _{p_{2}}\left( v_{m}\left(
	0,t\right) \right) \tilde{\phi}\left( 0\right) +\langle f_{2}\left(
	u_{m},v_{m}\right) ,\tilde{\phi}\rangle +\langle F_{2m}\left( t\right) ,%
	\tilde{\phi}\rangle ,%
	\end{cases}
	\tag{3.54}  \label{c54}
	\end{equation}%
	for all $(\phi ,\tilde{\phi})\in \mathbb{V}_{1}\times \mathbb{V}_{2}$,
	together with the initial conditions%
	\begin{equation}
	\left( u_{m}\left( 0\right) ,u_{m}^{\prime }\left( 0\right) \right) =\left(
	u_{0m},u_{1m}\right) ,\text{ \ }\left( v_{m}\left( 0\right) ,v_{m}^{\prime
	}\left( 0\right) \right) =\left( v_{0m},v_{1m}\right) .  \tag{3.55}
	\label{c55}
	\end{equation}
	
	Moreover, the smoothness of $\left( u_{m},v_{m}\right) $ on the interval $%
	[0,T_{\ast }]$ is said by (\ref{c53}) and we recall below the uniform
	boundedness (independent of $m$) of $\mathcal{S}_{m}\left( t\right) $ on $%
	[0,T_{\ast }]$ in (\ref{c5}) due to the same arguments derived above.%
	\begin{align}\tag{3.56}  \label{c56}
	 \mathcal{S}_{m}\left( t\right) &=\left\Vert u_{m}^{\prime }\left( t\right)
	\right\Vert ^{2}+\left\Vert v_{m}^{\prime }\left( t\right) \right\Vert
	^{2}+\left\Vert u_{mx}\left( t\right) \right\Vert ^{2}+\left\Vert
	v_{mx}\left( t\right) \right\Vert ^{2}+2\lambda _{1}\int_{0}^{t}\left\Vert
	u_{m}^{\prime }\left( s\right) \right\Vert _{L^{r_{1}}}^{r_{1}}ds 
	  \\
	&+2\lambda _{2}\int_{0}^{t}\left\Vert v_{m}^{\prime }\left( s\right)
	\right\Vert _{L^{r_{2}}}^{r_{2}}ds+2\mu _{1}\int_{0}^{t}\left\vert
	u_{m}^{\prime }\left( 1,s\right) \right\vert ^{2}ds+2\mu
	_{2}\int_{0}^{t}\left\vert v_{m}^{\prime }\left( 0,s\right) \right\vert
	^{2}ds\leq C_{T},  \notag 
	\end{align}%
	where $C_{T}$ denotes a positive constant independent of $m$ and $t$ and $%
	t$ moves along the interval $\left[ 0,T_{\ast }\right] $.
	
	Define $U_{m,k}:=u_{m}-u_{k}$ and $V_{m,k}:=v_{m}-v_{k}$, then these
	quantities satisfy%
	\begin{equation}
	\begin{cases}
	\langle U_{m,k}^{\prime \prime }\left( t\right) ,\phi \rangle +\langle
	\nabla U_{m,k}\left( t\right) ,\phi _{x}\rangle +\lambda _{1}\langle \Psi
	_{r_{1}}\left( u_{m}^{\prime }\left( t\right) \right) -\Psi _{r_{1}}\left(
	u_{k}^{\prime }\left( t\right) \right) ,\phi \rangle  \\
	+\mu
	_{1}U_{m,k}^{\prime }\left( 1,t\right) \phi \left( 1\right) 
	=K_{1}\left[ \Psi _{p_{1}}\left( u_{m}\left( 1,t\right) \right) -\Psi
	_{p_{1}}\left( u_{k}\left( 1,t\right) \right) \right] \phi \left( 1\right)\\
	+\langle f_{1}\left( u_{m},v_{m}\right) -f_{1}\left( u_{k},v_{k}\right)
	,\phi \rangle +\langle F_{1m}\left( t\right) -F_{1k}\left( t\right) ,\phi
	\rangle , \\ 
	\langle V_{m,k}^{\prime \prime }\left( t\right) ,\tilde{\phi}\rangle
	+\langle \nabla V_{m,k}\left( t\right) ,\tilde{\phi}_{x}\rangle +\lambda
	_{2}\langle \Psi _{r_{2}}\left( v_{m}^{\prime }\left( t\right) \right) -\Psi
	_{r_{2}}\left( v_{k}^{\prime }\left( t\right) \right) ,\tilde{\phi}\rangle
	 \\ +\mu _{2}V_{m,k}^{\prime }\left( 0,t\right) \tilde{\phi}\left( 0\right)
	=K_{2}\left[ \Psi _{p_{2}}\left( v_{m}\left( 0,t\right) \right) -\Psi
	_{p_{2}}\left( v_{k}\left( 0,t\right) \right) \right] \tilde{\phi}\left(
	0\right) \\+\langle f_{2}\left( u_{m},v_{m}\right) -f_{2}\left(
	u_{k},v_{k}\right) ,\tilde{\phi}\rangle +\langle F_{2m}\left( t\right)
	-F_{2k}\left( t\right) ,\tilde{\phi}\rangle ,%
	\end{cases}
	\tag{3.57}  \label{c57}
	\end{equation}%
	for all $(\phi ,\tilde{\phi})\in \mathbb{V}_{1}\times \mathbb{V}_{2}$ and
	the initial conditions are%
	\begin{align}	
	\left( U_{m,k}\left( 0\right) ,U_{m,k}^{\prime }\left( 0\right) \right)
	=\left( u_{0m}-u_{0k},u_{1m}-u_{1k}\right), \nonumber \\ \left( V_{m,k}\left(
	0\right) ,V_{m,k}^{\prime }\left( 0\right) \right) =\left(
	v_{0m}-v_{0k},v_{1m}-v_{1k}\right) .\tag{3.58}  \label{c58}
	\end{align}
	
	It is obvious to obtain the fact that%
	\begin{align}\tag{3.59}  \label{c59}
	\mathcal{S}_{m,k}\left( t\right) &=\mathcal{S}_{m,k}\left( 0\right)  
	+2\int_{0}^{t}\left( \langle f_{1}\left(
	u_{m},v_{m}\right) -f_{1}\left( u_{k},v_{k}\right) ,U_{m,k}^{\prime }\left(
	s\right) \rangle \right. \\&
	\left. +\langle f_{2}\left( u_{m},v_{m}\right) -f_{2}\left(
	u_{k},v_{k}\right) ,V_{m,k}^{\prime }\left( s\right) \rangle \right) ds \nonumber\\ 
	&+2\int_{0}^{t}\left[ \langle F_{1m}\left(
	s\right) -F_{1k}\left( s\right) ,U_{m,k}^{\prime }\left( s\right) \rangle
	+\langle F_{2m}\left( s\right) -F_{2k}\left( s\right) ,V_{m,k}^{\prime
	}\left( s\right) \rangle \right] ds \nonumber\\ 
	&+2K_{1}\int_{0}^{t}\left[ \Psi _{p_{1}}\left(
	u_{m}\left( 1,s\right) \right) -\Psi _{p_{1}}\left( u_{k}\left( 1,s\right)
	\right) \right] U_{m,k}^{\prime }\left( 1,s\right) ds \nonumber\\ 
	&+2K_{2}\int_{0}^{t}\left[ \Psi _{p_{2}}\left(
	v_{m}\left( 0,s\right) \right) -\Psi _{p_{2}}\left( v_{k}\left( 0,s\right)
	\right) \right] V_{m,k}^{\prime }\left( 0,s\right) ds \nonumber\\ 
	& =\mathcal{S}_{m,k}\left( 0\right) +\mathcal{Z}_{1}+%
	\mathcal{Z}_{2}+\mathcal{Z}_{3}+\mathcal{Z}_{4},\nonumber%
	\end{align}%
	where we have denoted by
	\begin{align}	\tag{3.60}  \label{c60}
	\mathcal{S}_{m,k}\left( t\right) &=\left\Vert U_{m,k}^{\prime }\left(
	t\right) \right\Vert ^{2}+\left\Vert V_{m,k}^{\prime }\left( t\right)
	\right\Vert ^{2}+\left\Vert \nabla U_{m,k}\left( t\right) \right\Vert
	^{2}+\left\Vert \nabla V_{m,k}\left( t\right) \right\Vert ^{2} \\ 
	&+2\lambda _{1}\int_{0}^{t}\langle \Psi
	_{r_{1}}\left( u_{m}^{\prime }\left( s\right) \right) -\Psi _{r_{1}}\left(
	u_{k}^{\prime }\left( s\right) \right) ,U_{m,k}^{\prime }\left( s\right)
	\rangle ds \nonumber \\ 
	&+2\lambda _{2}\int_{0}^{t}\langle \Psi
	_{r_{2}}\left( v_{m}^{\prime }\left( s\right) \right) -\Psi _{r_{2}}\left(
	v_{k}^{\prime }\left( s\right) \right) ,V_{m,k}^{\prime }\left( s\right)
	\rangle ds \nonumber \\ &+2\mu _{1}\int_{0}^{t}\left\vert
	U_{m,k}^{\prime }\left( 1,s\right) \right\vert ^{2}ds+2\mu
	_{2}\int_{0}^{t}\left\vert V_{m,k}^{\prime }\left( 0,s\right)
	\right\vert ^{2}ds,\nonumber%
	\end{align}%
	\begin{equation}
	\begin{tabular}{l}
	$\mathcal{S}_{m,k}\left( 0\right) =\left\Vert u_{1m}-u_{1k}\right\Vert
	^{2}+\left\Vert v_{1m}-v_{1k}\right\Vert ^{2}+\left\Vert
	u_{0mx}-u_{0kx}\right\Vert ^{2}+\left\Vert v_{0mx}-v_{0kx}\right\Vert ^{2}.$%
	\end{tabular}
	\tag{3.61}  \label{c61}
	\end{equation}
	
	The above calculations are done by a valid replacement of the test functions $%
	\phi $ and $\tilde{\phi}$ by $U_{m,k}^{\prime }$ and $V_{m,k}^{\prime }$,
	respectively, in (\ref{c57}) and then integrating with respect to $t$. By
	the same strategy and using (\ref{c38}), (\ref{c56}) and (\ref{c60}), we
	can estimate the terms on the right-hand side of (\ref{c59}) and obtain%
	\begin{equation}
	\mathcal{S}_{m,k}\left( t\right) \leq \mathcal{R}_{m,k}+2\left( 1+\eta
	\left( p_{1},p_{2}\right) +8\mathcal{R}_{T}\right) \int_{0}^{t}%
	\mathcal{S}_{m,k}\left( s\right) ds,%
	\tag{3.62}  \label{c62}
	\end{equation}%
	where the involved terms are appropriately defined as follows:
	\begin{align*}
	\mathcal{R}_{T}&:=\max_{i=1,2}\sup_{\left\vert y\right\vert ,\left\vert
		z\right\vert \leq \sqrt{C_{T}}}\left( \left\vert \frac{\partial f_{i}}{%
		\partial y}\left( y,z\right) \right\vert +\left\vert \frac{\partial f_{i}}{%
		\partial z}\left( y,z\right) \right\vert \right) , \\ 
	\eta \left( p_{1},p_{2}\right) &:=\frac{K_{1}^{2}}{\mu _{1}}\left(
	p_{1}-1\right) ^{2}C_{T}^{p_{1}-2}+\frac{K_{2}^{2}}{\mu _{2}}\left(
	p_{2}-1\right) ^{2}C_{T}^{p_{2}-2}, \\ 
	\mathcal{R}_{m,k}&:=2\mathcal{S}_{m,k}\left( 0\right) +2\left\Vert
	F_{1m}-F_{1k}\right\Vert _{L^{2}\left( Q_{T}\right) }^{2}+2\left\Vert
	F_{2m}-F_{2k}\right\Vert _{L^{2}\left( Q_{T}\right) }^{2}.%
	\end{align*}
	
	Here we remark that $\mathcal{R}_{m,k}$ approaches zero as $m$ and $k$ tend
	to infinity. By the aid of Gronwall's inequality, it follows from (\ref%
	{c62}) that for all $t\in \left[ 0,T_{\ast }\right] $%
	\begin{equation}
	\mathcal{S}_{m,k}\left( t\right) \leq \mathcal{R}_{m,k}\mbox{exp}\left(
	2T\left( 1+\eta \left( p_{1},p_{2}\right) +8\mathcal{R}_{T}\right) \right) .
	\tag{3.63}  \label{c63}
	\end{equation}
	
	Thus, one can show that the right-hand side of (\ref{c63}) goes to zero as
	the indexes $m$ and $k$ tend to infinity by the direct argument concerning
	convergences of $\left\{ \left( u_{0m},u_{1m}\right) \right\} $ and $\left\{
	\left( v_{0m},v_{1m}\right) \right\} $. As by-product, it gives us the
	following results%
	\begin{equation}
	\left( u_{m},v_{m}\right) \rightarrow \left( u,v\right) \;\mbox{strong in}%
	\;C\left( \left[ 0,T_{\ast }\right] ;\mathbb{V}_{1}\times \mathbb{V}%
	_{2}\right) \cap C^{1}\left( \left[ 0,T_{\ast }\right] ;L^{2}\times
	L^{2}\right) ,  \tag{3.64}  \label{c64}
	\end{equation}%
	\begin{equation}
	\left( u_{m}^{\prime },v_{m}^{\prime }\right) \rightarrow \left( u^{\prime
	},v^{\prime }\right) \;\mbox{strong in}\;L^{r_{1}}\left( Q_{T_{\ast
	}}\right) \times L^{r_{2}}\left( Q_{T_{\ast }}\right) ,  \tag{3.65}
	\end{equation}%
	\begin{equation}
	\left( u_{m}\left( 1,\cdot \right) ,v_{m}\left( 0,\cdot \right) \right)
	\rightarrow \left( u\left( 1,\cdot \right) ,v\left( 0,\cdot \right) \right)
	\;\mbox{strong in}\;[H^{1}\left( 0,T_{\ast }\right)]^2.  \tag{3.66}  \label{c66}
	\end{equation}
	
	An important point should be mentioned here is that by (\ref{c56}) we can
	extract a subsequence of $\left\{ \left( u_{m},v_{m}\right) \right\} $
	(still relabel with the old index $m$) which reads as%
	\begin{equation}
	\left( u_{m},v_{m}\right) \rightarrow \left( u,v\right) \;\mbox{weak-*
		in}\;L^{\infty }\left( 0,T_{\ast };\mathbb{V}_{1}\times \mathbb{V}%
	_{2}\right) ,  \tag{3.67}  \label{c67}
	\end{equation}%
	\begin{equation}
	\left( u_{m}^{\prime },v_{m}^{\prime }\right) \rightarrow \left( u^{\prime
	},v^{\prime }\right) \;\mbox{weak-* in}\;L^{\infty }\left( 0,T_{\ast
	};L^{2}\times L^{2}\right) ,  \tag{3.68}  \label{c68}
	\end{equation}%
	and inheriting from (\ref{c64})-(\ref{c66}), one deduces%
	\begin{equation}
	\left( f_{1}\left( u_{m},v_{m}\right) ,f_{2}\left( u_{m},v_{m}\right)
	\right) \rightarrow \left( f_{1}\left( u,v\right) ,f_{2}\left( u,v\right)
	\right) \;\mbox{strong in}\;[L^{2}\left( Q_{T_{\ast }}\right)]^2,  \tag{3.69}  \label{c69}
	\end{equation}%
	\begin{equation}
	\left( \Psi _{r_{1}}\left( u_{m}^{\prime }\right) ,\Psi _{r_{2}}\left(
	v_{m}^{\prime }\right) \right) \rightarrow \left( \Psi _{r_{1}}\left(
	u^{\prime }\right) ,\Psi _{r_{2}}\left( v^{\prime }\right) \right) \;%
	\mbox{strong in}\;[L^{2}\left( Q_{T_{\ast }}\right)]^2.  \tag{3.70}  \label{c70}
	\end{equation}
	
	Henceforward, it is advantageous to state our limit processing now. In fact, passing to the
	limit in (\ref{c54}) associated with (\ref{c55}) and evidenced by (\ref{c64}%
	)-(\ref{c70}), we obtain a couple of functions $\left( u,v\right) $
	satisfying the variational problem%
	\begin{equation*}
	\begin{cases}
	\frac{d}{dt}\left\langle u^{\prime }\left( t\right) ,\phi \right\rangle
	+\left\langle u_{x}\left( t\right) ,\phi _{x}\right\rangle +\lambda
	_{1}\left\langle \Psi _{r_{1}}\left( u^{\prime }\left( t\right) \right)
	,\phi \right\rangle +\mu _{1}\Psi _{q_{1}}\left( u^{\prime }\left(
	1,t\right) \right) \phi \left( 1\right) \\ 
	\text{ \ \ \ \ \ \ \ \ \ \ \ \ \ \ \ \ \ \ \ \ \ \ \ \ \ \ \ \ \ \ \ \ }%
	=K_{1}\Psi _{p_{1}}\left( u\left( 1,t\right) \right) \phi \left( 1\right)
	+\left\langle f_{1}\left( u,v\right) ,\phi \right\rangle +\left\langle
	F_{1}\left( t\right) ,\phi \right\rangle , \\ 
	\frac{d}{dt}\left\langle v^{\prime }\left( t\right) ,\tilde{\phi}%
	\right\rangle +\left\langle v_{x}\left( t\right) ,\tilde{\phi}%
	_{x}\right\rangle +\lambda _{2}\left\langle \Psi _{r_{2}}\left( v^{\prime
	}\left( t\right) \right) ,\tilde{\phi}\right\rangle +\mu _{2}\Psi
	_{q_{2}}\left( v^{\prime }\left( 0,t\right) \right) \tilde{\phi}\left(
	0\right) \\ 
	\text{ \ \ \ \ \ \ \ \ \ \ \ \ \ \ \ \ \ \ \ \ \ \ \ \ \ \ \ \ \ \ \ \ }%
	=K_{2}\Psi _{p_{2}}\left( v\left( 0,t\right) \right) \tilde{\phi}\left(
	0\right) +\left\langle f_{2}\left( u,v\right) ,\tilde{\phi}\right\rangle
	+\left\langle F_{2}\left( t\right) ,\tilde{\phi}\right\rangle ,%
	\end{cases}%
	\end{equation*}%
	for all $(\phi ,\tilde{\phi})\in \mathbb{V}_{1}\times \mathbb{V}_{2}$, and
	endowed with (\ref{b2}). This also indicates the existence of a local
	solution. The uniqueness of such a weak solution is directly obtained by
	using the regularization procedure investigated by Lions (see e.g. 
	\cite{12}). Hence, we end up with the proof of the theorem.
\end{proof}

\begin{remark}
	If one has $N=\frac{1}{2}\max \left\{ 2;\alpha ;\beta ;%
	\frac{q_{1}\left( p_{1}-1\right) }{q_{1}-1};\frac{q_{2}\left( p_{2}-1\right) 
	}{q_{2}-1}\right\} \leq 1$ (cf. Remark \ref{rem:3.2}) and considers the assumptions 
	\textbf{(}$\mathbf{B}_{\mathbf{1}}$\textbf{)} and \textbf{(}$\mathbf{B}_{2}$%
	\textbf{)}, the integral $\mathcal{S}_{m}\left( t\right) $ (which has been
	bounded by a constant $C_{T}$; see in (\ref{c10})) can be estimated
	globally (and uniformly as well) in time, i.e.%
	\begin{equation*}
	\begin{tabular}{l}
	$\mathcal{S}_{m}\left( t\right) \leq C_{T},$ $\forall m\in \mathbb{N},\text{ 
	}\forall t\in \left[ 0,T\right] ,\text{ }\forall T>0.$%
	\end{tabular}%
	\end{equation*}
	Consequently, all arguments used in the proof of Theorem \ref{thm:3.3} are applicable
	to prove that there exists a global weak solution $\left( u,v\right) $ to
	the problem $\left( P\right) $ and satisfying
	\begin{align}\tag{3.71} \label{3.71}
	&\left( u,v\right) \in L^{\infty }\left( 0,T;\mathbb{V}_{1}\times \mathbb{V}%
	_{2}\right), \left( u^{\prime },v^{\prime }\right) \in L^{\infty }\left(
	0,T;L^{2}\times L^{2}\right), \\&u\left( 1,\cdot \right),v\left(
	0,\cdot \right) \in H^{1}\left( 0,T\right).\nonumber 
	\end{align}
	Nevertheless, it should be noticed that the aforementioned case of $N$ does not imply the
	weak solution belongs to $C\left( \left[ 0,T\right] ;\mathbb{V}%
	_{1}\times \mathbb{V}_{2}\right) \cap C^{1}\left( \left[ 0,T\right]
	;L^{2}\times L^{2}\right) $ and one cannot also say any further statement
	for the uniqueness in this case.%
\end{remark}

\section{Finite time blow-up}
In principle, the blow-up phenomenon of a solution to a
time-dependent equation is devoted to the study of maximal time domain for
which it is defined by a finite length. At the endpoint of that interval,
the solution behaves in such a way that either it goes to infinity in some
specific senses, or it stops being smooth, and so forth. Our main objective
here is to show that if the initial energy is negative, then every weak
solution of $\left( P\right) $ blows up in finite time. The result here
draws from ideas in the treatment of a single wave equation \cite{6a,11,15} and
also, for example, the recent results for systems in \cite{6,13}, but our
proofs have to be radically altered. In addition, the technique we use here is an adaptation of an argument postulated in \cite{15a} to treat the boundary damping terms.

Let us first make a brief note concerning the so-called total energy. From
the mathematical point of view, the energy method plays a vital role in the
study of partial differential equations and the energy is usually used to
derive the well-posedness of such equations. For a very fundamental scalar
wave equation which was originally discovered by d'Alembert, one may easily
find the energy integral computed by%
\begin{equation*}
\frac{1}{2}\int_{\Omega }\left\vert u_{t}\right\vert ^{2}dx+\frac{1}{2}%
\int_{\Omega }\left\vert u_{x}\right\vert ^{2}dx,
\end{equation*}%
where we basically multiply the equation by $u_{t}$ and integrate over $%
\Omega $ and use the divergence theorem.

In the energy equation, the first term is typically called as the kinetic
energy, while the second term is referred to as the potential energy. If one
imposes mixed boundary conditions, says $au+bu_{x}=0$ on the boundary, where 
$a$ and $b$ have the same sign, the energy can only decrease. If they have
opposite sign, the problem is unstable in the sense that the energy will
increase. In this section, the total energy are computed in a common way and
then it would be a decreasing function along the trajectories, starting from
a negative initial value.

From here on, we aim to consider the problem $\left( P\right) $ in a
specific case: linear damping $r_{i}=2$ with $F_{i}=0$, and $q_{i}=2,$ $%
p_{i}>2,$ $K_{i}>0,$ $\lambda _{i}>0,$ $\mu _{i}>0$ for $i=1,2$. 
In this case, the total quadratic-type energy $E(t)$ associated with the solution $(u,v)$ is defined by
\begin{align*}
E\left( t\right) &:=\frac{1}{2}\left( \left\Vert u^{\prime }\left( t\right)
\right\Vert ^{2}+\left\Vert v^{\prime }\left( t\right) \right\Vert
^{2}+\left\Vert u_{x}\left( t\right) \right\Vert ^{2}+\left\Vert v_{x}\left(
t\right) \right\Vert ^{2}\right)  \\ &
 -\left( \frac{K_{1}}{p_{1}}\left\vert u\left( 1,t\right)
\right\vert ^{p_{1}}+\frac{K_{2}}{p_{2}}\left\vert v\left( 0,t\right)
\right\vert ^{p_{2}}\right) -\int_{0}^{1}\mathcal{F}\left( u\left(
x,t\right) ,v\left( x,t\right) \right) dx.%
\nonumber 
\end{align*}

To prove the solution in this case blows up in finite time, we consider the following assumptions on the interior sources and on the initial energy:

\textbf{(}$\mathbf{A}_{3}^{\prime }$\textbf{)} there exists a $C^{2}$-function $\mathcal{F}:%
\mathbb{R}^{2}\rightarrow \mathbb{R}$  such that%
\begin{equation*}
\frac{\partial \mathcal{F}}{\partial u}\left( u,v\right) =f_{1}\left(
u,v\right) ,\text{ \ }\frac{\partial \mathcal{F}}{\partial v}\left(
u,v\right) =f_{2}\left( u,v\right) ,
\end{equation*}%
and there also exists the constants $\alpha ,$ $\beta >2;$ $d_{1},$ $d_{2},$ 
$\bar{d}_{1},$ $\bar{d}_{2}>0$ such that%
\begin{equation}
d_{1}\mathcal{F}\left( u,v\right) \leq uf_{1}\left( u,v\right) +vf_{2}\left(
u,v\right) \leq d_{2}\mathcal{F}\left( u,v\right) \quad \mbox{for all}%
\;\left( u,v\right) \in \mathbb{R}^{2},  \tag{4.1}  \label{d1}
\end{equation}%
\begin{equation}
\bar{d}_{1}\left( \left\vert u\right\vert ^{\alpha }+\left\vert v\right\vert
^{\beta }\right) \leq \mathcal{F}\left( u,v\right) \leq \bar{d}_{2}\left(
\left\vert u\right\vert ^{\alpha }+\left\vert v\right\vert ^{\beta }\right)
\quad \mbox{for all}\;\left( u,v\right) \in \mathbb{R}^{2}.  \tag{4.2}
\label{d2}
\end{equation}

\textbf{(}$\mathbf{A}_{4}$\textbf{)} define $H\left( t\right) :=-E\left( t\right) $ and assume that
\begin{align}\tag{4.3}
-H\left( 0\right) &=\frac{1}{2}\left( \left\Vert \tilde{u}_{1}\right\Vert
^{2}+\left\Vert \tilde{v}_{1}\right\Vert ^{2}+\left\Vert \tilde{u}%
_{0x}\right\Vert ^{2}+\left\Vert \tilde{v}_{0x}\right\Vert ^{2}\right)  \\ 
&-\left( \frac{K_{1}}{p_{1}}\left\vert \tilde{u}%
_{0}\left( 1\right) \right\vert ^{p_{1}}+\frac{K_{2}}{p_{2}}\left\vert 
\tilde{v}_{0}\left( 0\right) \right\vert ^{p_{2}}\right)
-\int_{0}^{1}\mathcal{F}\left( \tilde{u}_{0}\left( x\right) ,%
\tilde{v}_{0}\left( x\right) \right) dx<0.\nonumber
\end{align}

It is worth noting that the time derivative of $H$ is non-negative
along the trajectories in a local time, namely%
\begin{equation}
\begin{tabular}{l}
$H^{\prime }\left( t\right) =\lambda _{1}\left\Vert u^{\prime }\left(
t\right) \right\Vert ^{2}+\lambda _{2}\left\Vert v^{\prime }\left( t\right)
\right\Vert ^{2}+\mu _{1}\left\vert u^{\prime }\left( 1,t\right)
\right\vert^{2} +\mu _{2}\left\vert v^{\prime }\left( 0,t\right) \right\vert
^{2}\geq 0,$%
\end{tabular}
\tag{4.4}  \label{d4}
\end{equation}%
for all $t\in \left[ 0,T_{\ast }\right) $ where we have multiplied (\ref%
{main1}) by $\left( u^{\prime }\left( x,t\right) ,v^{\prime }\left(
x,t\right) \right) $ and integrated the resulting equation over $\Omega $.
Together with the fact that $H\left( 0\right) >0$, we have
\begin{equation}
\begin{tabular}{l}
$0<H\left( 0\right) \leq H\left( t\right) ,\text{ }\forall t\in \left[
0,T_{\ast }\right) .$%
\end{tabular}
\tag{4.5}  \label{d5}
\end{equation}

Observe the final term in the total energy $E(t)$, it follows from (\ref%
{d2}) that%
\begin{equation}
\bar{d}_{1}\left( \left\Vert u\left( t\right) \right\Vert _{L^{\alpha
}}^{\alpha }+\left\Vert v(t)\right\Vert _{L^{\beta }}^{\beta }\right) \leq
\int_{0}^{1}\mathcal{F}\left( u\left( x,t\right) ,v\left(
x,t\right) \right) dx\leq \bar{d}_{2}\left( \left\Vert u\left( t\right)
\right\Vert _{L^{\alpha }}^{\alpha }+\left\Vert v(t)\right\Vert _{L^{\beta
}}^{\beta }\right) .%
\tag{4.6}  \label{d6}
\end{equation}

Combining (\ref{d5}) and (\ref{d6}) gives%
\begin{align}\tag{4.7}  \label{d7}
0 &<H\left( t\right) \leq \frac{K_{1}}{p_{1}}\left\vert u\left( 1,t\right)
\right\vert ^{p_{1}}+\frac{K_{2}}{p_{2}}\left\vert v\left( 0,t\right)
\right\vert ^{p_{2}}+\bar{d}_{2}\left( \left\Vert u\left( t\right)
\right\Vert _{L^{\alpha }}^{\alpha }+\left\Vert v(t)\right\Vert _{L^{\beta
}}^{\beta }\right)   \\
& \leq \bar{D}_{2}\left( \left\vert u\left( 1,t\right) \right\vert
^{p_{1}}+\left\vert v\left( 0,t\right) \right\vert ^{p_{2}}+\left\Vert
u\left( t\right) \right\Vert _{L^{\alpha }}^{\alpha }+\left\Vert
v(t)\right\Vert _{L^{\beta }}^{\beta }\right) ,  \notag 
\end{align}%
for all $t\in \left[ 0,T_{\ast }\right) $ and $\bar{D}_{2}=\max \left\{ 
\frac{K_{1}}{p_{1}},\frac{K_{2}}{p_{2}},\bar{d}_{2}\right\} $. It also
allows us to derive that for all $t\in \left[ 0,T_{\ast }\right) $%
\begin{align}\tag{4.8}  \label{d8}
& \frac{1}{2}\left( \left\Vert u_{x}\left( t\right) \right\Vert
^{2}+\left\Vert v_{x}\left( t\right) \right\Vert ^{2}\right) \\&\leq \frac{%
	K_{1}}{p_{1}}\left\vert u\left( 1,t\right) \right\vert ^{p_{1}}+\frac{K_{2}}{%
	p_{2}}\left\vert v\left( 0,t\right) \right\vert
^{p_{2}}+\int\nolimits_{0}^{1}\mathcal{F}\left( u\left( x,t\right) ,v\left(
x,t\right) \right) dx  \nonumber  \\
&\leq \bar{D}_{2}\left( \left\vert u\left( 1,t\right) \right\vert
^{p_{1}}+\left\vert v\left( 0,t\right) \right\vert ^{p_{2}}+\left\Vert
u\left( t\right) \right\Vert _{L^{\alpha }}^{\alpha }+\left\Vert
v(t)\right\Vert _{L^{\beta }}^{\beta }\right) .  \notag 
\end{align}

Now we construct the following functional%
\begin{equation}
\begin{tabular}{l}
$L\left( t\right) :=H^{1-\xi }\left( t\right) +\varepsilon \psi \left(
t\right) ,$%
\end{tabular}
\tag{4.9}  \label{d9}
\end{equation}%
where we have defined%
\begin{align}\tag{4.10}  \label{d10}
\psi \left( t\right) &:=\langle u\left( t\right) ,u^{\prime }\left( t\right)
\rangle +\langle v\left( t\right) ,v^{\prime }\left( t\right) \rangle \\&+\frac{%
	\lambda _{1}}{2}\left\Vert u\left( t\right) \right\Vert ^{2}+\frac{\lambda
	_{2}}{2}\left\Vert v\left( t\right) \right\Vert ^{2}+\frac{\mu _{1}}{2}%
u^{2}\left( 1,t\right) +\frac{\mu _{2}}{2}v^{2}\left( 0,t\right) ,
\nonumber 
\end{align}%
for $\varepsilon >0$ sufficiently small and $\xi \in \left( 0,\min \left\{ \frac{%
	\alpha -2}{2\alpha },\frac{\beta -2}{2\beta }\right\} \right] \subset \left(
0,\frac{1}{2}\right) $.

To show that one can choose $\varepsilon >0$ small enough such that $L\left(
t\right) $ is non-decreasing for all $t\in \left[ 0,T_{\ast }\right) $,
namely%
\begin{equation}
\begin{tabular}{l}
$L\left( t\right) \geq L\left( 0\right) >0,\text{ \ }\forall t\in \left[
0,T_{\ast }\right) ,$%
\end{tabular}
\tag{4.11}  \label{d11}
\end{equation}%
a very clear way is to consider the derivative of such a function. In fact,
let us prove the following lemma.

\begin{lemma}\label{lem:4.1}
	There exists a positive constant $\gamma $ such that
	\begin{align}\tag{4.12}  \label{d12}
	L^{\prime }\left( t\right) &\geq \gamma \left[ H\left( t\right) +\left\Vert
	u^{\prime }\left( t\right) \right\Vert ^{2}+\left\Vert v^{\prime }\left(
	t\right) \right\Vert ^{2}+\left\Vert u_{x}\left( t\right) \right\Vert
	^{2}+\left\Vert v_{x}\left( t\right) \right\Vert ^{2} \right. \\ & \left.+\left\Vert u\left(
	t\right) \right\Vert _{L^{\alpha }}^{\alpha }+\left\Vert v(t)\right\Vert
	_{L^{\beta }}^{\beta }+\left\vert u\left( 1,t\right) \right\vert
	^{p_{1}}+\left\vert v\left( 0,t\right) \right\vert ^{p_{2}}\right] .	\nonumber
	\end{align}
\end{lemma}

\begin{proof}[Proof of Lemma \ref{lem:4.1}]
	Multiplying (\ref{main1}) by $\left( u\left(
	x,t\right) ,v\left( x,t\right) \right) $ and then integrating the resulting
	equation over $\Omega $, the derivative of $\psi \left( t\right) $ can be
	defined as follows:
	\begin{align*}
		\psi ^{\prime }\left( t\right) &=\left\Vert u^{\prime }\left( t\right)
		\right\Vert ^{2}+\left\Vert v^{\prime }\left( t\right) \right\Vert
		^{2}-\left( \left\Vert u_{x}\left( t\right) \right\Vert ^{2}+\left\Vert
		v_{x}\left( t\right) \right\Vert ^{2}\right) +K_{1}\left\vert u\left(
		1,t\right) \right\vert ^{p_{1}}\\
		&+K_{2}\left\vert v\left( 0,t\right)
		\right\vert ^{p_{2}} +\langle f_{1}\left( u\left( t\right) ,v\left( t\right) \right) ,u\left(
		t\right) \rangle +\langle f_{2}\left( u\left( t\right) ,v\left( t\right)
		\right) ,v\left( t\right) \rangle .
	\end{align*}
	
	By using this formulation and the fact computed from (\ref{d9}) that%
	\begin{align*}
		L^{\prime }\left( t\right) &=\left( 1-\xi \right) H^{-\xi }\left( t\right)
		H^{\prime }\left( t\right) +\varepsilon \left( \left\Vert u^{\prime }\left(
		t\right) \right\Vert ^{2}+\left\Vert v^{\prime }\left( t\right) \right\Vert
		^{2}\right) \\&-\varepsilon \left( \left\Vert u_{x}\left( t\right) \right\Vert
		^{2}+\left\Vert v_{x}\left( t\right) \right\Vert ^{2}\right) +\varepsilon \left( K_{1}\left\vert u\left( 1,t\right) \right\vert
		^{p_{1}}+K_{2}\left\vert v\left( 0,t\right) \right\vert ^{p_{2}}\right)
		\\
		&+\varepsilon \left( \langle f_{1}\left( u\left( t\right) ,v\left( t\right)
		\right) ,u\left( t\right) \rangle +\langle f_{2}\left( u\left( t\right)
		,v\left( t\right) \right) ,v\left( t\right) \rangle \right) ,
	\end{align*}%
	in combination with (\ref{d4}), (\ref{d5}), (\ref{d8}) and the following
	inequality%
	\begin{align*}
	&\langle f_{1}\left( u\left( t\right) ,v\left( t\right) \right) ,u\left(
	t\right) \rangle +\langle f_{2}\left( u\left( t\right) ,v\left( t\right)
	\right) ,v\left( t\right) \rangle \\& \geq d_{1}\int_{0}^{1}\mathcal{F}%
	\left( u\left( x,t\right) ,v\left( x,t\right) \right) dx\geq d_{1}\bar{d}%
	_{1}\left( \left\Vert u\left( t\right) \right\Vert _{L^{\alpha }}^{\alpha
	}+\left\Vert v(t)\right\Vert _{L^{\beta }}^{\beta }\right) ,
	\end{align*}%
	we get%
	\begin{align}\tag{4.13}
	\label{d13}
	&L^{\prime }\left( t\right) \\&\geq \varepsilon \left( \left\Vert u^{\prime
	}\left( t\right) \right\Vert ^{2}+\left\Vert v^{\prime }\left( t\right)
	\right\Vert ^{2}\right) -2\varepsilon \bar{D}_{2}\left( \left\vert u\left(
	1,t\right) \right\vert ^{p_{1}}+\left\vert v\left( 0,t\right) \right\vert
	^{p_{2}}+\left\Vert u\left( t\right) \right\Vert _{L^{\alpha }}^{\alpha
	}+\left\Vert v(t)\right\Vert _{L^{\beta }}^{\beta }\right)  \nonumber \\
	&+\varepsilon \left[ K_{1}\left\vert u\left( 1,t\right) \right\vert
	^{p_{1}}+K_{2}\left\vert v\left( 0,t\right) \right\vert ^{p_{2}}+d_{1}\bar{d}%
	_{1}\left( \left\Vert u\left( t\right) \right\Vert _{L^{\alpha }}^{\alpha
	}+\left\Vert v\left( t\right) \right\Vert _{L^{\beta }}^{\beta }\right) %
	\right]  \notag \\
	&\geq \varepsilon \left( \left\Vert u^{\prime }\left( t\right) \right\Vert
	^{2}+\left\Vert v^{\prime }\left( t\right) \right\Vert ^{2}\right) \nonumber\\&
	+\varepsilon \left( \bar{D}_{1}-2\bar{D}_{2}\right) \left( \left\vert
	u\left( 1,t\right) \right\vert ^{p_{1}}+\left\vert v\left( 0,t\right)
	\right\vert ^{p_{2}}+\left\Vert u\left( t\right) \right\Vert _{L^{\alpha
	}}^{\alpha }+\left\Vert v(t)\right\Vert _{L^{\beta }}^{\beta }\right) , \nonumber
	\end{align}%
	where we have put $\bar{D}_{1}=\min \left\{ K_{1},K_{2},d_{1}\bar{d}%
	_{1}\right\} $.
	
	If one assumes that $2\max \left\{ \frac{K_{1}}{p_{1}},\frac{%
		K_{2}}{p_{2}},\bar{d}_{2}\right\} <\min \left\{ K_{1},K_{2},d_{1}\bar{d}%
	_{1}\right\} $, we deduce%
	\begin{align}\tag{4.14}  \label{d14}
	0<\bar{D}_{3}&=\bar{D}_{1}-2\bar{D}_{2}=\min \left\{ K_{1},K_{2},d_{1}\bar{d}%
	_{1}\right\} -2\max \left\{ \frac{K_{1}}{p_{1}},\frac{K_{2}}{p_{2}},\bar{d}%
	_{2}\right\} \\&<\min \left\{ K_{1},K_{2},d_{1}\bar{d}_{1}\right\} .%
	\nonumber
	\end{align}
	
	Recall from (\ref{d7}) that%
	\begin{equation*}
	\begin{tabular}{l}
	$\bar{D}_{2}\left( \left\vert u\left( 1,t\right) \right\vert
	^{p_{1}}+\left\vert v\left( 0,t\right) \right\vert ^{p_{2}}+\left\Vert
	u\left( t\right) \right\Vert _{L^{\alpha }}^{\alpha }+\left\Vert
	v(t)\right\Vert _{L^{\beta }}^{\beta }\right) \geq H\left( t\right) ,$%
	\end{tabular}%
	\end{equation*}%
	and also thanks to (\ref{d8}) and (\ref{d13}) with the assumption to get (%
	\ref{d14}), it is sufficient to show that there exists $\gamma >0$ such that
	(\ref{d12}) holds. Hence, we complete the proof of the lemma.
\end{proof}

From Lemma \ref{lem:4.1}, we obtain (\ref{d11}). The assumption $2\max \left\{ \frac{%
	K_{1}}{p_{1}},\frac{K_{2}}{p_{2}},\bar{d}_{2}\right\} <\min \left\{
K_{1},K_{2},d_{1}\bar{d}_{1}\right\} $ is additionally necessary to approach
the blow-up result. Before going to the main result, let us
consider the following supplementary inequalities, whose proof can be found in Appendix \ref{app1}.

\begin{lemma}\label{lem:4.2}
	Let $\xi >0$ such that
	\begin{equation*}
	\begin{cases}
	2\leq 2/\left( 1-2\xi \right) \leq \min \left\{ \alpha ,\beta \right\} , \\ 
	2\leq 2/\left( 1-\xi \right) \leq \min \left\{ \alpha ,\beta
	,p_{1},p_{2}\right\} ,%
	\end{cases}%
	\end{equation*}%
	then the following inequalities hold
	\begin{align}	\tag{4.15}  \label{d15}
	\left\Vert u \right\Vert _{L^{\alpha }}^{2/\left( 1-2\xi
		\right) }+\left\Vert u \right\Vert ^{2/\left( 1-\xi \right)
	}+\left\vert u\left(1\right) \right\vert ^{2/\left( 1-\xi \right) } \leq
	3\left( \left\Vert u_{x} \right\Vert ^{2}+\left\Vert u \right\Vert _{L^{\alpha }}^{\alpha }+\left\vert u\left( 1\right)
	\right\vert ^{p_{1}}\right) ,\text{ }\forall u\in \mathbb{V}_{1},
	\end{align}
	\begin{align}\tag{4.16}  \label{d16}
	\left\Vert v \right\Vert _{L^{\beta }}^{2/\left( 1-2\xi
		\right) }+\left\Vert v \right\Vert ^{2/\left( 1-\xi \right)
	}+\left\vert v\left( 0\right) \right\vert ^{2/\left( 1-\xi \right) } \leq
	3\left( \left\Vert v_{x} \right\Vert ^{2}+\left\Vert v \right\Vert _{L^{\beta }}^{\beta }+\left\vert v\left( 0\right)
	\right\vert ^{p_{2}}\right) ,\text{ }\forall v\in \mathbb{V}_{2}.
	\end{align}
\end{lemma}

\begin{theorem}\label{thm:4.3}
	Consider $r_i = 2$ with $F_i =0$ and $q_i=2,p_i > 2$ for $i=1,2$ in the problem $(P)$. Suppose that \emph{\textbf{(}$\mathbf{A}%
	_{1}$\textbf{)}}, \emph{\textbf{(}$\mathbf{A}%
	_{3}^{\prime }$\textbf{)}} and \emph{\textbf{(}$\mathbf{A}%
_{4}$\textbf{)}} hold. If $2\max \left\{ \frac{K_{1}}{p_{1}},\frac{K_{2}}{p_{2}},\bar{d}%
_{2}\right\} <\min \left\{ K_{1},K_{2},d_{1}\bar{d}_{1}\right\} $ holds in \emph{\textbf{(}$\mathbf{A}_{3}^{\prime }$\textbf{)}}, then any weak solution $(u,v)$ of $(P)$ blows up in a finite time in the sense  that there exists $T_{\ast}>0$ such that
\[
\lim_{t\to T_{\ast}^{-}}\left( \left\Vert u^{\prime }\left(
t\right) \right\Vert ^{2}+\left\Vert v^{\prime }\left( t\right) \right\Vert
^{2}+\left\Vert u_{x}\left( t\right) \right\Vert ^{2}+\left\Vert v_{x}\left(
t\right) \right\Vert ^{2}\right) =+\infty.
\]
\end{theorem}

\begin{proof}[Proof of Theorem \ref{thm:4.3}]
	By using an elementary inequality
	\begin{equation*}
	\left( \sum_{i=1}^{7}z_{i}\right)^{r}\leq
	7^{r-1}\sum_{i=1}^{7}z_{i}^{r},\text{ }\forall r>1,\text{ }x_{i}\geq 0,%
	\text{ }i=\overline{1,7},%
	\end{equation*}%
	and thanks to (\ref{d9}), (\ref{d10}), one can show that
	\begin{align}\tag{4.17}  \label{d17}
	L^{1/\left( 1-\xi \right) }\left( t\right) &\leq C\left( H\left( t\right)
	+\left\vert \langle u\left( t\right) ,u^{\prime }\left( t\right) \rangle
	\right\vert ^{1/\left( 1-\xi \right) }+\left\vert \langle v\left( t\right)
	,v^{\prime }\left( t\right) \rangle \right\vert ^{1/\left( 1-\xi \right)
	}\right.   \nonumber \\
	&\left. +\left\Vert u\left( t\right) \right\Vert ^{2/\left( 1-\xi \right)
	}+\left\Vert v\left( t\right) \right\Vert ^{2/\left( 1-\xi \right)
	}+u^{2/\left( 1-\xi \right) }\left( 1,t\right) +v^{2/\left( 1-\xi \right)
	}\left( 0,t\right) \right) ,  \notag
	\end{align}%
	where
	\begin{align*}
	C=2^{-1/\left( 1-\xi \right) }\max &\left\{ 2^{1/\left( 1-\xi \right)
	},2^{1/\left( 1-\xi \right) }\varepsilon ^{1/\left( 1-\xi \right) },%
	\left( \lambda _{1}\varepsilon \right) ^{1/\left( 1-\xi \right) }, \right. \\ & \left.
	\left( \lambda _{2}\varepsilon \right) ^{1/\left( 1-\xi \right) },%
	\left( \mu _{1}\varepsilon \right) ^{1/\left( 1-\xi \right) },\left( \mu _{2}\varepsilon \right) ^{1/\left( 1-\xi \right) }\right\}>0. 
	\end{align*}
	
	Moreover, we find%
	\begin{equation}
	\begin{tabular}{l}
	$\left\vert \langle u\left( t\right) ,u^{\prime }\left( t\right) \rangle
	\right\vert ^{1/\left( 1-\xi \right) }\leq \left\Vert u\left( t\right)
	\right\Vert ^{1/\left( 1-\xi \right) }\left\Vert u^{\prime }\left( t\right)
	\right\Vert ^{1/\left( 1-\xi \right) }\leq \left\Vert u\left( t\right)
	\right\Vert _{L^{\alpha }}^{1/\left( 1-\xi \right) }\left\Vert u^{\prime
	}\left( t\right) \right\Vert ^{1/\left( 1-\xi \right) },$%
	\end{tabular}
	\tag{4.18}  \label{d18}
	\end{equation}%
	where we have used the Cauchy-Schwartz inequality and the H\"{o}lder inequality.
	
	According to Young's inequality introduced in (\ref{c6}), by choosing $%
	\delta =1,$ $q=\frac{2\left( 1-\xi \right) }{1-2\xi },$ $q^{\prime }=2\left(
	1-\xi \right) $ and letting $a=\left\Vert u\left( t\right) \right\Vert
	_{L^{\alpha }}^{1/\left( 1-\xi \right) },$ $b=\left\Vert u^{\prime }\left(
	t\right) \right\Vert ^{1/\left( 1-\xi \right) }$ one easily obtains from (%
	\ref{d18}) that%
	\begin{equation}
	\begin{tabular}{l}
	$\left\vert \langle u\left( t\right) ,u^{\prime }\left( t\right) \rangle
	\right\vert ^{1/\left( 1-\xi \right) }\leq c_{1}\left( \left\Vert u\left(
	t\right) \right\Vert _{L^{\alpha }}^{2/\left( 1-2\xi \right) }+\left\Vert
	u^{\prime }\left( t\right) \right\Vert ^{2}\right) ,$%
	\end{tabular}
	\tag{4.19}  \label{d19}
	\end{equation}%
	where $c_{1}=\max \left\{ \frac{1-2\xi }{2\left( 1-\xi \right) },\frac{1}{%
		2\left( 1-\xi \right) }\right\} <1$. Similarly, we have for a constant $%
	c_{2}\in \left( 0,1\right) $%
	\begin{equation}
	\begin{tabular}{l}
	$\left\vert \langle v\left( t\right) ,v^{\prime }\left( t\right) \rangle
	\right\vert ^{1/\left( 1-\xi \right) }\leq c_{2}\left( \left\Vert v\left(
	t\right) \right\Vert _{L^{\beta }}^{2/\left( 1-\xi \right) }+\left\Vert
	v^{\prime }\left( t\right) \right\Vert ^{2}\right) .$%
	\end{tabular}
	\tag{4.20}  \label{d20}
	\end{equation}
	
	Therefore, there always exists a positive constant (here we reuse the
	notation $C$ for simplicity) such that%
	\begin{align}\tag{4.21}  \label{d21}
	L^{1/\left( 1-\xi \right) }\left( t\right) &\leq C\left( H\left( t\right)
	+\left\Vert u^{\prime }\left( t\right) \right\Vert ^{2}+\left\Vert v^{\prime
	}\left( t\right) \right\Vert ^{2}+\left\Vert u_{x}\left( t\right)
	\right\Vert ^{2}+\left\Vert v_{x}\left( t\right) \right\Vert ^{2}\right. 
	 \\
	&\left. +\left\Vert u\left( t\right) \right\Vert _{L^{\alpha }}^{\alpha
	}+\left\Vert v\left( t\right) \right\Vert _{L^{\beta }}^{\beta }+\left\vert
	u\left( 1,t\right) \right\vert ^{p_{1}}+\left\vert v\left( 0,t\right)
	\right\vert ^{p_{2}}\right) ,\text{ }\forall t\in \left[ 0,T_{\ast }\right) .
	\notag  
	\end{align}
	
	In addition, the estimate (\ref{d12}) together with (\ref{d21}) allows us to
	take a positive constant $\bar{C}$ such that%
	\begin{equation}
	\begin{tabular}{l}
	$L^{\prime }\left( t\right) \geq \bar{C}L^{1/\left( 1-\xi \right) }\left(
	t\right) ,\text{ \ }\forall t\in \left[ 0,T_{\ast }\right) .$%
	\end{tabular}
	\tag{4.22}  \label{d22}
	\end{equation}
	
	Now, integrating (\ref{d22}) over $\left( 0,t\right) $ one has%
	\begin{equation*}
	L^{\xi /\left( 1-\xi \right) }\left( t\right) \geq \frac{1}{L^{-\xi /\left(
			1-\xi \right) }\left( 0\right) -\frac{\bar{C}\xi }{1-\xi }t},\; t\in \left[
	0,\frac{1-\xi }{\bar{C}\xi }L^{-\xi /\left( 1-\xi \right) }\left( 0\right)
	\right) ,
	\end{equation*}%
	which yields that $L\left( t\right) $ definitely blows up in a finite time 
	$T_{\ast }=\frac{1-\xi }{\bar{C}\xi }L^{-\xi /\left( 1-\xi \right) }\left(
	0\right)$. Hence, we complete the proof of the theorem.
\end{proof}

\section{Exponential decay}
While the previous section employs the total energy to find the
finite time blow-up result, the exponential decay is proved in this section.
In particular, we study the global solution $\left( u,v\right) $ of $\left(
P\right) $ satisfying (\ref{3.71}) with $r_{1}=r_{2}=q_{1}=q_{2}=2$ and $%
p_{1},$ $p_{2}>2$. Like the blow-up phenomenon, this sort of results can be
seen as an extension of many previous works from the single wave equation
in, for example, \cite{7,16} to the system of equations.

Our result here relies on the construction of a Lyapunov functional by
performing a suitable modification of the energy. To this end, for $\delta
>0 $ being chosen later, we define
\begin{equation}
\begin{tabular}{l}
$\mathcal{L}\left( t\right) :=E\left( t\right) +\delta \psi \left( t\right) ,$%
\end{tabular}
\tag{5.1}  \label{e1}
\end{equation}%
where we have recalled the function $\psi \left( t\right) $ in (\ref{d10}).

Under some additional conditions, it is sufficient to see that $\mathcal{L}%
(t)$ and $E(t)$ are equivalent in the sense that there exist two positive
constants $\beta _{1}$ and $\beta _{2}$ depending on $\delta $ such that for 
$t\geq 0$,%
\begin{equation}
\begin{tabular}{l}
$\beta _{1}E\left( t\right) \leq \mathcal{L}\left( t\right) \leq \beta
_{2}E\left( t\right) .$%
\end{tabular}
\tag{5.2}  \label{e2}
\end{equation}

Before explicitly providing the statement of this equivalence, let us first
consider the time derivative of the total energy, which can be defined in
the same way as (\ref{d4}), by the following lemma.

\begin{lemma}\label{lem:5.1}
	The time derivative of the total energy satisfies
	\begin{align}	\tag{5.3}  \label{e3}
	E^{\prime }\left( t\right) &\leq -\lambda _{\ast }\left( \left\Vert
	u^{\prime }\left( t\right) \right\Vert ^{2}+\left\Vert v^{\prime }\left(
	t\right) \right\Vert ^{2}\right) -\mu _{\ast }\left( \left\vert u^{\prime
	}\left( 1,t\right) \right\vert ^{2}+\left\vert v^{\prime }\left( 0,t\right)
	\right\vert ^{2}\right)  \\ &
	+\frac{1}{2}\left( \left\Vert F_{1}\left( t\right)
	\right\Vert +\left\Vert F_{2}\left( t\right) \right\Vert \right) +\frac{1}{2}%
	\left( \left\Vert F_{1}\left( t\right) \right\Vert +\left\Vert F_{2}\left(
	t\right) \right\Vert \right) \left( \left\Vert u^{\prime }\left( t\right)
	\right\Vert ^{2}+\left\Vert v^{\prime }\left( t\right) \right\Vert
	^{2}\right) , \nonumber 
	\end{align}%
	\begin{align}\tag{5.4}  \label{e4}
	E^{\prime }\left( t\right) &\leq -\left( \lambda _{\ast }-\frac{\varepsilon
		_{1}}{2}\right) \left( \left\Vert u^{\prime }\left( t\right) \right\Vert
	^{2}+\left\Vert v^{\prime }\left( t\right) \right\Vert ^{2}\right) -\mu
	_{\ast }\left( \left\vert u^{\prime }\left( 1,t\right) \right\vert
	^{2}+\left\vert v^{\prime }\left( 0,t\right) \right\vert ^{2}\right)  \\ &
	+\frac{1}{2\varepsilon _{1}}\left( \left\Vert F_{1}\left(
	t\right) \right\Vert ^{2}+\left\Vert F_{2}\left( t\right) \right\Vert
	^{2}\right) ,\nonumber 
	\end{align}%
\end{lemma}

\begin{proof}[Proof of Lemma \ref{lem:5.1}]
	Multiplying (\ref{main1}) by $\left(
	u^{\prime }\left( x,t\right) ,v^{\prime }\left( x,t\right) \right) $,
	integrating over $\Omega $ and using the integration by parts, we obtain%
	\begin{align}\tag{5.5}
	E^{\prime }(t)&=-\lambda _{1}\left\Vert u^{\prime }(t)\right\Vert
	^{2}-\lambda _{2}\left\Vert v^{\prime }(t)\right\Vert ^{2}-\mu
	_{1}\left\vert u^{\prime }(1,t)\right\vert ^{2}\\&-\mu _{2}\left\vert v^{\prime
	}(0,t)\right\vert ^{2}+\langle F_{1}(t),u^{\prime }(t)\rangle +\langle
	F_{2}(t),v^{\prime }(t)\rangle .\nonumber 
	\end{align}
	
	To obtain (\ref{e3}), we only need to estimate the last two terms. It is
	straightforward that by using the standard inequalities which read as%
	\begin{align*}
	\langle F_{1}(t),u^{\prime }(t)\rangle \leq \frac{1}{2}\left\Vert
	F_{1}(t)\right\Vert +\frac{1}{2}\left\Vert F_{1}(t)\right\Vert \left\Vert
	u^{\prime }(t)\right\Vert ^{2}, \\ 
	\langle F_{2}(t),v^{\prime }(t)\rangle \leq \frac{1}{2}\left\Vert
	F_{2}(t)\right\Vert +\frac{1}{2}\left\Vert F_{2}(t)\right\Vert \left\Vert
	v^{\prime }(t)\right\Vert ^{2},
	\end{align*}%
	it yields%
	\begin{align*}
	&\langle F_{1}(t),u^{\prime }(t)\rangle +\langle F_{2}(t),v^{\prime
	}(t)\rangle \\&\leq \frac{1}{2}\left( \left\Vert F_{1}(t)\right\Vert
	+\left\Vert F_{2}(t)\right\Vert \right) +\frac{1}{2}\left( \left\Vert
	F_{1}(t)\right\Vert +\left\Vert F_{2}(t)\right\Vert \right) \left(
	\left\Vert u^{\prime }(t)\right\Vert ^{2}+\left\Vert v^{\prime
	}(t)\right\Vert ^{2}\right) .
	\end{align*}
	
	For the second estimate (\ref{e4}), we also use the same approach to prove.
	Indeed, by Young's inequality one deduces the following estimate%
	\begin{align*}
	&\langle F_{1}(t),u^{\prime }(t)\rangle +\langle F_{2}(t),v^{\prime
	}(t)\rangle \\&\leq \frac{1}{2\varepsilon _{1}}\left( \left\Vert
	F_{1}(t)\right\Vert ^{2}+\left\Vert F_{2}(t)\right\Vert ^{2}\right) +\frac{%
		\varepsilon _{1}}{2}\left( \left\Vert u^{\prime }(t)\right\Vert
	^{2}+\left\Vert v^{\prime }(t)\right\Vert ^{2}\right) .
	\end{align*}
	
	Hence, the proof of the lemma is complete.
\end{proof}

Let us next define the following functions $I_{i}\left( t\right)
:=I_{i}\left( u\left( t\right) \right) $ for $i=1,2$ and $J\left( t\right)
:=J\left( u\left( t\right) \right) $ by rewriting the total energy $E\left(
t\right) $:%
\begin{equation}
E\left( t\right) =\frac{1}{2}\left( \left\Vert u^{\prime }\left( t\right)
\right\Vert ^{2}+\left\Vert v^{\prime }\left( t\right) \right\Vert
^{2}\right) +J\left( t\right) ,%
\tag{5.6}  \label{e6}
\end{equation}%
\begin{equation}
J\left( t\right) =\frac{1}{2}\left( 1-\frac{1}{p_{1}}-\frac{1}{p_{2}}%
\right) \left( \left\Vert u_{x}\left( t\right) \right\Vert ^{2}+\left\Vert
v_{x}\left( t\right) \right\Vert ^{2}\right) +\frac{I_{1}\left( t\right) }{%
	p_{1}}+\frac{I_{2}\left( t\right) }{p_{2}},
\tag{5.7}  \label{e7}
\end{equation}%
\begin{equation}
I_{1}\left( t\right) =\frac{1}{2}\left( \left\Vert u_{x}\left( t\right)
\right\Vert ^{2}+\left\Vert v_{x}\left( t\right) \right\Vert ^{2}\right)
-K_{1}\left\vert u\left( 1,t\right) \right\vert ^{p_{1}}-\frac{p_{1}}{2}%
\int_{0}^{1}\mathcal{F}\left( u\left( x,t\right) ,v\left(
x,t\right) \right) dx,
\tag{5.8}  \label{e8}
\end{equation}%
\begin{equation}
I_{2}\left( t\right) =\frac{1}{2}\left( \left\Vert u_{x}\left( t\right)
\right\Vert ^{2}+\left\Vert v_{x}\left( t\right) \right\Vert ^{2}\right)
-K_{2}\left\vert v\left( 0,t\right) \right\vert ^{p_{2}}-\frac{p_{2}}{2}%
\int_{0}^{1}\mathcal{F}\left( u\left( x,t\right) ,v\left(
x,t\right) \right) dx.
\tag{5.9}  \label{e9}
\end{equation}

Furthermore, let us provide the following assumption:

\textbf{(}$\mathbf{B}_{2}^{\prime }$\textbf{)} $F_{1},$ $F_{2}\in L^{\infty
}\left( \mathbb{R}_{+};L^{2}\right) \cap L^{1}\left( \mathbb{R}%
_{+};L^{2}\right) $.

From now on, our main result in this section is established where the proof
of the equivalence between $\mathcal{L}\left( t\right) $ and $E\left(
t\right) $ is also included. It says that if there is an exponential rate of
energy decay for the external functions $F_{1}$ and $F_{2}$ and influenced
by such functions, the initial-related energy function, denoted by $E_{\ast }
$, is properly suited in a particular set, then the quadratic-type total
energy decays exponentially.

\begin{theorem}\label{thm:5.2}
	Consider $\left( u,v\right) $ of the problem $\left(
	P\right) $ satisfying (\ref{3.71}) with $r_{1}=r_{2}=q_{1}=q_{2}=2$ and $%
	p_{1},$ $p_{2}>2$. Suppose that \emph{\textbf{(}$\mathbf{A}_{\mathbf{1%
		}}$\textbf{)}}, \emph{\textbf{(}$\mathbf{A}_{3}^{\prime }$\textbf{)}} and \emph{\textbf{(}$\mathbf{B}_{2}^{\prime }$\textbf{)}} hold, together with $d_{2}<\min \left\{ p_{1},p_{2}\right\} $ in \emph{\textbf{(}$%
	\mathbf{A}_{3}^{\prime }$\textbf{)}}.  Assume that $I_{1}\left(
0\right), I_{2}\left( 0\right) >0$ and the initial energy $E(0)$ satisfies
\begin{equation}
\eta _{\ast }=\frac{p_{1}+p_{2}}{2}\bar{d}_{2}\left[ \left( p_{\ast
}E_{\ast }\right) ^{\frac{\alpha }{2}-1}+\left( p_{\ast }E_{\ast }\right) ^{%
	\frac{\beta }{2}-1}\right] +K_{1}\left( p_{\ast }E_{\ast }\right) ^{\frac{%
		p_{1}}{2}-1}+K_{2}\left( p_{\ast }E_{\ast }\right) ^{\frac{p_{2}}{2}-1}<1,
\tag{5.10}  \label{e10}
\end{equation}%
where $p_{\ast }=\frac{2p_{1}p_{2}}{p_{1}p_{2}-p_{1}-p_{2}}$, $%
E_{\ast }=\left[ E\left( 0\right) +\rho \right] \mbox{exp}\left( 2\rho
\right) $, $\rho =\frac{1}{2}\int_{0}^{\infty }\left( \left\Vert
F_{1}\left( s\right) \right\Vert +\left\Vert F_{2}\left( s\right)
\right\Vert \right) ds$. Moreover, assume the external
functions decays exponentially in the sense that
\begin{equation}
\left\Vert F_{1}\left( t\right) \right\Vert ^{2}+\left\Vert F_{2}\left(
t\right) \right\Vert ^{2}\leq \eta _{1}\mbox{exp}\left( -\eta _{2}t\right) ,%
\text{ \ }\forall t\geq 0,  \tag{5.11}  \label{e11}
\end{equation}
where $\eta _{1}$ and $\eta _{2}$ are two positive constants. Then there exist positive constants $C$ and $\gamma$ such that
\begin{equation}
E\left( t\right) \leq C\mbox{exp}\left( -\gamma t\right) ,\text{ \ }\forall
t\geq 0.  \tag{5.12}  \label{e12}
\end{equation}
\end{theorem}

\begin{proof}[Proof of Theorem \ref{thm:5.2}]
	First, we claim that $I_{i}\left( t\right) \geq 0$ for $i=1,2
	$ and for all $t\geq 0$. In fact, since $I_{i}\left( t\right) ,$ $i=1,2$ is
	continuous and their initial values are positive, thus there exist two positive
	constants $T_{1}$ and $T_{2}$ such that%
	\begin{equation*}
	\begin{tabular}{l}
	$I_{1}\left( t\right) \geq 0,\quad \forall t\in \left[ 0,T_{1}\right] ,\;%
	\mbox{and}\;I_{2}\left( t\right) \geq 0,\quad \forall t\in \left[ 0,T_{2}%
	\right] $%
	\end{tabular}%
	\end{equation*}%
	which also leads to the fact that%
	\begin{equation*}
	J\left( t\right) \geq \frac{1}{2}\left( 1-\frac{1}{p_{1}}-\frac{1}{p_{2}}%
	\right) \left( \left\Vert u_{x}\left( t\right) \right\Vert ^{2}+\left\Vert
	v_{x}\left( t\right) \right\Vert ^{2}\right) , \; \forall t\in \left[ 0,%
	\bar{T}\right] ,%
	\end{equation*}%
	where $\bar{T}=\min \left\{ T_{1},T_{2}\right\} >0$. In other words, we can say that
	\begin{equation}
	\left\Vert u_{x}\left( t\right) \right\Vert ^{2}+\left\Vert v_{x}\left(
	t\right) \right\Vert ^{2}\leq \frac{2p_{1}p_{2}}{p_{1}p_{2}-p_{1}-p_{2}}%
	J\left( t\right) \leq \frac{2p_{1}p_{2}}{p_{1}p_{2}-p_{1}-p_{2}}E\left(
	t\right) ,\; \forall t\in \left[ 0,\bar{T}\right] .%
	\tag{5.13}  \label{e13}
	\end{equation}
	
	By (\ref{e3}) and thanks to $I_{i}\left( 0\right) >0$ for $i=1,2$ associated
	with (\ref{e6}), we deduce
	\begin{align*}
	E\left( t\right) &\leq E\left( 0\right) +\frac{1}{2}\int_{0}^{%
		\infty }\left( \left\Vert F_{1}\left( s\right) \right\Vert +\left\Vert
	F_{2}\left( s\right) \right\Vert \right) ds\\&+\int_{0}^{t}\left(
	\left\Vert F_{1}\left( s\right) \right\Vert +\left\Vert F_{2}\left( s\right)
	\right\Vert \right) E\left( s\right) ds,\; \forall t\in \left[ 0,\bar{T}%
	\right] ,
	\end{align*}%
	then using Gronwall's inequality leads to the following
	\begin{align}\tag{5.14}  \label{e14}
	E\left( t\right) &\leq \left[ E\left( 0\right) +\frac{1}{2}\int_{0}^{\infty
	}\left( \left\Vert F_{1}\left( s\right) \right\Vert +\left\Vert F_{2}\left(
	s\right) \right\Vert \right) ds\right] \mbox{exp}\left( \int_{0}^{t}\left(
	\left\Vert F_{1}\left( s\right) \right\Vert +\left\Vert F_{2}\left( s\right)
	\right\Vert \right) ds\right)   \\
	&\leq \left[ E\left( 0\right) +\rho \right] \mbox{exp}\left( 2\rho \right).
	\notag  
	\end{align}
	
	Therefore, with $E_{\ast }=\left[
	E\left( 0\right) +\rho \right] \mbox{exp}\left( 2\rho \right) ,$ $p_{\ast }=%
	\frac{2p_{1}p_{2}}{p_{1}p_{2}-p_{1}-p_{2}}$ we combine (\ref{e14}) with (\ref{e13}) to obtain
	\begin{equation}
	\begin{tabular}{l}
	$\left\Vert u_{x}\left( t\right) \right\Vert ^{2}+\left\Vert v_{x}\left(
	t\right) \right\Vert ^{2}\leq p_{\ast }E_{\ast },\text{ \ }\forall t\in %
	\left[ 0,\bar{T}\right] .$%
	\end{tabular}
	\tag{5.15}  \label{e15}
	\end{equation}
	
	By the assumption (\ref{d2}) in ($\mathbf{A}_{3}^{\prime }$) and also (\ref%
	{e15}), it yields
	\begin{align}	\tag{5.16}  \label{e16}
	K_{1}\left\vert u\left( 1,t\right) \right\vert ^{p_{1}}&+K_{2}\left\vert
	v\left( 0,t\right) \right\vert ^{p_{2}}
	\leq K_{1}\left\Vert u_{x}\left(
	t\right) \right\Vert ^{p_{1}}+K_{2}\left\Vert v_{x}\left( t\right)
	\right\Vert ^{p_{2}} \\ &
	\leq K_{1}\left( p_{\ast }E_{\ast }\right) ^{\frac{p_{1}}{2}-1}\left\Vert
	u_{x}\left( t\right) \right\Vert ^{2}+K_{2}\left( p_{\ast }E_{\ast }\right)
	^{\frac{p_{2}}{2}-1}\left\Vert v_{x}\left( t\right) \right\Vert ^{2},%
	\nonumber 
	\end{align}%
	and
	\begin{align}	\tag{5.17}  \label{e17}
	\left( p_{1}+p_{2}\right)& \int_{0}^{1}\mathcal{F}\left( u\left(
	x,t\right) ,v\left( x,t\right) \right) dx\leq \left( p_{1}+p_{2}\right) \bar{%
		d}_{2}\left( \left\Vert u\left( t\right) \right\Vert _{L^{\alpha }}^{\alpha
	}+\left\Vert v\left( t\right) \right\Vert _{L^{\beta }}^{\beta }\right)  \\ &
	\leq \left( p_{1}+p_{2}\right) \bar{d}_{2}%
	\left[ \left( p_{\ast }E_{\ast }\right) ^{\frac{\alpha }{2}-1}\left\Vert
	u_{x}\left( t\right) \right\Vert ^{2}+\left( p_{\ast }E_{\ast }\right) ^{%
		\frac{\beta }{2}-1}\left\Vert v_{x}\left( t\right) \right\Vert ^{2}\right] .\nonumber %
	\end{align}
	
	Thus, it follows from (\ref{e16}) and (\ref{e17}) that
	\begin{align*}
	&K_{1}\left\vert u\left( 1,t\right) \right\vert ^{p_{1}}+K_{2}\left\vert
	v\left( 0,t\right) \right\vert ^{p_{2}}+\frac{p_{1}+p_{2}}{2}%
	\int_{0}^{1}\mathcal{F}\left( u\left( x,t\right) ,v\left(
	x,t\right) \right) dx
	\\&\leq \eta _{\ast }\left( \left\Vert u_{x}\left(
	t\right) \right\Vert ^{2}+\left\Vert v_{x}\left( t\right) \right\Vert
	^{2}\right)  
	 <\left\Vert
	u_{x}\left( t\right) \right\Vert ^{2}+\left\Vert v_{x}\left( t\right)
	\right\Vert ^{2}, \;\forall t\in \left[ 0,\bar{T}\right] ,%
	\end{align*}%
	where $\eta _{\ast }$ is given by (\ref{e10}).
	
	So we now claim that $I_{1}\left( t\right) $ and $I_{2}\left( t\right) $ are
	positive for all $t\in \left[ 0,\bar{T}\right] $. If we put%
	\begin{equation*}
	\begin{tabular}{l}
	$T_{\ast }=\sup \left\{ T>0:I_{1}\left( t\right) \;\mbox{and}\;I_{2}\left(
	t\right) \;\mbox{are positive}\;\forall t\in \left[ 0,T\right] \right\} ,$%
	\end{tabular}%
	\end{equation*}%
	and if $T_{\ast }<\infty $, then (by the continuity of $I_{1}\left( t\right) 
	$ and $I_{2}\left( t\right) $) we have $I_{1}\left( T_{\ast }\right) $ and $%
	I_{2}\left( T_{\ast }\right) $ are non-negative. Therefore, by the same
	arguments it is possible to show that there exists $\bar{T}_{\ast }>T_{\ast
	}$ such that $I_{i}\left( t\right) >0$ for $i=1,2$ for all $t\in \left[ 0,%
	\bar{T}_{\ast }\right] $. Hence, we can conclude that $I_{i}\left( t\right)
	\geq 0\;\left( i=1,2\right) $ for all $t\geq 0$.
	
	Next, we end up with the equivalence of $\mathcal{L}\left( t\right) $ and $%
	E\left( t\right) $. Let us recall that by (\ref{e1}), (\ref{e3}), (\ref{e4}%
	), (\ref{e6}), (\ref{e7}) and (\ref{d10}) the function $\mathcal{L}\left(
	t\right) $ can be rewritten as%
	\begin{align*}
	\mathcal{L}\left( t\right) &=\frac{1}{2}\left( \left\Vert u^{\prime }\left(
	t\right) \right\Vert ^{2}+\left\Vert v^{\prime }\left( t\right) \right\Vert
	^{2}\right) +\frac{1}{2}\left( 1-\frac{1}{p_{1}}-\frac{1}{p_{2}}\right)
	\left( \left\Vert u_{x}\left( t\right) \right\Vert ^{2}+\left\Vert
	v_{x}\left( t\right) \right\Vert ^{2}\right) \\ & +\frac{I_{1}\left( t\right) }{%
		p_{1}}+\frac{I_{2}\left( t\right) }{p_{2}} 
	+\delta \langle u\left( t\right) ,u^{\prime }\left( t\right)
	\rangle +\delta \langle v\left( t\right) ,v^{\prime }\left( t\right) \rangle
	\\&+\frac{\delta \lambda _{1}}{2}\left\Vert u\left( t\right) \right\Vert ^{2}+%
	\frac{\delta \lambda _{2}}{2}\left\Vert v\left( t\right) \right\Vert ^{2}+%
	\frac{\delta \mu _{1}}{2}u^{2}\left( 1,t\right) +\frac{\delta \mu _{2}}{2}%
	v^{2}\left( 0,t\right) .
	\end{align*}
	
	It is straightforward to see that%
	\begin{align}	\tag{5.18}  \label{e18}
	\mathcal{L}\left( t\right) &\leq \frac{1}{2}\left( 1+\delta \right) \left(
	\left\Vert u^{\prime }\left( t\right) \right\Vert ^{2}+\left\Vert v^{\prime
	}\left( t\right) \right\Vert ^{2}\right) +\frac{I_{1}\left( t\right) }{p_{1}}%
	+\frac{I_{2}\left( t\right) }{p_{2}} \\ &
	+\frac{1}{2}\left( 1-\frac{1}{p_{1}}-\frac{1}{p_{2}}+\delta
	\left( \lambda _{1}+\lambda _{2}+\mu _{1}+\mu _{2}\right) \right) \left(
	\left\Vert u_{x}\left( t\right) \right\Vert ^{2}+\left\Vert v_{x}\left(
	t\right) \right\Vert ^{2}\right), \nonumber 
	\end{align}%
	where we have used an elementary inequality%
	\begin{equation*}
	\langle u\left( t\right) ,u^{\prime }\left( t\right) \rangle +\langle
	v\left( t\right) ,v^{\prime }\left( t\right) \rangle \leq \frac{1}{2}\left(
	\left\Vert u_{x}\left( t\right) \right\Vert ^{2}+\left\Vert u_{x}^{\prime
	}(t)\right\Vert ^{2}\right) +\frac{1}{2}\left( \left\Vert v_{x}\left(
	t\right) \right\Vert ^{2}+\left\Vert v_{x}^{\prime }\left( t\right)
	\right\Vert ^{2}\right) .
	\end{equation*}
	
	Therefore, one can choose from (\ref{e18}) that%
	\begin{align*}
	\beta _{2}&=\max \left\{1+\delta , \frac{1-\frac{1}{p_{1}}-\frac{1}{p_{2}}%
		+\delta \left( \lambda _{1}+\lambda _{2}+\mu _{1}+\mu _{2}\right) }{1-\frac{1%
		}{p_{1}}-\frac{1}{p_{2}}}\right\}
	\\&=\max \left\{1+\delta , 1+\frac{\delta \left(
		1+\lambda _{1}+\lambda _{2}+\mu _{1}+\mu _{2}\right) }{1-\frac{1}{p_{1}}-%
		\frac{1}{p_{2}}}\right\},
	\end{align*}%
	to obtain $\mathcal{L}\left( t\right) \leq \beta _{2}E\left( t\right) $.
	
	Furthermore, one easily finds that%
	\begin{align*}
	\mathcal{L}\left( t\right) &\geq \frac{1}{2}\left( 1-\delta \right) \left(
	\left\Vert u^{\prime }\left( t\right) \right\Vert ^{2}+\left\Vert v^{\prime
	}\left( t\right) \right\Vert ^{2}\right) \\&+\frac{I_{1}}{2p_{1}}+\frac{I_{2}}{%
		2p_{2}}+\frac{1}{2}\left( 1-\frac{1}{p_{1}}-\frac{1}{p_{2}}-\delta \right)
	\left( \left\Vert u_{x}\left( t\right) \right\Vert ^{2}+\left\Vert
	v_{x}\left( t\right) \right\Vert ^{2}\right) .
	\end{align*}
	
	Thus, choosing $\delta >0$ small enough for which $\beta _{1}=\min
	\left\{1-\delta ,1-\frac{\delta }{1-\frac{1}{p_{1}}-\frac{1}{p_{2}}}\right\}>0,$ most
	likely we choose $\delta \in \left( 0,1-\frac{\delta }{1-\frac{1}{p_{1}}-%
		\frac{1}{p_{2}}}\right) $ to obtain $\mathcal{L}\left( t\right) \leq \beta
	_{2}E\left( t\right) $. Hereby, our Lyapunov functional $\mathcal{L}\left(
	t\right) $ is definitely equivalent to the total energy $E\left( t\right) $
	for all $t\geq 0$.
	
	It remains to consider the functional $\psi \left( t\right) $ defined in (%
	\ref{d10}). The time derivative of $\psi \left( t\right) $ can be found by
	multiplying (\ref{main1}) by $\left( u\left( x,t\right) ,v\left( x,t\right)
	\right) $ and then integrating the resulting equation over $\Omega $. It
	therefore has the following expression%
	\begin{align}	\tag{5.19}  \label{e19}
	\psi ^{\prime }\left( t\right) &=\left\Vert u^{\prime }\left( t\right)
	\right\Vert ^{2}+\left\Vert v^{\prime }\left( t\right) \right\Vert
	^{2}-\left\Vert u_{x}\left( t\right) \right\Vert ^{2}-\left\Vert v_{x}\left(
	t\right) \right\Vert ^{2}+K_{1}\left\vert u\left( 1,t\right) \right\vert
	^{p_{1}} \\&+K_{2}\left\vert v\left( 0,t\right) \right\vert ^{p_{2}}
	+\langle f_{1}\left( u\left( t\right) ,v\left( t\right)
	\right) ,u\left( t\right) \rangle +\langle f_{2}\left( u\left( t\right)
	,v\left( t\right) \right) ,v\left( t\right) \rangle \nonumber \\  & +\langle F_{1}\left(
	t\right) ,u\left( t\right) \rangle +\langle F_{2}\left( t\right) ,v\left(
	t\right) \rangle .\nonumber 
	\end{align}
	
	On the one hand, we have%
	\begin{equation*}
	\langle f_{1}\left( u\left( t\right) ,v\left( t\right) \right) ,u\left(
	t\right) \rangle +\langle f_{2}\left( u\left( t\right) ,v\left( t\right)
	\right) ,v\left( t\right) \rangle \leq d_{2}\int_{0}^{1}\mathcal{F}%
	\left( u\left( x,t\right) ,v\left( x,t\right) \right) dx,
	\end{equation*}%
	and by (\ref{e8}) and (\ref{e9}), we continue to estimate the above
	inequality by%
	\begin{align}\tag{5.20}  \label{e20}
	&\langle f_{1}\left( u\left( t\right) ,v\left( t\right) \right) ,u\left(
	t\right) \rangle +\langle f_{2}\left( u\left( t\right) ,v\left( t\right)
	\right) ,v\left( t\right) \rangle \\ & \leq d_{2}\left[ \frac{1}{2}\left( \frac{1%
	}{p_{1}}+\frac{1}{p_{2}}\right) \left( \left\Vert u_{x}\left( t\right)
	\right\Vert ^{2}+\left\Vert v_{x}\left( t\right) \right\Vert ^{2}\right) \right. \nonumber  \\ & \left. 
	 -\left( \frac{K_{1}}{p_{1}}\left\vert u\left(
	1,t\right) \right\vert ^{p_{1}}+\frac{K_{2}}{p_{2}}\left\vert v\left(
	0,t\right) \right\vert ^{p_{2}}\right) -\left( \frac{I_{1}\left( t\right) }{%
		p_{1}}+\frac{I_{2}\left( t\right) }{p_{2}}\right) \right] .
	\nonumber 
	\end{align}
	
	On the other hand, one easily has%
	\begin{align}	\tag{5.21}  \label{e21}
	\langle F_{1}\left( t\right) ,u\left( t\right) \rangle +\langle F_{2}\left(
	t\right) ,v\left( t\right) \rangle &\leq \frac{\varepsilon _{2}}{2}\left(
	\left\Vert u_{x}\left( t\right) \right\Vert ^{2}+\left\Vert v_{x}\left(
	t\right) \right\Vert ^{2}\right) \\ &
	+\frac{1}{2\varepsilon _{2}}\left( \left\Vert F_{1}\left( t\right)
	\right\Vert ^{2}+\left\Vert F_{2}\left( t\right) \right\Vert ^{2}\right) ,%
	\; \forall \varepsilon _{2}>0.\nonumber
	\end{align}
	
	Thus, we obtain from (\ref{e19})-(\ref{e21}) that%
	\begin{align}	\tag{5.22}  \label{e22}
	\psi ^{\prime }\left( t\right) &\leq \left\Vert u^{\prime }\left( t\right)
	\right\Vert ^{2}+\left\Vert v^{\prime }\left( t\right) \right\Vert
	^{2}-\left( 1-\frac{\varepsilon _{2}}{2}-\frac{d_{2}}{2}\left( \frac{1}{p_{1}%
	}+\frac{1}{p_{2}}\right) \right) \left( \left\Vert u_{x}\left( t\right)
	\right\Vert ^{2}+\left\Vert v_{x}\left( t\right) \right\Vert ^{2}\right) 
	\\ 
	&+\left( 1-\frac{d_{2}}{p_{1}}\right) K_{1}\left\vert
	u\left( 1,t\right) \right\vert ^{p_{1}}+\left( 1-\frac{d_{2}}{p_{2}}\right)
	K_{2}\left\vert v\left( 0,t\right) \right\vert ^{p_{2}} \nonumber  \\ 
	& -d_{2}\left( \frac{I_{1}\left( t\right) }{p_{1}}+\frac{%
		I_{2}\left( t\right) }{p_{2}}\right) +\frac{1}{2\varepsilon _{2}}\left(
	\left\Vert F_{1}\left( t\right) \right\Vert ^{2}+\left\Vert F_{2}\left(
	t\right) \right\Vert ^{2}\right) .\nonumber
	\end{align}
	
	By (\ref{e4}), (\ref{e22}) and thanks to (\ref{e10}) and (\ref{e16}), the
	time derivative of our Lyapunov functional $\mathcal{L}\left( t\right) $ can
	be estimated as follows:%
	\begin{align}	\tag{5.23}  \label{e23}
	\mathcal{L}^{\prime }\left( t\right) & \leq -\left( \lambda _{\ast }-\delta -%
	\frac{\varepsilon _{1}}{2}\right) \left( \left\Vert u^{\prime }\left(
	t\right) \right\Vert ^{2}+\left\Vert v^{\prime }\left( t\right) \right\Vert
	^{2}\right) -\delta d_{2}\left( \frac{I_{1}\left( t\right) }{p_{1}}+\frac{%
		I_{2}\left( t\right) }{p_{2}}\right)  \\ &
	 +\frac{1}{2}\left( \frac{1}{\varepsilon _{1}}+\frac{\delta }{%
		\varepsilon _{2}}\right) \left( \left\Vert F_{1}\left( t\right) \right\Vert
	^{2}+\left\Vert F_{2}\left( t\right) \right\Vert ^{2}\right) \nonumber  \\& -\delta \left[ \left( 1-\eta _{\ast }\right) \left( 1-\frac{%
		d_{2}}{\max \left\{ p_{1},p_{2}\right\} }\right) -\frac{\varepsilon _{2}}{2}%
	\right] \left( \left\Vert u_{x}\left( t\right) \right\Vert ^{2}+\left\Vert
	v_{x}\left( t\right) \right\Vert ^{2}\right) 
	\nonumber 
	\end{align}%
	for all $\delta ,$ $\varepsilon _{1},$ $\varepsilon _{2}>0$. Here we imply%
	\begin{equation*}
	d_{2}<\min \left\{ p_{1},p_{2}\right\} ,\text{ \ }0<\varepsilon
	_{2}<2\left( 1-\eta _{\ast }\right) \left( 1-\frac{d_{2}}{\max \left\{
		p_{1},p_{2}\right\} }\right) .
	\end{equation*}
	
	Then for $\delta $ sufficiently small, satisfying $0<\delta <\lambda _{\ast }$,
	and let $\varepsilon _{1}>0$ such that $0<\varepsilon _{1}<2\left( \lambda
	_{\ast }-\delta \right) $, the equivalence of $\mathcal{L}\left( t\right) $
	and $E\left( t\right) $ together with (\ref{e23}) and (\ref{e11}) says that
	there exists a constant $\gamma \in \left( 0,\eta _{2}\right) $ such that%
	\begin{equation*}
	\begin{tabular}{l}
	$\mathcal{L}^{\prime }\left( t\right) \leq -\gamma \mathcal{L}\left(
	t\right) +\bar{\eta}_{1}\exp \left( -\eta _{2}t\right) ,$ $\ \forall t\geq
	0, $%
	\end{tabular}%
	\end{equation*}%
	and by Gronwall's inequality, we have%
	\begin{equation*}
	\begin{tabular}{l}
	$\mathcal{L}\left( t\right) \leq \left( \mathcal{L}\left( 0\right) +\frac{%
		\bar{\eta}_{1}}{\eta _{2}}\right) \exp \left( -\gamma t\right) ,$ $\forall
	t\geq 0,$%
	\end{tabular}%
	\end{equation*}%
	which leads to (\ref{e12}) by the equivalence of $\mathcal{L}\left( t\right) 
	$ and $E\left( t\right) $. Hence, we complete the proof of the theorem.
\end{proof}

\section{A numerical example}
The one-dimensional linear damped system of nonlinear wave equations
has been qualitatively investigated in the sense of exponential decays. Therefore, in this section the emphasis is put on the illustrative
framework. Particularly, our problem here reduces to the following form:
\begin{equation*}
\begin{cases}
u_{tt}-u_{xx}+\lambda _{1}u_{t}=f_{1}\left( u,v\right) +F_{1}\left(
x,t\right) , \\ 
v_{tt}-v_{xx}+\lambda _{2}v_{t}=f_{2}\left( u,v\right) +F_{2}\left(
x,t\right) ,%
\end{cases}%
\end{equation*}%
along with the boundary conditions (\ref{eq:main2}) and initial conditions (%
\ref{main3}). Notice that the interior sources $f_{1}$ and $f_{2}$ have
been introduced in Remark \ref{rem:3.1}. In fact, the functional $\mathcal{F}\in
C^{2}\left( \mathbb{R}^{2};\mathbb{R}\right) $ (cf. (\ref{b5}%
) as an example) can be established by%
\begin{align*}
&f_{1}\left( u,v\right) =\alpha \left( \gamma _{1}\left\vert u\right\vert
^{\alpha -2}+\frac{\gamma _{2}}{2}\left\vert u\right\vert ^{\frac{\alpha }{2}%
	-2}\left\vert v\right\vert ^{\frac{\beta }{2}}\right) u,\text{ \ }\\&
f_{2}\left( u,v\right) =\beta \left( \gamma _{1}\left\vert v\right\vert
^{\beta -2}+\frac{\gamma _{2}}{2}\left\vert u\right\vert ^{\frac{\alpha }{2}%
}\left\vert v\right\vert ^{\frac{\beta }{2}-2}\right) v.
\end{align*}

By virtue of these explicit expressions, (\ref{b3}) obviously holds. Moreover, we
deduce the fact that for all $\left( u,v\right) \in \mathbb{R}^{2}$%
\begin{equation*}
\left( \gamma _{1}-\frac{\gamma _{2}}{2}\right) \left( \left\vert
u\right\vert ^{\alpha }+\left\vert v\right\vert ^{\beta }\right) \leq 
\mathcal{F}\left( u,v\right) \leq \left( \gamma _{1}+\frac{\gamma _{2}}{2}%
\right) \left( \left\vert u\right\vert ^{\alpha }+\left\vert v\right\vert
^{\beta }\right) ,
\end{equation*}%
\begin{equation*}
\min \left\{ \alpha ,\beta \right\} \mathcal{F}\left( u,v\right) \leq
uf_{1}\left( u,v\right) +vf_{2}\left( u,v\right) \leq \max \left\{ \alpha
,\beta \right\} \mathcal{F}\left( u,v\right) ,
\end{equation*}%
which imply $d_{1}=\min \left\{ \alpha ,\beta \right\} ,$ $d_{2}=\max
\left\{ \alpha ,\beta \right\} ,$ $\bar{d}_{1}=\gamma _{1}-\frac{\gamma _{2}%
}{2}$, and $\bar{d}_{2}=\gamma _{1}+\frac{\gamma _{2}}{2}$ in (\ref{d1}), (%
\ref{d2}). If one can choose an appropriate set of constants $\alpha ,$ $%
\beta ,$ $\gamma _{1},$ and $\gamma _{2}$ such that $\gamma _{2}<2\gamma
_{1} $ and $d_{2}<\min \left\{ p_{1},p_{2}\right\} $, \textbf{(}$\mathbf{A}%
_{3}^{\prime }$\textbf{)} is clearly valid.

In this example, we take $\alpha =\beta =4,$ $p_{1}=p_{2}=6,$ $q_{1}=q_{2}=2$
and $K_{i}=\mu _{i}=\lambda _{i}=1$ for $i=1,2$ with $\gamma _{1}=\frac{3}{4}%
,$ $\gamma _{2}=\frac{1}{2}$ (so, $d_{2}=4$ and $\bar{d}_{2}=1$). Next, the
initial conditions are given by%
\begin{equation*}
\begin{tabular}{l}
$\tilde{u}_{0}\left( x\right) =x\left( e^{9}+1\right) ^{\frac{-1}{4}},$ $\ 
\tilde{u}_{1}\left( x\right) =-xe^{9}\left( e^{9}+1\right) ^{\frac{-5}{4}},$
\\ 
$\tilde{v}_{0}\left( x\right) =\left( 1-x\right) \left( e^{9}+1\right) ^{%
	\frac{-1}{4}},$ $\ \tilde{v}_{1}\left( x\right) =\left( x-1\right)
e^{9}\left( e^{9}+1\right) ^{\frac{-5}{4}}.$%
\end{tabular}%
\end{equation*}
Then our external functions are%
\begin{align*}
&F_{1}\left( x,t\right) =-\frac{4x^{3}-2x^{2}+x}{\left( e^{9+4t}+1\right)
	^{3/4}}-\frac{5e^{9+4t}x}{\left( e^{9+4t}+1\right) ^{9/4}},
\\& F_{2}\left(
x,t\right) =\frac{\left( x-1\right) \left( 4x^{2}-6x+3\right) }{\left(
	e^{9+4t}+1\right) ^{3/4}}-\frac{5e^{9+4t}\left( 1-x\right) }{\left(
	e^{9+4t}+1\right) ^{9/4}}.
\end{align*}

This way we can find the exact solutions in the following form:
\begin{equation}
u_{ex}\left( x,t\right) =\frac{x}{\sqrt[4]{e^{9+4t}+1}},\; v_{ex}\left(
x,t\right) =\frac{1-x}{\sqrt[4]{e^{9+4t}+1}},
\tag{6.1}  \label{f1}
\end{equation}%
but they shall be neglected since we want to consider the illustrative
approximation of solutions.

\begin{figure}[htp]
	\centering
	\subcaptionbox{Exact solution $u_{ex}(x,t)$}{\includegraphics[width=0.50\textwidth]{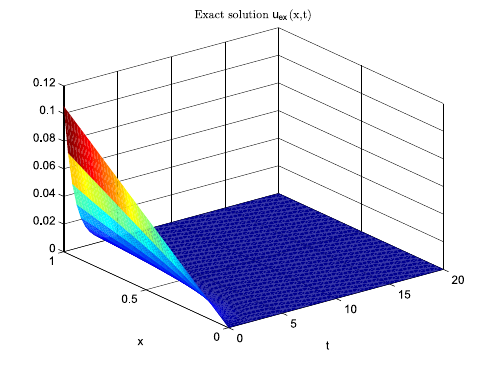}}%
	\subcaptionbox{Exact solution $v_{ex}(x,t)$}{\includegraphics[width=0.50\textwidth]{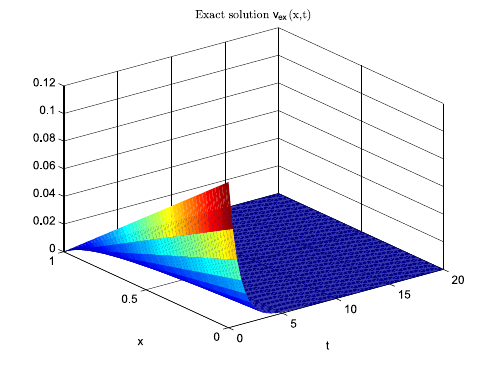}}%
	\caption{Exact solutions.}\label{Fig:1}
\end{figure}

In order to fulfill all assumptions of Theorem \ref{thm:5.2}, it remains to check (\ref%
{e10}) and (\ref{e11}). Obviously, one has $p_{\ast }=3$ and $\rho $ is
computed by the following%
\begin{align}\tag{6.2}  \label{f2}
\rho &=\frac{1}{2}\int_{0}^{\infty }\left( \left\Vert F_{1}\left(
s\right) \right\Vert +\left\Vert F_{2}\left( s\right) \right\Vert \right)
ds\leq \frac{1}{2}\int_{0}^{\infty }\left( \frac{12}{\left(
	e^{9+4s}+1\right) ^{3/4}}+\frac{18}{\left( e^{9+4s}+1\right) ^{3/4}}\right)
ds \\ 
& =15\int_{0}^{\infty }\frac{1}{\left( e^{9+4s}+1\right) ^{3/4}}%
ds\leq 15\int_{0}^{\infty }\frac{1}{e^{\frac{27}{4}+3s}}ds=15e^{%
	\frac{-27}{4}}\int_{0}^{\infty }e^{-3s}ds \nonumber \\ 
& =5e^{\frac{-27}{4}}\approx 0.58\times 10^{-4}\approx 0,\nonumber 
\end{align}%
and the initial energy can be estimated by
\begin{align}\tag{6.3}
\label{f3}
E\left( 0\right) &=\frac{1}{2}\left( \left\Vert \tilde{u}_{1}\right\Vert
^{2}+\left\Vert \tilde{v}_{1}\right\Vert ^{2}+\left\Vert \tilde{u}%
_{0x}\right\Vert ^{2}+\left\Vert \tilde{v}_{0x}\right\Vert ^{2}\right)
-\left( \left\vert \tilde{u}_{0}\left( 1\right) \right\vert ^{6}+\left\vert 
\tilde{v}_{0}\left( 0\right) \right\vert ^{6}\right)   \\
&-\frac{3}{4}\int_{0}^{1}\left( \left\vert \tilde{u}_{0}\left( x\right)
\right\vert ^{4}+\left\vert \tilde{v}_{0}\left( x\right) \right\vert
^{4}\right) dx-\frac{1}{2}\int_{0}^{1}\left\vert \tilde{u}_{0}\left(
x\right) \right\vert ^{2}\left\vert \tilde{v}_{0}\left( x\right) \right\vert
^{2}dx  \notag \\
&=\frac{1}{3}e^{18}\left( e^{9}+1\right) ^{-\frac{5}{2}}+\left(
e^{9}+1\right) ^{\frac{-1}{2}}-2\left( e^{9}+1\right) ^{-\frac{3}{2}}  \notag \\&
-\frac{3%
}{10}\left( e^{9}+1\right) ^{-1}-\frac{1}{60}\left( e^{9}+1\right)
^{-1}<0.015.\nonumber   
\end{align}

Combining (\ref{f2}) and (\ref{f3}), the checking argument for (\ref%
{e11}) is made by (\ref{f2}) and we obtain that%
\begin{equation*}
E_{\ast }=\left[ E\left( 0\right) +\rho \right] \mbox{exp}\left( 2\rho
\right) \approx E\left( 0\right) <0.015,
\end{equation*}%
which leads to $\eta _{\ast }=12p_{\ast }E_{\ast }+2\left( p_{\ast }E_{\ast
}\right) ^{2}<1.$

Therefore, the assumptions needed to check are all satisfied.

At the discretization level for this problem, a uniform grid of mesh-points $%
\left( x_{k},t_{n}\right) $ is used. Here $x_{k}=k\Delta x$ and $%
t_{n}=n\Delta t$ where $k$ and $n$ are integers and $\Delta x=\frac{1}{K},$ $%
\Delta t=\frac{T}{N}$ the equivalent mesh-widths in space $x$ and time $t$,
respectively. We first consider the following differential system for
the unknowns $\left( U_{k}\left( t\right) ,V_{k}\left( t\right) \right)
\equiv \left( u\left( x_{k},t\right) ,v\left( x_{k},t\right) \right) $ for $%
k=\overline{0,K}$:%
\begin{align}\tag{6.4}  \label{f4}
\frac{d\bar{U}_{1}}{dt}\left( t\right) &=K^{2}\left( U_{2}\left( t\right)
-2U_{1}\left( t\right) \right) -\bar{U}_{1}\left( t\right) \nonumber \\&+3U_{1}^{3}\left(
t\right) +V_{1}^{2}\left( t\right) U_{1}\left( t\right) +F_{1}\left(
x_{1},t\right) ,   \nonumber\\ 
\frac{d\bar{U}_{k}}{dt}\left( t\right) &=K^{2}\left( U_{k-1}\left( t\right)
-2U_{k}\left( t\right) +U_{k+1}\left( t\right) \right) -\bar{U}_{k}\left(
t\right) \nonumber  \\ 
&+3U_{k}^{3}\left( t\right) +V_{k}^{2}\left( t\right) U_{k}\left(
t\right) +F_{1}\left( x_{k},t\right) ,\;k=\overline{2,K-1}, \nonumber\\ 
\frac{d\bar{U}_{K}}{dt}\left( t\right) &=K\left( U_{K}^{5}\left( t\right) -%
\bar{U}_{K}\left( t\right) \right) -K^{2}\left( U_{K}\left( t\right)
-U_{K-1}\left( t\right) \right) \nonumber \\ &-\bar{U}_{K}\left( t\right)
+3U_{K}^{3}\left( t\right) +F_{1}\left( x_{K},t\right) ,
\nonumber 
\end{align}%
and
\begin{align}\tag{6.5}  \label{f5}
\frac{d\bar{V}_{0}}{dt}\left( t\right) &=K\left( V_{0}^{5}\left( t\right) -%
\bar{V}_{0}\left( t\right) \right) +K^{2}\left( V_{1}\left( t\right)
 -V_{0}\left( t\right) \right) \\ &-\bar{V}_{0}\left( t\right) +3V_{0}^{3}\left(
t\right) +F_{2}\left( x_{0},t\right) , \nonumber\\ 
\frac{d\bar{V}_{k}}{dt}\left( t\right) &=K^{2}\left( V_{k-1}\left( t\right)
-2V_{k}\left( t\right) +V_{k+1}\left( t\right) \right) -\bar{V}_{k}\left(
t\right)  \nonumber\\ 
&+3V_{k}^{3}\left( t\right) +U_{k}^{2}\left( t\right) V_{k}\left(
t\right) +F_{2}\left( x_{k},t\right) ,\;k=\overline{1,K-2}, \nonumber \\ 
\frac{d\bar{V}_{K-1}}{dt}\left( t\right) &=K^{2}\left( V_{K-2}\left(
t\right) -2V_{K-1}\left( t\right) \right) -\bar{V}_{K-1}\left( t\right)
\nonumber \\ &+3V_{K-1}^{3}\left( t\right) +V_{K-1}^{2}\left( t\right) U_{K-1}\left(
t\right) +F_{2}\left( x_{K-1},t\right) ,\nonumber
\end{align}%
where $\left( \bar{U}_{k}\left( t\right) ,\bar{V}_{k}\left( t\right) \right) 
$ identically stands for $\left( \frac{dU_{k}}{dt}\left( t\right) ,\frac{%
	dV_{k}}{dt}\left( t\right) \right) $ and the initial conditions are%
\begin{equation}\tag{6.6}  \label{f6}
\left( U_{k}\left( 0\right) ,V_{k}\left( 0\right) \right) =\left( \tilde{u}%
_{0}\left( x_{k}\right) ,\tilde{v}_{0}\left( x_{k}\right) \right) , 
\left( \bar{U}_{k}\left( 0\right) ,\bar{V}_{k}\left( 0\right) \right)
=\left( \tilde{u}_{1}\left( x_{k}\right) ,\tilde{v}_{1}\left( x_{k}\right)
\right) ,\;k=\overline{0,K}.
\end{equation}
Here we also notice that the values of $U_{0}\left(t\right)$ and $%
V_{K}\left(t\right)$ are known, so those cases are not considered.

Below a basic numerical approach is achieved by using the linear recursive
method where the nonlinear terms are linearized. Accordingly, after some
rearrangements we rewrite the linearized differential system to be a
differential equation where the solution includes all discrete solutions of
the linearized system. Such a unifying way also guarantees that the
numerical solution uniquely exists. When doing this, at the $m$-th
iterative stage ($m\geq 1$) the linearized differential system of (\ref{f4}%
)-(\ref{f6}) one by one becomes
\begin{align}\tag{6.7}  \label{f7}
\frac{d\bar{U}_{1}^{\left( m\right) }}{dt}\left( t\right) &=K^{2}\left(
U_{2}^{\left( m\right) }\left( t\right) -2U_{1}^{\left( m\right) }\left(
t\right) \right) -\bar{U}_{1}^{\left( m\right) }\left( t\right)  \\ 
& +3\left( U_{1}^{\left( m-1\right) }\left(
t\right) \right) ^{3}+\left( V_{1}^{\left( m-1\right) }\left( t\right)
\right) ^{2}U_{1}^{\left( m\right) }\left( t\right) +F_{1}\left(
x_{1},t\right) ,\nonumber  \\ 
\frac{d\bar{U}_{k}^{\left( m\right) }}{dt}\left( t\right) &=K^{2}\left(
U_{k-1}^{\left( m\right) }\left( t\right) -2U_{k}^{\left( m\right) }\left(
t\right) +U_{k+1}^{\left( m\right) }\left( t\right) \right) -\bar{U}%
_{k}^{\left( m\right) }\left( t\right) \nonumber \\ 
&+3\left( U_{k}^{\left( m-1\right) }\left(
t\right) \right) ^{3}+\left( V_{k}^{\left( m-1\right) }\left( t\right)
\right) ^{2}U_{k}^{\left( m\right) }\left( t\right) +F_{1}\left(
x_{k},t\right) ,\;k=\overline{2,K-1}, \nonumber  \\ 
\frac{d\bar{U}_{K}^{\left( m\right) }}{dt}\left( t\right) &=K\left[ \left(
U_{K}^{\left( m-1\right) }\left( t\right) \right) ^{5}-\bar{U}_{K}^{\left(
	m\right) }\left( t\right) \right] -K^{2}\left( U_{K}^{\left( m\right)
}\left( t\right) -U_{K-1}^{\left( m\right) }\left( t\right) \right) \nonumber \\ 
& -\bar{U}_{K}^{\left( m\right) }\left( t\right)
+3\left( U_{K}^{\left( m-1\right) }\left( t\right) \right) ^{3}+F_{1}\left(
x_{K},t\right) ,
\nonumber 
\end{align}%
and
\begin{align}\tag{6.8}  \label{f8}
\frac{d\bar{V}_{0}^{\left( m\right) }}{dt}\left( t\right) &=K\left[ \left(
V_{0}^{\left( m-1\right) }\left( t\right) \right) ^{5}-\bar{V}_{0}^{\left(
	m\right) }\left( t\right) \right] +K^{2}\left( V_{1}^{\left( m\right)
}\left( t\right) -V_{0}^{\left( m\right) }\left( t\right) \right)  \\ 
& -\bar{V}_{0}^{\left( m\right) }\left( t\right)
+3\left( V_{0}^{\left( m-1\right) }\left( t\right) \right) ^{3}+F_{2}\left(
x_{0},t\right) , \nonumber\\ 
\frac{d\bar{V}_{k}^{\left( m\right) }}{dt}\left( t\right) &=K^{2}\left(
V_{k-1}^{\left( m\right) }\left( t\right) -2V_{k}^{\left( m\right) }\left(
t\right) +V_{k+1}^{\left( m\right) }\left( t\right) \right) -\bar{V}%
_{k}^{\left( m\right) }\left( t\right) \nonumber \\ 
& +3\left( V_{k}^{\left( m-1\right) }\left(
t\right) \right) ^{3}+\left( U_{k}^{\left( m-1\right) }\left( t\right)
\right) ^{2}V_{k}^{\left( m\right) }\left( t\right) +F_{2}\left(
x_{k},t\right) ,\;k=\overline{1,K-2}, \nonumber \\ 
\frac{d\bar{V}_{K-1}^{\left( m\right) }}{dt}\left( t\right) &=K^{2}\left(
V_{K-2}^{\left( m\right) }\left( t\right) -2V_{K-1}^{\left( m\right) }\left(
t\right) \right) -\bar{V}_{K-1}^{\left( m\right) }\left( t\right) \nonumber \\ 
& +3\left( V_{K-1}^{\left( m-1\right) }\left(
t\right) \right) ^{3}+\left( U_{K-1}^{\left( m-1\right) }\left( t\right)
\right) ^{2}V_{K-1}^{\left( m\right) }\left( t\right) +F_{2}\left(
x_{K-1},t\right) ,\nonumber
\end{align}%
\begin{align}\tag{6.9}  \label{f9}
&\left( U_{k}^{\left( m\right) }\left( 0\right) ,V_{k}^{\left( m\right)
}\left( 0\right) \right) =\left( \tilde{u}_{0}\left( x_{k}\right) ,
\tilde{v}%
_{0}\left( x_{k}\right) \right) , \\ &\left( \bar{U}_{k}^{\left( m\right)
}\left( 0\right) ,\bar{V}_{k}^{\left( m\right) }\left( 0\right) \right)
=\left( \tilde{u}_{1}\left( x_{k}\right) ,\tilde{v}_{1}\left( x_{k}\right)
\right) ,k=\overline{0,K}.
\nonumber
\end{align}
Then, a matrix-type differential system derived from (\ref{f7})-(\ref{f9})
can be expressed as follows:
\begin{equation}
\frac{d\mathbb{U}^{\left( m\right) }}{dt}\left( t\right) =\mathbb{A}%
^{\left( m\right) }\left( t\right) \mathbb{U}^{\left( m\right) }\left(
t\right) +\mathbb{F}_{1}^{\left( m\right) }\left( t\right) ,\; \frac{d%
	\mathbb{V}^{\left( m\right) }}{dt}\left( t\right) =\mathbb{B}^{\left(
	m\right) }\left( t\right) \mathbb{V}^{\left( m\right) }\left( t\right) +%
\mathbb{F}_{2}^{\left( m\right) }\left( t\right) ,
\tag{6.10}  \label{f10}
\end{equation}%
where the solutions $\mathbb{U}^{\left( m\right) },$ $\mathbb{V}^{\left(
	m\right) }\in \mathbb{R}^{2K},$ and the block matrices $\mathbb{A}^{\left(
	m\right) },$ $\mathbb{B}^{\left( m\right) }\in \mathbb{R}^{2K}\times \mathbb{%
	R}^{2K},$ and the functions $\mathbb{F}_{1}^{\left( m\right) },$ $\mathbb{F}%
_{2}^{\left( m\right) }\in \mathbb{R}^{2K}$ are defined by%
\begin{equation*}
\mathbb{U}^{\left( m\right) }=%
\begin{bmatrix}
U_{1}^{\left( m\right) } \\ 
\vdots \\ 
U_{K}^{\left( m\right) } \\ 
\bar{U}_{1}^{\left( m\right) } \\ 
\vdots \\ 
\bar{U}_{K-1}^{\left( m\right) } \\ 
\bar{U}_{K}^{\left( m\right) }%
\end{bmatrix}%
,\text{ \ }\mathbb{F}_{1}^{\left( m\right) }\left( t\right) =%
\begin{bmatrix}
0 \\ 
\vdots \\ 
0 \\ 
\alpha _{1}^{\left( m\right) }\left( t\right) \\ 
\vdots \\ 
\alpha _{K-1}^{\left( m\right) }\left( t\right) \\ 
\alpha _{K}^{\left( m\right) }\left( t\right) +K\left( U_{K}^{\left(
	m-1\right) }\left( t\right) \right) ^{5}%
\end{bmatrix}%
,
\end{equation*}%
\begin{equation*}
\mathbb{V}^{\left( m\right) }=%
\begin{bmatrix}
V_{0}^{\left( m\right) } \\ 
V_{1}^{\left( m\right) } \\ 
\vdots \\ 
V_{K-1}^{\left( m\right) } \\ 
\bar{V}_{0}^{\left( m\right) } \\ 
\bar{V}_{1}^{\left( m\right) } \\ 
\vdots \\ 
\bar{V}_{K-1}^{\left( m\right) }%
\end{bmatrix}%
,\text{ \ }\mathbb{F}_{2}^{\left( m\right) }\left( t\right) =%
\begin{bmatrix}
0 \\ 
0 \\ 
\vdots \\ 
0 \\ 
\beta _{0}^{\left( m\right) }\left( t\right) +K\left( V_{0}^{\left(
	m-1\right) }\left( t\right) \right) ^{5} \\ 
\beta _{1}^{\left( m\right) }\left( t\right) \\ 
\vdots \\ 
\beta _{K-1}^{\left( m\right) }\left( t\right)%
\end{bmatrix}%
,
\end{equation*}%
\begin{equation*}
\begin{tabular}{l}
$\alpha _{k}^{\left( m\right) }\left( t\right) =3\left( U_{k}^{\left(
	m-1\right) }\left( t\right) \right) ^{3}+F_{1}\left( x_{k},t\right) ,\text{ }%
k=\overline{1,K},$%
\end{tabular}%
\end{equation*}%
\begin{equation*}
\begin{tabular}{l}
$\beta _{k}^{\left( m\right) }\left( t\right) =3\left( V_{k}^{\left(
	m-1\right) }\left( t\right) \right) ^{3}+F_{2}\left( x_{k},t\right) ,\text{ }%
k=\overline{0,K-1},$%
\end{tabular}%
\end{equation*}%
\begin{equation*}
\mathbb{A}^{\left( m\right) }=%
\begin{bmatrix}
{\LARGE O} & 
\begin{array}{ccccc}
1 & 0 & \cdots & \cdots & 0 \\ 
0 & 1 &  &  & \vdots \\ 
\vdots &  & \ddots &  & \vdots \\ 
\vdots &  &  & 1 & 0 \\ 
0 & \cdots \text{ \ } & \cdots \text{ \ } & 0 & 1%
\end{array}%
\text{ \ } \\ 
\begin{array}{ccccc}
a_{1}^{\left( m\right) } & K^{2} & \cdots & \cdots & 0 \\ 
K^{2} & a_{2}^{\left( m\right) } & K^{2} &  & \vdots \\ 
\vdots & \ddots & \ddots & \ddots & \vdots \\ 
\vdots &  & K^{2} & a_{K-1}^{\left( m\right) } & K^{2} \\ 
0 & \cdots & \cdots & K^{2} & -K^{2}%
\end{array}
& 
\begin{array}{ccccc}
-1 & 0 & \cdots & \cdots & 0 \\ 
0 & -1 &  &  & \vdots \\ 
\vdots &  & \ddots &  & \vdots \\ 
\vdots &  &  & -1 & 0 \\ 
0 & \cdots & \cdots & 0 & -1-K%
\end{array}%
\end{bmatrix}%
,
\end{equation*}%
\begin{equation*}
\mathbb{B}^{\left( m\right) }=%
\begin{bmatrix}
{\LARGE O} & \text{ \ }%
\begin{array}{ccccc}
1 & 0 & \cdots & \cdots & 0 \\ 
0 & 1 &  &  & \vdots \\ 
\vdots &  & \ddots &  & \vdots \\ 
\vdots &  &  & 1 & 0 \\ 
0 & \cdots \text{ \ \ } & \cdots \text{ \ \ } & 0 & 1%
\end{array}
\\ 
\begin{array}{ccccc}
-K^{2} & K^{2} & \cdots & \cdots & 0 \\ 
K^{2} & b_{1}^{\left( m\right) } & K^{2} &  & \vdots \\ 
\vdots & \ddots & \ddots & \ddots & \vdots \\ 
\vdots &  & K^{2} & b_{K-2}^{\left( m\right) } & K^{2} \\ 
0 & \cdots & \cdots & K^{2} & b_{K-1}^{\left( m\right) }%
\end{array}
& 
\begin{array}{ccccc}
-1-K & 0 & \cdots & \cdots & 0 \\ 
0 & -1 &  &  & \vdots \\ 
\vdots &  & \ddots &  & \vdots \\ 
\vdots &  &  & -1 & 0 \\ 
0 & \cdots & \cdots & 0 & -1%
\end{array}%
\end{bmatrix}%
,
\end{equation*}%
\begin{align*}
&a_{k}^{\left( m\right) }:=a_{k}^{\left( m\right) }\left( t\right)
=-2K^{2}+\left( V_{k}^{\left( m-1\right) }\left( t\right) \right)
^{2},\\&b_{k}^{\left( m\right) }:=b_{k}^{\left( m\right) }\left( t\right)
=-2K^{2}+\left( U_{k}^{\left( m-1\right) }\left( t\right) \right) ^{2},\;k=%
\overline{1,K-1},
\end{align*}%
associated with the initial data%
\begin{equation}
\mathbb{U}^{\left( m\right) }\left( 0\right) =%
\begin{bmatrix}
\tilde{u}_{0}\left( x_{1}\right) \\ 
\vdots \\ 
\tilde{u}_{0}\left( x_{K}\right) \\ 
\tilde{u}_{1}\left( x_{1}\right) \\ 
\vdots \\ 
\tilde{u}_{1}\left( x_{K}\right)%
\end{bmatrix}%
,\text{ \ \ \ }\mathbb{V}^{\left( m\right) }\left( 0\right) =%
\begin{bmatrix}
\tilde{v}_{0}\left( x_{0}\right) \\ 
\vdots \\ 
\tilde{v}_{0}\left( x_{K-1}\right) \\ 
\tilde{v}_{1}\left( x_{0}\right) \\ 
\vdots \\ 
\tilde{v}_{1}\left( x_{K-1}\right)%
\end{bmatrix}%
.  \tag{6.11}  \label{f11}
\end{equation}

Here the initial guess for our linearization method is simply the initial
values, i.e. $\left( \mathbb{U}^{\left( 0\right) }\left( t_{n}\right) ,%
\mathbb{V}^{\left( 0\right) }\left( t_{n}\right) \right) =$ $\left( \mathbb{U%
}^{\left( m\right) }\left( 0\right) ,\mathbb{V}^{\left( m\right) }\left(
0\right) \right) $ for all $n=\overline{0,N}$. In general, it is remarkable
that the system (\ref{f10})-(\ref{f11}) might be stiff, then it is not good
for our implementation since the objective is considering the system in a
large time. Understanding the basic instability coming from stiff systems,
we therefore apply the well-known implicit Euler method which reads as
\begin{align*}
& \mathbb{U}_{n+1}^{\left( m\right) }=\left( \mathbb{I}-\frac{T}{N}\mathbb{A}%
^{\left( m\right) }\left( t_{n+1}\right) \right) ^{-1}\left( \mathbb{U}%
_{n}^{\left( m\right) }+\frac{T}{N}\mathbb{F}_{1}^{\left( m\right) }\left(
t_{n+1}\right) \right) , \\ &
\mathbb{V}_{n+1}^{\left( m\right) }=\left( \mathbb{I}-\frac{T}{N}\mathbb{B}%
^{\left( m\right) }\left( t_{n+1}\right) \right) ^{-1}\left( \mathbb{V}%
_{n}^{\left( m\right) }+\frac{T}{N}\mathbb{F}_{2}^{\left( m\right) }\left(
t_{n+1}\right) \right) ,
\end{align*}%
for $n=\overline{0,N-1}$, where $\mathbb{I}$ stands for the $2K$-by-$2K$
identity matrix. This discrete system is endowed with the conditions (\ref{f11}).

\begin{figure}[htp]
	\centering
	\subcaptionbox{Approximate solution $u(x,t)$}{\includegraphics[width=0.50\textwidth]{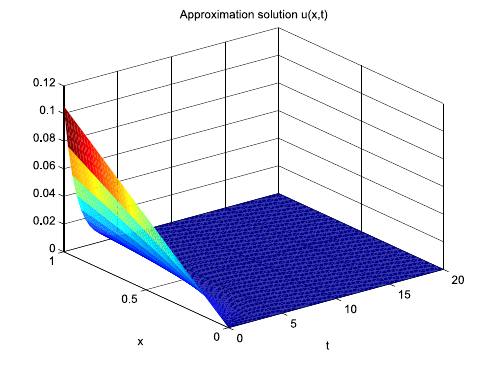}}%
	\subcaptionbox{Approximate solution $v(x,t)$}{\includegraphics[width=0.50\textwidth]{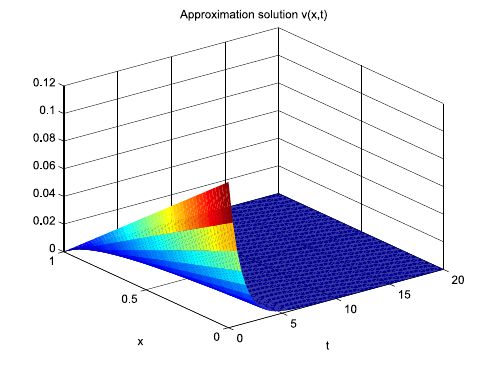}}%
	\caption{Approximate solutions.}\label{Fig:2}
\end{figure}

Choosing $m=5,$ $T=20$, and $K=N=50$, we have plotted in Figure \ref{Fig:2} the approximation
function $\left( u\left( x,t\right) ,v\left( x,t\right) \right) $ solutions
to our problem $\left( P\right) $ considered in this section. Compared to Figure \ref{Fig:1}, we observe that such functions not only behave like the exact solutions (\ref{f1}%
) (decay exponentially in time), but also completely have the same shapes
corresponding to each exact solution.

\begin{table}[h!]
	\caption{Numerical results at nodes $\left( \frac{4}{5}%
		,t_{n}\right)$ for $n\in\left\{10,20,30\right\}.$}\label{table:1}
	\begin{center}
		\begin{tabular}{|c|c|c|c|}
			\hline
			$n$ & $u_{ex}\left( \frac{4}{5},t_{n}\right) $ & $u\left( \frac{4}{5}%
			,t_{n}\right) $ & $\left\vert u_{ex}\left( \frac{4}{5},t_{n}\right) -u\left( 
			\frac{4}{5},t_{n}\right) \right\vert $ \\ \hline
			${\small 10}$ & \multicolumn{1}{|c|}{${\small 1.54436330E-03}$} & 
			\multicolumn{1}{|c|}{${\small 2.91855517E-03}$} & ${\small 1.37419186E-03}$
			\\ \hline
			${\small 20}$ & \multicolumn{1}{|c|}{${\small 2.82860006E-05}$} & 
			\multicolumn{1}{|c|}{${\small 7.20712002E-05}$} & ${\small 4.37851996E-05}$
			\\ \hline
			${\small 30}$ & \multicolumn{1}{|c|}{${\small 5.18076174E-07}$} & 
			\multicolumn{1}{|c|}{${\small 1.77972692E-06}$} & ${\small 1.26165074E-06}$
			\\ \hline
		\end{tabular}%
	\end{center}%
	\begin{center}
		\begin{tabular}{|c|c|c|c|}
			\hline
			$n$ & $v_{ex}\left( \frac{4}{5},t_{n}\right) $ & $v\left( \frac{4}{5}%
			,t_{n}\right) $ & $\left\vert v_{ex}\left( \frac{4}{5},t_{n}\right) -v\left( 
			\frac{4}{5},t_{n}\right) \right\vert $ \\ \hline
			${\small 10}$ & \multicolumn{1}{|c|}{${\small 3.86090827E-04}$} & 
			\multicolumn{1}{|c|}{${\small 7.29514168E-04}$} & ${\small 3.43423340E-04}$
			\\ \hline
			${\small 20}$ & \multicolumn{1}{|c|}{${\small 7.07150017E-06}$} & 
			\multicolumn{1}{|c|}{${\small 1.80147701E-05}$} & ${\small 1.09432699E-05}$
			\\ \hline
			${\small 30}$ & \multicolumn{1}{|c|}{${\small 1.29519043E-07}$} & 
			\multicolumn{1}{|c|}{${\small 6.22799676E-06}$} & ${\small 4.93280633E-07}$
			\\ \hline
		\end{tabular}%
	\end{center}%
\end{table}

\begin{table}[h!]
	\caption{Numerical results for the $l_{\infty } $ norm error $\mathcal{E}_{N,K}$.}\label{table:2}
	\begin{tabular}{|l|l|c|c|}
		\hline
		$K$ & $N$ & $\mathcal{E}_{N,K}\left( u\right) $ & $\mathcal{E}_{N,K}\left(
		v\right) $ \\ \hline
		${\small 50}$ & ${\small 50}$ & ${\small 6.68545424E-03}$ & ${\small %
			6.68150701E-03}$ \\ \hline
		${\small 100}$ & ${\small 100}$ & ${\small 3.59475057E-03}$ & ${\small %
			3.59201931E-03}$ \\ \hline
		${\small 150}$ & ${\small 150}$ & ${\small 2.45841870E-03}$ & ${\small %
			2.45632948E-03}$ \\ \hline
		${\small 200}$ & ${\small 200}$ & ${\small 1.86793338E-03}$ & ${\small %
			1.86628504E-03}$ \\ \hline
	\end{tabular}
\end{table}

Furthermore, numerical results of the solutions $\left( u,v\right) $
together with the exact solutions $\left( u_{ex},v_{ex}\right) $ at nodes $%
\left( \frac{4}{5},t_{n}\right) $ for $n\in \left\{10,20,30\right\}$ and various values
of error in the entry-wise norm%
\begin{align*}
&\mathcal{E}_{N,K}\left( u\right) =\max_{1\leq k\leq K}\max_{1\leq n\leq
	N}\left\vert u_{ex}\left( x_{k},t_{n}\right) -u\left( x_{k},t_{n}\right)
\right\vert , \\ 
&\mathcal{E}_{N,K}\left( v\right) =\max_{1\leq k\leq K}\max_{0\leq n\leq
	N-1}\left\vert v_{ex}\left( x_{k},t_{n}\right) -v\left( x_{k},t_{n}\right)
\right\vert ,
\end{align*}%
are all given in Table \ref{table:1}-Table \ref{table:2}, respectively. To demonstrate the fact that
the $l_{\infty }$ norm error $\mathcal{E}_{N,K}$ decreases (and obviously
tends to zero) as $K,N$ increase without any attention to the
discretization level, we show in Table \ref{table:2} the error values
when such constants go from $50$ to $200$. As expected from the analysis, our
numerical method used here is reasonable and efficient.

\begin{appendices}\section{Proofs of auxiliary results}\label{app1}

	\begin{proof}[Proof of Lemma \ref{lem:4.2}]
		We only need to prove the first inequality since (\ref{d15})
		and (\ref{d16}) are obviously the same. Put $s_{1}=2/\left( 1-2\xi \right) $
		and $s_{2}=2/\left( 1-\xi \right) $. Our strategy is to independently
		consider each terms on the left-hand side of (\ref{d15}) and furthermore, to
		estimate those quantities we mainly divide into two cases. For clarity, let
		us first consider two cases for $\left\Vert u\right\Vert _{L^{\alpha }}$.
		
		$(a_{1})$ $\ \left\Vert u\right\Vert _{L^{\alpha }}\leq 1$: One easily
		deduces from $2\leq s_{1}\leq \alpha $ that%
		\begin{equation*}
		\begin{tabular}{l}
		$\left\Vert u \right\Vert _{L^{\alpha }}^{s_{1}}\leq
		\left\Vert u \right\Vert _{L^{\alpha }}^{2}\leq \left\Vert
		u_{x} \right\Vert ^{2}+\left\Vert u
		\right\Vert _{L^{\alpha }}^{\alpha }+\left\vert u\left( 1\right)
		\right\vert ^{p_{1}}.$%
		\end{tabular}%
		\end{equation*}
		
		$(a_{2})$ $\ \left\Vert u\right\Vert _{L^{\alpha }}\geq 1$: Similarly, we get%
		\begin{equation*}
		\begin{tabular}{l}
		$\left\Vert u \right\Vert _{L^{\alpha }}^{s_{1}}\leq
		\left\Vert u \right\Vert _{L^{\alpha }}^{\alpha }\leq
		\left\Vert u_{x} \right\Vert ^{2}+\left\Vert u
		\right\Vert _{L^{\alpha }}^{\alpha }+\left\vert u\left( 1\right)
		\right\vert ^{p_{1}}.$%
		\end{tabular}%
		\end{equation*}
		
		Next, we consider two cases for $\left\Vert u\right\Vert $:
		
		$(b_{1})$ $\ \left\Vert u\right\Vert \leq 1$: Since $2\leq s_{2}\leq \alpha $%
		,%
		\begin{equation*}
		\begin{tabular}{l}
		$\left\Vert u \right\Vert ^{s_{2}}\leq \left\Vert u \right\Vert ^{2}\leq \left\Vert u_{x} \right\Vert
		^{2}+\left\Vert u \right\Vert _{L^{\alpha }}^{\alpha
		}+\left\vert u\left( 1\right) \right\vert ^{p_{1}}.$%
		\end{tabular}%
		\end{equation*}
		
		$(b_{2})$ $\ \left\Vert u\right\Vert \geq 1$: In the same manner, we
		have
		\begin{equation*}
		\begin{tabular}{l}
		$\left\Vert u \right\Vert ^{s_{2}}\leq \left\Vert u \right\Vert ^{\alpha }\leq \left\Vert u \right\Vert
		_{L^{\alpha }}^{\alpha }\leq \left\Vert u_{x} \right\Vert
		^{2}+\left\Vert u \right\Vert _{L^{\alpha }}^{\alpha
		}+\left\vert u\left( 1\right) \right\vert ^{p_{1}}.$%
		\end{tabular}%
		\end{equation*}
		
		Finally, we go through the last term $\left\vert u\left( 1\right)
		\right\vert $:
		
		$(c_{1})$ $\ \left\vert u\left( 1\right) \right\vert \leq 1$:%
		\begin{equation*}
		\begin{tabular}{l}
		$\left\vert u\left( 1\right) \right\vert ^{s_{2}}\leq \left\vert u\left(
		1\right) \right\vert ^{2}\leq \left\Vert u \right\Vert
		_{C\left( \overline{\Omega }\right) }^{2}\leq \left\Vert u_{x} \right\Vert
		^{2}+\left\Vert u \right\Vert _{L^{\alpha }}^{\alpha
		}+\left\vert u\left( 1\right) \right\vert ^{p_{1}}.$%
		\end{tabular}%
		\end{equation*}
		
		$(c_{2})$ $\ \left\vert u\left( 1\right) \right\vert \geq 1$: It is
		natural to say $\left\vert u\left( 1\right) \right\vert ^{s_{2}}\leq
		\left\vert u\left( 1\right) \right\vert ^{p_{1}}$.
		
		Combining all of the above inequalities completes the proof of the lemma.
	\end{proof}
\end{appendices}

\section*{Acknowledgments} The authors would like to thank the anonymous referees for fruitful comments through the improvement of this paper. V.A.K thanks Prof. Pierangelo Marcati (L'Aquila, Italy) for the trust and mathematical training that the Gran Sasso Science Institute has invested in V.A.K through the PhD time. V.A.K gratefully acknowledges Dr. Nguyen Thai Ngoc (G\"ottingen, Germany) for his unconditional helps during the time V.A.K shortly occupied a research associate position at the institute.


\medskip
Received xxxx 20xx; revised xxxx 20xx.
\medskip

\end{document}